\documentstyle[leqno]{article}

\newcommand{\carre}{$\rule{6pt}{6pt}$}

\newcommand{\bl}{\left( \begin{array}}
\newcommand{\er}{\end{array} \right)}
\newcommand{\n}{\noindent}
\newcommand{\spd}{\vspace{0.3cm}}
\newcommand{\spf}{\vspace{0.5cm}}

\newcommand{\hsps}{\hspace{0.7cm}}

\newcommand{\bbc}{\hbox{\rm\rlap {I}\kern-.125em{C}}}
\newcommand{\bbr}{\hbox{\rm\rlap {I}\thinspace {R}}}
\newcommand{\bbz}{\hbox{\rm\rlap {Z}\thinspace {Z}}}

\newtheorem{lemma}{Lemma}[section]

\newtheorem{prop}[lemma]{Proposition}
\newtheorem{definition}[lemma]{Definition}
\newtheorem{exam}[lemma]{Examples}
\newtheorem{theor}{Theorem}[section]
\newtheorem{rema}[lemma]{Remark}

\begin{document}

\thispagestyle{empty}

\title{Relative torsion}

\author{D. Burghelea (Ohio State University)\\ 
L. Friedlander (University of Arizona)\\ 
T. Kappeler (University of Z\"urich)}

\maketitle

\n {\bf Abstract:}  This paper achieves, among other things, the following:

\begin{itemize}
\item It frees the main result of [BFKM] from the hypothesis
of determinant class and extends this result from unitary 
to arbitrary representations.

\item It extends (and at the same times provides a new proof
of) the main result of Bismut and Zhang [BZ] from finite dimensional representations of $\Gamma$
to representations on an ${\cal A}-$Hilbert module of finite type ( ${\cal A}$ a finite von Neumann
algebra). The result of [BZ] corresponds to  ${\cal A}=\bbc.$

\item It provides interesting real valued functions on the 
space of representations 
of the fundamental group $\Gamma$ of a closed manifold $M.$ 
These functions  might be a useful source of topological and 
geometric invariants of $M$.

\end{itemize}

These objectives are achieved with the help of the relative torsion $\cal R $, first  introduced 
by Carey, Mathai
and Mishchenko [CMM]  in  special cases. The main result of this paper 
calculates explicitly this relative torsion
(cf Theorem 0.1).

\newpage

\n {\bf Table of contents}

\begin{itemize}
\item[0.] Introduction
\begin{itemize}
\item[0.1] Summary of the results
\item[0.2] Definitions and outline of the paper
\end{itemize}
\item[1.] Operators and complexes
\begin{itemize}
\item[1.1] Operators
\item[1.2] Complexes
\end{itemize}
\item[2.] Relative torsion and its Witten deformation
\begin{itemize}
\item[2.1] Relative torsion 
\item[2.2] Witten deformation of the relative torsion
\end{itemize}
\item[3.] Proof of Proposition 0.1 in the case $g=g'$
\begin{itemize}
\item[3.1] Outline of the proof 
\item[3.2] Asymptotic expansion of $\log {\cal R}_{sm} (t)$ 
\item[3.3] Asymptotic expansion of $\log T_{la} (t)$
\end{itemize}
\item[4.] Anomalies of the relative torsion
\begin{itemize}
\item[4.1] Metric and Hermitian anomalies
\item[4.2] Proof of Proposition 0.1
\item[4.3] Anomaly with respect to triangulation
\end{itemize}
\item[5.] Proof of Theorem 0.1
\begin{itemize}
\item[5.1] Proof of Theorem 0.1
\item[5.2] 
An invariant for odd dimensional manifolds.
\end{itemize}

\end{itemize}

Appendix A \  \ Lemma of Carey-Mathai-Mishchenko

\vspace{0.3cm}

Appendix B \  \ Determinant class property

\spf

\n {\Large{\bf 0. Introduction}}

\vspace{0.5cm}

\n {\large{\bf 0.1 Summary of the results}}

\vspace{0.3cm}

Let $(M, x_0)$ be a closed, connected manifold with base point $x_0$ and $\rho : \Gamma \rightarrow 
Gl_{\cal A} ({\cal W})$ a representation of $\Gamma := \pi_1 (M; x_0)$ on an 
${\cal A}$-Hilbert module ${\cal W}$ of finite type 
with ${\cal A}$ a finite von Neumann algebra. 
Denote by ${\cal E} \rightarrow M$ the vector bundle, associated to 
$M$ and $\rho$, equipped with the canonical flat connection $\nabla.$
Let $\mu$ be a Hermitian structure of ${\cal E}$ 
and $\tau = (h,g')$  a (generalized) triangulation of $M$
(cf subsection 0.2).

\n The relative torsion is a numerical invariant 
associated to ${\cal F}=(M,\rho ,\mu ,g,\tau ).$  Unlike the analytic torsion $T_{an}$
associated to $(M, \rho ,\mu ,g),$ or the Reidemeister torsion $T_{Reid}$
associated to ${\cal F}$  (cf subsection 0.2),
both of which are defined only when  $(M, \rho)$ is  of determinant
class (cf Appendix B) , ${\cal R}$ is always
defined. If $(M,\rho )$ is of determinant class the relative torsion 
turns out to be the quotient of the analytic and Reidemeister torsion (cf (2.15)).
However, as shown in Appendix B, there are many pairs $(M,\rho)$ which are 
not of determinant class.

\vspace{0.1cm}

\begin{prop} Assume $\mu_0$ is a Hermitian structure on the bundle induced by $\rho$ 
which is parallel in the neighborhood of $Cr(h).$

\n (i) $\log {\cal R}= 
\int_{M\setminus Cr(h)} \alpha_{\cal F}$ where $ \alpha_{\cal F}$ is a density on 
$M\setminus Cr(h),$
which vanishes on a neighborhood of $Cr(h)$ and
is a local quantity (cf subsection 0.2), depending on $\rho, \mu, g,$ and $X:= - grad_{g'} h.$ 

\n (ii) If $\mu$ is parallel then ${\cal R} = 1.$
\end {prop}

\n Following the work of Bismut Zhang [BZ], we  consider the closed one-form 
$\theta(\rho, \mu) \in  \Omega^1 (M),$  the form 
$\Psi(TM,g) \in \Omega^{n-1}(TM\setminus M; {\cal O}_{TM})$ 
(${\cal O}_{TM}$ the orientation bundle of $TM$),
and for  two Hermitian structures $\mu_1$ and $\mu_2$
on the bundle ${\cal E}\to M$  induced from $\rho$, the smooth function 
$V(\rho, \mu_1,\mu_2) \in \Omega^0 (M)$ (cf section 4). When $\mu$ is parallel with respect to $\nabla$ on an open set $U,$ $U\subset M,$ then 
$\theta(\rho, \mu)=0$ on $U.$  Denote by $e(M,g)$ 
the Euler form associated with the Riemannian metric $g.$
The triangulation $\tau =(h,g')$ provides  the vector field $X = -\rm {grad}_{g'} h$
which will be  regarded  as a smooth map $X: M\setminus Cr(h) \to TM\setminus M.$
The above quantities are used in formulae expressing  the change of the relative torsion 
when one varies the Riemannian metric $g$, the Hermitian 
structure $\mu$, or the generalized triangulation $\tau $ (cf Propositions 4.1-4.3).
These formulae are used to derive the following

\begin {theor} Assume $\mu_0$ is a Hermitian structure on the bundle induced by $\rho$ 
which is parallel in the neighborhood of $Cr(h)$ and $\mu$ is an arbitrary Hermitian structure.
Then
${\cal R} = {\cal R} (M, \rho,\mu, g, \tau)$ is given by

\[\log {\cal R}= (-1)^{n+1}\int_{M\setminus Cr(h)} \theta(\rho, \mu_0)
\wedge X^{\ast}(\Psi(TM,g)) 
+ \int_M V(\rho, \mu, \mu_0) e(M,g)-
\]
\[-\sum_{x \in Cr(h)}(-1)^{\rm{ind} (x)}V(\rho, \mu, \mu_0)(x).\]
\n In the case $M$ is of odd dimension, the above formula for $\log {\cal R}$ simplifies,

\[\log {\cal R}= \int_{M\setminus Cr(h)} \theta(\rho, \mu_0)
\wedge X^{\ast}(\Psi(TM,g))-\sum_{x\in Cr(h)}(-1)^{\rm{ind} (x)}V(\rho, \mu, 
\mu_0)(x).\]

\end {theor}

\n   
Since $\mu_0$ is parallel in the neighborhood of the critical points the first integral in the formula evaluating ${\cal}R$ is convergent, and this is the 
reasons $\mu_0$ is necessary.  The formulas in Theorem 0.1 have been 
discovered by Bismut-Zhang in the case  ${\cal A}=\bbc$ (we choose a slightly different definition of $\theta(\rho,\mu)$ which is half of Bismut-Zhang's).
 Theorem 0.1, applied in this case,  implies the main result in [BZ],
as in this case  ${\cal R} = T_{an} / T_{Reid}.$
The proof of Theorem 0.1 provides a conceptually new proof of the 
formula of [BZ]. In particular Theorem 0.1 also 
shows that the formula of Bismut-Zhang for the quotient of analytic 
and Reidemeister torsion  remains 
true for the quotient of $L_2-$
analytic and Reidemeister torsions, cf subsection 5.1, 
Remark to Theorem 0.1. This answers  
a question, raised by several experts in the field.

\n In [BFKM], we have shown that for a unitary representations $\rho$ of $\Gamma$ with $(M,\rho)$ 
of determinant class,
$T_{an} = T_{Reid}.$  Theorem 0.1 extends this result to an arbitrary representation  and
reformulates this result in a way that it holds without the assumption of $(M, \rho)$
being of determinant class.

\vspace{0.1cm}

\n Denote by $Rep(\Gamma ; \cal W)$ the space of representations of $\Gamma$ on $\cal W.$ The group
$\Gamma$ is a finitely presented, let us say with generators $\gamma_1,...,\gamma_{N_1}$ and relations
$R_j(\gamma_1,...,\gamma_{N_1}) = e_{\Gamma}$  $(1 \le j \le N_2) .$  Hence $Rep(\Gamma; \cal W)$
is a level set of the holomorphic map induced by the relations, 
$Gl_{\cal A}({\cal W} )^{N_1} \to Hom_{\cal A}({\cal W} )^{N_2}$
and thus a complex analytic space (of infinite dimension if 
$dim_{\bbc}\cal W$ is infinite).  
Using ${\cal R} (M, \rho, \mu , g, \tau)$ we will construct in section 5.2
a smooth invariant for a pair of a closed odd dimensional manifold with fundamental group $\Gamma$ and an Euler structure (cf section 5.2 for the definition 
of  Euler structure). This invariant will be a real analytic function on 
$Rep(\Gamma; {\cal W}).$
It will be shown in a forthcoming paper that this function and other 
similarly defined functions associated with even dimensional manifolds, 
are locally the real part of holomorphic functions on  
$Rep(\Gamma; {\cal W}).$ 
\n As an example (cf section 5.2), we will calculate this function in 
the case 
$M=N\times S^1$ with $N$ a simply connected manifold of even dimension
and $E$ the canonical Euler structure on $N\times S^1$. 
In this case, $Rep(\Gamma; {\cal W}) \simeq
Gl_{\cal A} ({\cal W} )$ and the function is equal to (cf Proposition 5.1) 
$-\frac {1}{2} 
\chi (N)\log det (A^*A)^{1/2}$
where $A \in Gl_{\cal A} (\cal W ),$  $det$ is the determinant in the sense of 
Fuglede- Kadison (cf.[BFKM]),
and $\chi(N)$ is the Euler-Poincar\'e characteristic of $N.$
This calculation shows that the invariant is nontrivial. In [BFK 5] we extend 
the results of this paper to the relative torsion for a 
homotopy triangulation.

\vspace{0.5cm}

\n {\large{\bf 0.2 Definitions and outline of the paper}}

\vspace{0.3cm}

\n Let $(M, x_0)$ be a closed manifold with base point $x_0,$ $\rho : \Gamma \rightarrow 
Gl_{\cal A} ({\cal W})$ a representation of $\Gamma := \pi_1 (M; x_0)$ on an 
${\cal A}$-Hilbert module ${\cal W}$,
${\cal E} \rightarrow M$ the vector bundle associated to 
 $\rho$, equipped with the canonical flat connection $\nabla,$
and  $\mu$ a Hermitian structure of ${\cal E}$.

\vspace{0.3cm}

\n {\it deRham complex:} 
Denote by $(\Omega (M; {\cal E}),
d)$ the deRham complex of smooth forms with values in ${\cal E}.$ Here 
$d_k : \Omega^k \rightarrow \Omega^{k + 1}$ is the exterior differential 
determined
by $\nabla.$ A Riemannian metric $g$ on $M$ together with the Hermitian
structure $\mu$ on ${\cal E}$ allow to introduce an inner product on 
$\Omega^k$ and therefore
determine an adjoint $d_k^{\ast}$ of $d_k$. 
Denote by $\Delta_k =
d_k^{\ast} d_k + d_{k-1} d_{k-1}^{\ast}$ the Laplacian on $k$-forms.
These Laplacians can be 
used to introduce $s-$ Sobolev norms on $\Omega.$ We denote by $H_s (\Lambda^k
(M; {\cal E}))$ the completion of $\Omega^k$ with respect to the s-Sobolev norm. This leads
to a family of complexes $H_s(\Lambda (M; {\cal E}))$ of ${\cal A}$-Hilbert 
modules.

\vspace{0.3cm}

\n {\it Analytic torsion} $T_{an}:$
The Laplacians $\Delta_k$ admit a regularized determinant $\det \Delta_k$ in the von Neumann
sense. If  $(M, \rho)$ is of
determinant class (cf Appendix B),  these determinants do not vanish so that  the analytic 
torsion can be defined by
\[\log T_{an}(M,\rho,g,\mu):= 1/2 \sum_q (-1)^{q+1} q \log \det \Delta_q.\] 

(The concept of determinant class was introduced in 
[BFKM] Definition 4.1  and will be reviewed  in Appendix B.)

\vspace{0.3cm}

\n {\it Combinatorial complex:}
A generalized triangulation of $M$ is a pair $\tau = (h,g')$ with
the following properties:

\n {\em (T1) \hsps $h: M \rightarrow \bbr$ is a smooth Morse function;}
\footnote{In [BFKM], the Morse function h was supposed to be selfindexing
$(h(x) = index(x)$ for any critical point $x$ of $h)$. For our proof of Theorem 0.1 (step 2),
it is convenient to drop this condition (cf [BFK2]). }

\n {\em (T2) \hsps $g'$ is a Riemannian metric so that  $-{\rm grad}_{g'} h$ 
satisfies
the Morse - Smale condition (for any two distinct critical points  $x$ and
$y$ of $h,$ the stable manifold $W^+_x$ and the unstable manifold $W^-_y,$ with
respect to $-{\rm grad}_{g'} h,$ intersect transversely);}

\n {\em (T3) \hsps in a neighborhood of any critical point of $h$ one can
introduce local coordinates such that, with $q$ denoting the index of this 
critical point,

\[h(x) = h(0) - (x^2_1 + \ldots + x^2_q) / 2 + (x^2_{q+1} + \ldots + x^2_{d}) 
/ 2 \]

\n and the metric $g'$ is Euclidean in these coordinates.}

\vspace{0.3cm}

\n A generalized triangulation $\tau$  and a 
Let $p : \tilde{M} \rightarrow M$ be the universal covering of $M$ and 
denote by $\tilde{g}$ and $\tilde{\tau} = (\tilde{h}, \tilde{g}')$ the lifts
of $g$ and $\tau$ on $\tilde{M}.$ Denote by $Cr_q (\tilde{h}) \subset 
\tilde{M},$ resp. $Cr_q (h) \subset M,$ the set of critical 
points of index $q$ of
$\tilde{h},$ resp. $h$, and let $Cr(\tilde{h}) = \cup_q Cr_q (\tilde{h}).$ 
Notice that the group $\Gamma$ acts freely on $Cr_q (\tilde{h})$ and the 
quotient set can be identified with $Cr_q(h).$ 

An orientation 
\[{\cal O}_h := \{{\cal O}_{\tilde{x}} | \tilde{x} \in Cr(\tilde{h}) \}\]
is provided 
by a collection of orientations
${\cal O}_{\tilde{x}}$ for
the unstable manifolds $W^-_{\tilde{x}}$ that 
are $\Gamma$-invariant.

\n To the system $(M,\rho,\mu, \tau, {\cal O}_h, )$, we associate a cochain 
complex of finite type over the von Neumann algebra ${\cal A}, \ {\cal C} 
\equiv {\cal C}
(M, \rho, \mu, \tau, {\cal O}_{h,} ) = \{{\cal C}^q, \delta_q \},$ 
where ${\cal C}^q = 
\oplus_{x \in Cr_q (h)} {\cal E}_x,$ which can be identified with the module of
$\Gamma$-equivariant maps $f : Cr_q(\tilde{h}) \rightarrow {\cal W},$ 
and the maps
$\delta_q : {\cal C}^q \rightarrow {\cal C}^{q+1}$ are given by 

\[\delta_q (f) (\tilde{x}) := \sum_{\tilde{y} \in Cr_q (\tilde{h})}
\nu_{q+1}
(\tilde{x}, \tilde{y}) f(\tilde{y}) \]

\n where $\nu_q : Cr_q (\tilde{h}) \times Cr_{q-1} (\tilde{h}) \rightarrow \bbz$
is defined by $\nu_q (\tilde{x}, \tilde{y}) :=$ intersection number 
$(W^-_{\tilde{x}} \cap V, W^+_{\tilde{y}} \cap V)$ with $V := \tilde{h}^{-1}
(q - \frac{1}{2}).$ 
The cochain complex ${\cal C}
(M, \rho, \mu, \tau, {\cal O}_{h} )$ depends on $\mu$ only via $\mu_x,  x\in Cr(h).$
Let $\Delta^{comb}_q := \delta^{\ast}_q \delta_q +
\delta_{q-1} \delta^{\ast}_{q-1}$ denote the combinatorial Laplacian.

\vspace{0.3cm}

\n {\it Reidemeister torsion} $T_{Reid}:$
$\Delta_q^{comb}$ admits a regularized determinant $\det \Delta_q^{comb}$ 
in the von Neumann
sense which, under the additional hypothesis of $(M, \rho)$  being of
determinant class (cf Appendix B), does not vanish. Under this hypothesis, 
the combinatorial torsion can be defined

\[\log T_{comb}(M,\rho,\mu, \tau):= 1/2 \sum_q (-1)^{q+1} q \log \det \Delta_q^{comb}.\]
 It can be shown that the combinatorial torsion is independent of the choice of 
${\cal O}_h$.
Given a Riemannian metric $g$, the metric part $T_{met}$ of the Reidemeister torsion is defined as in [BFKM]
and the Reidemeister torsion $T_{Reid}$ is then given by
 \[ \log T_{Reid} := \log T_{comb} + \log T_{met} \] . 

\vspace{0.3cm}

\n {\it Relative torsion:}  
Notice that the 
integration on the $q$-cells of the generalized triangulation $\tau$ which are
given by the unstable manifolds of $-{\rm grad}_{g'}h,$ defines a morphism 

\[Int :(\Omega, d) \rightarrow ({\cal C}, \delta),\]

\n i.e. $Int_q : \Omega^q \rightarrow {\cal C}^q$ is an ${\cal A}$-linear map so that
$\delta_q Int_q = Int_{q+1} d_q$ (cf Appendix by F. Laudenbach in [BZ] or 
[B1]). We
would like to define the relative torsion of the system $(M, \rho, \mu,  g, 
\tau)$ as the torsion of the mapping cone defined by the integration
map. Unfortunately, the torsion of this mapping cone cannot be defined (at least 
not in an obvious way);  a first difficulty comes from the fact that the 
integration maps $Int_k$ do not extend to closed maps, defined on a dense 
domain of $H_0(\Lambda^k(M; {\cal E})).$ However, for $s$ sufficiently large,
the integration morphism has an extension $Int_s$

\[Int_s : H_s(\Lambda (M; {\cal E})) \rightarrow {\cal C}\]

\n to a bounded morphism. We then consider the composition $g_s$

\[g_s : H_0 (\Lambda(M; {\cal E}) )
\stackrel{(\Delta + Id)^{-s/2}}{\longrightarrow}
H_s (\Lambda (M; {\cal E})) 
\stackrel{Int_s}{\longrightarrow} {\cal C}\]

\n and prove that  $g_s$ induces an
isomorphism in algebraic cohomology

\n $ Ker (d_k) / {\rm range} (d_{k-1}) $
as well as in reduced co\-ho\-mo\-lo\-gy (cf
Pro\-po\-si\-tion 2.5)

$$Ker (d_k) / \overline{{\rm range}
(d_{k-1})}.$$ 

\vspace{0.3cm}

\n This implies that the mapping cone ${\cal C}(g_s),$ defined 
by 
$${\cal C}(g_s)_k := {\cal C}^{k-1}
\oplus H_0(\Lambda^k(M; {\cal E})); \,\,\,d(g_s)_k := \bl {cc} 
- \delta_{k-1} & g_{k;s} \\
0 & d_k \er,
$$ 
is an algebraically acyclic cochain complex (i.e. 
$Ker d(g_s)_k / {\rm range}(d(g_s)_{k-1}) = 0).$
In particular the complex 
${\cal C}(g_s)$ is a cochain complex of ${\cal A}$-Hilbert modules (cf Lemma
1.11) whose
Laplacians $\Delta (g_s)_k$ 
are unbounded operators which
admit a nonvanishing 
regularized determinant $\det \Delta
(g_s)_k$ in the von Neumann  sense. One of the purposes of section 1 is to
establish this result.

\vspace{0.3cm}

\n The relative torsion ${\cal R}_s$ is defined as the torsion 
$T({\cal C}(g_s))$ of the
mapping cone ${\cal C}(g_s),$

\[\log {\cal R}_s := \log T({\cal C}(g_s)) := \frac{1}{2}
 \sum_{k} (-1)^{k+1} k \ \log 
\det(\Delta(g_s)_k).\]

\n In section 2, we show that $\log{\cal R}_s$ is independent of $s$ 
($s$ sufficiently large) and thus provide a well defined number which we denote by
$\log {\cal R}$.

\vspace{0.3cm}

\n As first shown in [CMM] in a slightly different and less general 
version, one verifies that if $(M, \rho)$ is of determinant class, 
then $\log {\cal R} = \log T_{an} - \log T_{Reid}$ (cf (2.15))
where $T_{an},$ resp. $T_{Reid}$, denotes the analytic 
torsion resp. Reidemeister torsion. 
The main result of [BFKM] can thus be stated as follows: When $(M,
\rho)$ is of determinant class and $\mu$ is parallel then 
${\cal R} = 1.$

\vspace{0.3cm}

\n{\it Representations:} Let $\rho: \Gamma \to Gl_{\cal A}({\cal W})$  be a representation
of $\Gamma$ on a ${\cal A}-$Hilbert module of finite type. The representation 
$\rho^{\sharp} : \Gamma \to Gl_{\cal A}({\cal W})$ {\it dual} to $\rho$ is defined by
$\rho^{\sharp}( \gamma ) := \rho( \gamma ^{-1})^{\ast} \, (\gamma \in  \Gamma)$
where $^\ast$ denotes the adjoint in  $Gl_{\cal A}({\cal W})$ with respect to the inner product
of ${\cal W}.$ The representation $\rho$ is said to be unitary if $\rho^{\sharp} = \rho.$
Notice that the flat bundle ${\cal E} \to M$ induced by a unitary representation admits 
a parallel Hermitian structure.

\n A representation  $\rho ': \Gamma \to Gl_{\cal A'}({\cal W '})$ (${\cal W '}$ a ${\cal A}'-$Hilbert module
of finite type) is said to be {\it isomorphic} to $\rho$ if there exist an isomorphism 
$\theta : \cal A \to {\cal A}'$ ( with $tr_{{\cal A}'} (\theta a) = tr_{\cal A} (a) , \forall a \in {\cal A} $) and an isomorphism
$\Theta : {\cal W} \to {\cal W}'$, such that $\Theta$ is $\theta-$linear. As an example, we mention that
$\rho \otimes \rho^{\sharp}$ and  $\rho^{\sharp} \otimes \rho$ are isomorphic.

\n We denote by $Rep (\Gamma, {\cal W})$ the complex analytic space of representations (cf subsection 0.1).

\vspace{0.3cm}

\n{\it Densities and local quantities:} A {\it density}  $\alpha$ on a manifold $M$ of dimension $n$ is a smooth $n-$form
with coefficients in the orientation bundle ${\cal O}_M \to M.$ Densities can be integrated on relative compact
open subsets $U \subset M,$ and the integral is denoted by $\int_U \alpha.$

\n In this paper, we use the word {\it local quantity} for an assignment which to any set of data 
${\cal G} := ( {\cal E} \to M, \nabla, \mu , g, X)$ consisting of a bundle ${\cal E} \to M,$
a flat connection $\nabla,$ a Hermitian structure $\mu,$ a Riemannian metric $g$, and a nonzero
vectorfield $X$ on $M$, provides a density $\alpha ({\cal G})$ with the following properties:

\begin{itemize}
\item[(P1)]  For any open set $U \subset M,$ $\alpha ({\cal G} \vert_U) = \alpha ({\cal G}) \vert_U.$ 

\item[(P2)] If $(\phi, \Phi )$ is an isomorphism between ${\cal G}_1$ and ${\cal G}_2$ (i.e. $\Phi : {\cal E}_1
\to {\cal E}_2$ is a bundle isomorphism above the diffeomorphism $\phi: M_1 \to M_2$ such that
$\phi^{\ast} \nabla _2 = \nabla _1$, $\phi^{\ast} \mu _2 = \mu _1,$ $ \phi^{\ast} g_2 = g_1,$ 
$ \phi^{\ast} X_2 = X_1$)  then  $\phi ^{\ast} \alpha ({\cal G}_2) = \alpha ({\cal G}_1).$


\end{itemize}

To a system ${\cal F} = (M, \rho , \mu , g , \tau )$ we associate the set of data

\n ${\cal G} _{\cal F} = ( {\cal E} \vert _{M \setminus Cr(h)} \to M\setminus Cr(h), 
\nabla , \mu , X$) where $({\cal  E } \to  M , \nabla)$ is the bundle with flat connection $\nabla$,
canonically associated to $\rho$ and $X = - grad_{g'} h.$

\vspace{0.3cm}

Let us briefly comment on the proofs of the results stated in subsection 0.1:
The proof of Proposition 0.1 (cf section 3 and subsection 4.2) uses the Witten deformation
of the deRham complex, given by the differentials $d_k (t) := e^{-th} d_k e^{th}.$
We consider a deformation ${\cal R} (t)$ of the relative torsion ${\cal R}
= {\cal R} (0)$, induced by the Witten deformation of the deRham complex. From
many different possibilities for defining the deformation ${\cal R} (t)$ we choose 
one for which the variation 
$\frac{d}{dt} \log {\cal R}(t)$ can be computed. We define
${\cal R}(t)$ as the torsion of the mapping cone associated to

\[\Big(H_0 (\Lambda (M; {\cal E})), d(t)\Big) 
\stackrel{e^{th}}{\rightarrow} 
\Big(H_0 (\Lambda
(M; {\cal E})), d\Big) \stackrel{g_s}{\rightarrow} ({\cal C}, \delta),\]

\n which will be checked to be independent of $s,$ cf Proposition 2.7 and 
Definition 2.8. 
The variation of $\log {\cal R}(t)$ can be computed to be (cf Theorem 2.1)

\spf

\n {\bf (0.1)} \hsps $\frac{d}{dt} \log {\cal R} (t) = \dim {\cal E}_{x} \cdot
 \int_M h \ e(M,g),$

\spd

\n where $e(M,g)$ is the Euler form of the tangent bundle equipped with the
Levi-Civit\`a connection induced from $g.$ It follows that $\log {\cal R} (t)$
is given by $\log {\cal R} + t \int_M h \ e(M,g).$

\spd

\n To continue our discussion we say 
that a function $G: \bbr\rightarrow \bbr$ admits an asymptotic 
expansion for $t \rightarrow \infty$ if there exist a sequence $i_1 > i_2 > 
\ldots > i_N = 0$ and constants $(a_k)_{1 \le k \le N}, (b_k)_{1 \le k \le N}$
such that, as $t \rightarrow \infty,$

\spf

\n {\bf (0.2)} \hsps $G(t) = \sum^N_1 a_k t^{i_k} + \sum^N_1 b_k t^{i_k} \log t + 
o(1).$

\spd

\n For convenience, we denote by  $FT(G(t))$ or $FT_{t=\infty} (G)$ 
the coefficient $a_N$ in the
asymptotic expansion of $G(t)$ and refer to it as the free term of the 
expansion. Notice that the equality (0.1) implies that $\log {\cal R} (t)$
admits an asymptotic expansion of the form (0.2) with $\log {\cal R}$ as the free
term.

\spd

\n By  different considerations outlined below, we derive the existence of an
asymptotic expansion for $\log {\cal R} (t)$ of the form (0.2) and calculate
the free term of this expansion as the integral on $M\setminus Cr(h)$
of a local density expressed in terms of $g,
\mu$ and $\tau.$  The proof of Proposition 0.1 will then 
be derived from the comparison 
of these two calculations.
Using Witten deformation  we
decompose $\log {\cal R} (t)$ into two parts
\[ \log {\cal R} (t) = \log
{\cal R}_{sm} (t) + \log T_{la} (t),\] and show that each of the two parts
admit an asymptotic expansion of the form (0.2) whose free terms can be 
computed. 
Precisely, we consider the Witten deformation
$(\Omega (M; {\cal E}), d(t))$ of the deRham com\-plex, $d(t) = e^{-th} d e^{th}.$ 
For $t$ sufficiently large, the deformed deRham complex can be decomposed,

\[\Big( \Omega (M, {\cal E}), d(t)\Big) = \Big( \Omega_{sm} (t), d(t)\Big) \oplus
\Big( \Omega_{la} (t), d(t)\Big),\]

\n corresponding to the small and the large part of the spectrum of the 
deformed Laplacians $\Delta_k (t).$ This decomposition induces the decomposition

\[(H_0 (\Lambda (M; {\cal E})), d(t)) = (\Omega_{sm} (t), d(t)) 
\oplus (H_{0, la} (t), d(t)).\]

\n The complex $(\Omega_{sm} (t),
d(t))$ is of finite type. Denote by $g_{s, sm} (t)$  the restriction of 
$g_s \cdot e^{th}$ to $\Omega_{sm} (t)$, where 
$e^{th} : \Omega (M, {\cal E}) \rightarrow 
\Omega (M, {\cal E})$ denotes the mul\-ti\-pli\-ca\-tion by $e^{th}$.
The mapping cone 
${\cal C}(g_{s,sm} (t))$ 
is a cochain complex of ${\cal A}$-Hilbert modules of finite type which is
algebraically acyclic and has a well defined torsion
denoted by ${\cal R}_{sm} (t).$ As in [BFKM, section 6],
one verifies that $(H_{0, la} (t), d(t))$
is algebraically acyclic and
has a well defined torsion, denoted by
$T_{la} (t).$ We prove that 
$\log {\cal R} (t) = \log
{\cal R}_{sm} (t) + \log T_{la} (t)$  (cf subsection 3.1, statement $(B)$).

\spd

\n To prove that $\log {\cal R}_{sm} (t)$ has an asymptotic expansion, 
we show that
${\cal R}_{sm} (t) = T({\cal C}(f(t)))$ (cf Proposition 3.1), the torsion 
of the mapping cone induced by

\[f(t) : \Big( \Omega_{sm} (t), d(t)\Big) \stackrel{e^{th}}{\rightarrow} 
(\Omega, d) 
\stackrel{Int}{\rightarrow} ({\cal C}, \delta)\]

\n and use the Witten-Helffer-Sj\"ostrand theory 
as presented in [BFKM]
to prove that 
$\log T
({\cal C}(f(t)))$ admits an asymptotic expansion of the type (0.2) and to
compute its free term (Proposition 3.1).
To analyze $\log T_{la} (t)$, we proceed as in [BFK1] or [BFKM].
We derive the existence of an asymptotic expansion for $\log T_{la} (t)$ and 
a closed 
formula for its free term from Proposition 3.2 iii). This result  states the 
existence of an 
asymptotic expansion for the difference $\log T_{la}(t)- \log \tilde{T}_{la}(t)$
where  $(M, \rho, \mu, g, \tau)$ and 
$(\tilde{M}, \tilde{\rho}, \tilde{\mu}, \tilde{g}, \tilde{\tau})$ 
are two systems whose Morse functions 
$h$ and $\tilde h$ have the same number of critical points in each index, 
and provides a formula 
for
its free term. 

\vspace{0.1cm}

\n The proof of Theorem 0.1  ( cf subsection 5.1) uses  
the anomalies of the relative torsion (Propositions 4.1- 4.3) and
Proposition 0.1.  The locality of the density given Proposition 0.1 is essential for this proof. 

\spd

\n The paper is organized in 4 sections and two appendices.

\vspace{0.1cm}

\noindent Section 1: In subsection 1.1 we single out a  class of unbounded
${\cal A}$-linear operators on ${\cal A}$-Hilbert modules  ( ${\cal A}$ is a 
von Neu\-mann algebra with finite trace )
for which the regularized determinant can be defined and  in subsection 1.2 
a class of complexes of ${\cal A}$-Hilbert modules for which the
Laplacians are operators in the above class. Therefore the torsion of an
(algebraically) acyclic 
complex of such type can be defined. These abstract results are needed 
because  
the mapping cone of (a regularized version of) the integration
map  between the deRham complex and the combinatorial complex 
has as its components  direct sums of  ${\cal A}$-Hilbert modules of finite type and  
completions (with respect to a Sobolev norm) of  spaces of smooth sections in 
smooth  bundles of ${\cal A}$-Hilbert modules of finite type (cf. [BFKM] for 
definitions). We  show that the mapping cone of a regularized
version of the integration map is  a complex of the type introduced in 
subsection 1.2.

\vspace{0.1cm}

\n Section 2: Using the results of section 1, we introduce in subsection 2.1 the notion
of relative torsion ${\cal R} $ and in subsection 2.2 
the Witten deformation  ${\cal R} (t)$ of the relative torsion ${\cal R}$ and calculate 
its variation. 

\vspace{0.1cm}

\n Section 3: Proof of  Proposition  0.1 in the case $g=g'.$

\vspace{0.1cm}

\n Section 4: We investigate how the relative torsion
${\cal R}(M,\rho , \mu , g, \tau)$ varies with respect to the Riemannian 
metric $g$, the Hermitian structure $\mu$ and the 
triangulation $\tau.$  Further we prove Proposition 0.1 in full generality.

\vspace{0.1cm}

\n Section 5: Proof of  Theorem  0.1 and the construction of an invariant for odd dimensional manifolds.

\vspace{0.1cm}

\n Appendix A: For the
convenience of the reader we present  a proof of a slightly
stronger version of a Lemma due to Carey-Mathai-Mishchenko.

\vspace{0.1cm}

\n  Appendix B:  We review the concept
 of determinant class and 
provide a simple example of a pair $(M,\rho )$, $M:=S^1$ and $\rho$ 
a representation of $\pi_1 (S^1) = \bbz$ on $l_2(\bbz )$, which is not of 
determinant class. The existence of such pairs motivates in part the concept 
of relative torsion.

\vspace{0.1cm}

\n Throughout the paper the notions of trace (denoted by Tr), 
dimen\-sion, de\-ter\-mi\-nant etc. are always understood in the von Neu\-mann 
sense.

\section{Operators and complexes}

\subsection{Operators}

\vspace{0.3cm}

\noindent Let ${\cal W}_1, {\cal W}_2$ be ${\cal A}$-Hilbert modules and 
$\varphi : {\cal W}_1 \rightarrow
{\cal W}_2$ an operator. Denote its domain by 
${\rm domain}(\varphi) \subseteq {\cal W}_1$ and its
{\rm range} by ${\rm range}(\varphi) = 
\varphi \Big({\rm domain}(\varphi)\Big) \subseteq {\cal W}_2.$ 

\spd

\n Introduce the following properties of $\varphi.$

{\it \n $Op (1) \hsps \varphi$ is ${\cal A}$-linear.

\n $Op (2) \hsps \varphi$ is densely defined, i.e. $\overline{{\rm domain} (\varphi)} =
{\cal W}_1.$

\n $Op (3) \hsps \varphi$ is closed (i.e. the graph $\Gamma (\varphi)$ of $\varphi$ is
closed in ${\cal W}_1 \times {\cal W}_2$).}

\spd

\n If, instead of $Op(3), \varphi$ satisfies

\spf

\n {\em Op(3)' \hsps $\varphi$ is closable}

\spd

\n we consider the minimal extension of $\varphi$ and denote it again by 
$\varphi$.

\n In the sequel, we will not distinguish between property $Op(3)$
and $Op(3)'.$

\spd

\n An operator $\varphi$ satisfying $Op(1)-Op(3),$ admits an adjoint, 
$\varphi^{\ast}$ (cf. e.g. [RS, p 316]) and we can consider the nonnegative, 
selfadjoint operator $\varphi^{\ast} \varphi.$ Functional calculus permits
to define the square root $|\varphi|:= (\varphi^{\ast} \varphi)^{1/2}.$ Then,
for any $0 < t < \infty,$ the heat evolution operator $e^{-t|\varphi|}$ is 
bounded, nonnegative and selfadjoint.
The operators $|\varphi|$ and $e^{-t|\varphi|}$ satisfy the properties
$Op(1) - Op(3).$
In the sequel, if not mentioned otherwise, we assume that the operators 
con\-si\-de\-red sa\-tis\-fy $Op(1) - Op(3).$

\spd

\n Denote by  $dP(\lambda) \equiv dP_{|\varphi|} (\lambda)$ the (operator 
valued) spectral measure associated to $|\varphi|$ defined by the orthogonal 
projectors $P_{|\varphi|} ([0, \lambda]) \
(0 \le \lambda \le \infty)$ 

\[P_{|\varphi|} ([0, \lambda])  = \int^{\lambda_+}_{0_-} 
dP_{|\varphi|} (\lambda).
\]

\n Most of the operators $\varphi$ considered in this paper will satisfy the
fol\-lo\-wing important property

\spf

\n {\em $Op (4) \hsps Tr P_{|\varphi|} ([a,b]) < \infty  \ \ \ \mbox{(for any} \ 0 < 
a \le b < \infty \mbox{)}$}

\spd

\n where

\spf

\n {\bf (1.1)} \hsps $
P_{|\varphi|} ([a,b]) := \int^{b_+}_{a_-} dP_{|\varphi|} 
(\lambda).$

\spd

\n If $\varphi$ satisfies $Op(4)$ one can define the Stiltjes measure 
$dF_{|\varphi|} (\lambda)$ on the half line $(0, \infty),$ given by $(0 < a <
b < \infty)$

\spf

\n {\bf (1.2)} \hsps $
\int^{b_+}_{a_-} dF_{|\varphi|} (\lambda) = Tr
P_{|\varphi|} ([a,b]).$

\spd

\n The next two properties concern the asymptotic behavior of $TrP([a,b])$ for
$b \nearrow \infty$ and $a \searrow 0:$

\spf

\n {\em $Op(5)$ \hsps There exists $\alpha \ge 0$ such that}

\[Tr P_{|\varphi|} ([1, \lambda]) = 0 (\lambda^{\alpha}) \ \mbox{for} \ \lambda
\nearrow \infty\]

\n and

\spf

\n {\em $Op(6) \hsps Tr P_{|\varphi|} 
([\lambda, 1]) = 0(1) \ \mbox{for} \ \lambda \searrow
0.$}

\spd

\n If $\varphi$ satisfies properties $Op(4)$ and $Op(6)$ one can define the
spectral distribution functions, associated to $|\varphi|,$

\spf

\n {\bf (1.3)} \hsps $F_{|\varphi|}^+ (\lambda) = Tr P_{|\varphi|}( ( 0,
\lambda ))$; $F_{\vert\varphi\vert}(\lambda ):= TrP_{\vert\varphi\vert} 
((-\infty ,\lambda ))$.

\spd

\n The spectral distribution function $F^+_{\vert\varphi\vert}(\lambda )$ can 
be described variationally as follows,
$(\lambda > 0)$ 

\spf

\n {\bf (1.4)} \hsps
$F_{|\varphi|}^+ (\lambda) = sup \{ \dim L | L \subset \  {\rm domain}
(|\varphi|),  \ 
\perp  Ker\varphi,$
is an ${\cal A}$-Hilbert 
\linebreak
\hspace*{3.2cm} sub\-mo\-du\-le 
with $||  |\varphi| (x)
|| < \lambda ||x|| \ \forall x \in L \}$

\spd

\n (cf e.g. [GS], [BFKM]).

\spd

\n If, in addition, $\varphi$ satisfies also $Op(5)$, one can define the heat
trace $\theta_{|\varphi |} (t)$ associated to $\varphi$ (excluding zero modes)

\spf

\n {\bf (1.5)} \hsps $
\theta_{|\varphi|} (t) := \int^{\infty}_{0_+} e^{-t\lambda}
dF_{|\varphi|}^+ (\lambda) \ (t > 0).$

\spd

\n Using integration by parts in the Stiltjes integral (1.5) one concludes
from $Op(5)$ that

\spf

\n {\bf (1.6)} \hsps $
\theta_{|\varphi|} (t) = 0 (t^{-\alpha}) \
(t \searrow 0). $

\spd

\n Next, let us introduce the `partial' zeta function
$\zeta^{\rm{I}}_{|\varphi|} (s)$ associated to the heat trace $\theta_{|\varphi |}
(t)$ and
defined for $s \in \bbc$ with $Res > \alpha$ (with $\alpha$ given by 
Op(5))

\spf

\n {\bf (1.7)} \hsps $
\zeta^{\rm{I}}_{|\varphi|} (s) := \frac{1}{\Gamma (s)}
\int^{1}_{0} t^{s-1} \theta_{|\varphi|}(t) dt$

\spd

\n where $\Gamma (s)$ denotes the gamma function. Notice that
$\zeta^{\rm{I}}_{|\varphi|} (s)$ is holo\-mor\-phic in the halfplane $Res >
\alpha.$

\spd

\n The following property is needed to define the notion of regularized
determinant:

\spd

\n {\em $Op(7) \hsps \zeta^{\rm{I}}_{|\varphi|} (s)$ has an analytic 
continuation at $s =
0.$}

\spd

\n A sufficient condition for $Op(7)$ to hold is the exi\-sten\-ce of an
asymp\-to\-tic ex\-pan\-sion of $\theta_{|\varphi|} (t)$ near $t=0$ of 
the form 

\spd

\begin{itemize}
\item[(Asy)] $\theta_{|\varphi|} = \sum^{m-1}_{j=0} a_j t^{-\alpha_j} + a_m +
R(t)$ where $m \in \bbz_{\ge 0}, $ $ 0 < \alpha_{m-1} 
\linebreak
< \ldots < \alpha_0, \ a_0,
\ldots , a_m \in \bbr$ and where $R(t)$ is a bounded, continuous function
with $R(t) = 0(t^{\rho})$ for some $\rho \in \bbr_{>0}.$
\end{itemize}
\spd

\n The fact that (Asy) implies $Op(7)$ follows from the following 
lemma, using integration by parts.

\spf

\begin{lemma} Let $f: (0,1] \rightarrow \bbr$ be a continuous function of the
form

\[f(t) = \sum^{m-1}_{j=0} a_j t^{-\alpha_j} + a_m + R(t)\]

\n with $a_0, \ldots , a_m, \alpha_0, \ldots , \alpha_{m-1} $ and $R(t)$ as
in (Asy).

\spd

\n Then, the holomorphic function
$\xi(x) := \frac{1}{\Gamma(s)} \int^1_0 t^{s-1} f(t) dt,$ defined for $Res
> \alpha_0,$ has a meromorphic con\-tinua\-tion to the half plane $Res > - 
\rho$ with $\rho >0$ as in $(Asy)$. It
has only simple poles, located at $s = \alpha_j \ (0 \le j \le m -1)$. In
particular, $\xi (s)$ is regular at $s = 0$ and $\xi (0) = a_m$.
\end{lemma}

\spd

\n Some of the operators $\varphi$ we will consider have the property that 
$0$ is not in the spectrum, i.e.,

\spf

\n {\em $Op(8)$ \hsps There exists $\varepsilon > 0$ with the property 
that $Tr P_{|\varphi|} ([0, \varepsilon]) = 0.$}

\spd

\n If $\varphi$ satisfies the properties $Op(1)-Op(8),$ then

\spf

\n {\bf (1.8)} \hsps $
\theta_{|\varphi|} (t) = 0(e^{-t\varepsilon/2}) \ \ \mbox{for} \ t
\rightarrow \infty.$

\spd

\n Thus $\int^{\infty}_1 \ t^{s-1} \theta_{|\varphi|} (t) dt$ is an entire function of
$s$ and, as a consequence,

\spf

\n {\bf (1.9)} \hsps
$\zeta^{\rm{II}}_{|\varphi|} (s) := \frac{1}{\Gamma(s)}
\int^{\infty}_1 t^{s-1} \theta_{|\varphi|} (t) dt$

\spd

\n is a holomorphic function in  $s \in \bbc$ with

\spf

\n {\bf (1.10)} \hsps
$\zeta^{\rm{II}}_{|\varphi|} (0) = 0.$

\spd

\n (Actually for $\zeta^{{\rm II}}_{|\varphi|} (s),$ no zeta regularization is
needed (cf [BFKM], subsection 6.1).)

\spd

\n For operators $\varphi$ satisfying properties $Op(1)-Op(8),$ one can
introduce the zeta function $\zeta_{|\varphi|} (s)$ 

\[\zeta_{|\varphi|} (s) := \zeta^{\rm{I}}_{|\varphi|} (s) +
\zeta^{\rm{II}}_{|\varphi|} (s) = \frac{1}{\Gamma (s)} \int^{\infty}_0 t^{s-1}
\theta_{|\varphi|} (t) dt\]

\n for $Res > \alpha$ which has an analytic continuation at $s=0$. This allows
to define the volume of the operator $\varphi,\  Vol \varphi,$ by

\spf

\n {\bf (1.11)} \hsps $
log Vol (\varphi) := log det |\varphi| := -
\frac{d}{ds}|_{s=0} \ \zeta_{|\varphi|} (s).$

\spf

\begin{definition} 

\begin{enumerate}
\item An operator $\varphi: {\cal W}_1 \rightarrow {\cal W}_2$
is said to be of sF type (strong Fredholm type) if
\begin{itemize}
\item[(i)] $\varphi$ satisfies $Op(1) - Op(6)$ and
\item[(ii)] $dim (Ker \varphi) := Tr P_{|\varphi|} (\{0\}) < \infty.$
\end{itemize}
\item An operator $\varphi$ is a bounded operator of trace class if
\begin{itemize}
\item[(i)] $\varphi$ is bounded (hence satisfies $Op(1) - Op(3))$;
\item[(ii)] $\varphi$ satisfies $Op(4)$;
\item[(iii)] $|| \varphi ||_{tr} := \int^{\infty}_{0_+} \lambda
dF_{|\varphi|} (\lambda) < \infty$.
\end{itemize}
\item An operator $\varphi$ is said to be $\zeta$-regular if $\varphi$
satisfies $Op(1)-Op(7).$
\item An operator $\varphi$ is a bounded operator of finite rank if
\begin{itemize}
\item[(i)] $\varphi$ is bounded; 
\item[(ii)] $dim\Big(\overline{{\rm range}(\varphi)}\Big) < \infty.$
\end{itemize}
\n (As a consequence, $\varphi$ satisfies $Op(4) - Op(6)$ and (Asy) 
with $m=0$ and $\rho = 1$ and
thus, $Op(1)-Op(7)$ hold.)
\end{enumerate}
\end{definition}

\n Notice that an operator $\varphi : {\cal W}_1 \rightarrow {\cal W}_2$ of
sF type might not be bounded in the case where $dim {\cal W}_1 = \infty.$

\spf

\begin{prop} (cf [Di]) Let $u : {\cal W}_1 \rightarrow {\cal W}_2$ and $v
: {\cal W}_2 \rightarrow {\cal W}_3$ be bounded operators with $u$,
respectively $v$, of trace class. Then the composition $v u$ is a bounded
operator of trace class and $|| v u ||_{tr} \ \le \ || v
|| \ ||u||_{tr}, respectively,
||vu|| \
\le \ ||v||_{tr} \ ||u||.$
\end{prop}

\n As the proof of Proposition 1.3 is fairly standard, we omit it. The
following result will be used in the proof of Proposition 1.5, stated below:

\spd

\begin{lemma} Assume that $\varphi : {\cal W} \rightarrow {\cal W}$ is a 
nonnegative, selfadjoint operator (i.e $\varphi = |\varphi|$) of sF
type. Let $0<a<b$ and assume that $\varepsilon >0$ satisfies 
$\varepsilon < F_{\varphi} (a)$.  Then one can choose a smooth function $f \in
C^{\infty}_0 ({\bbr}_{\ge 0} ; {\bbr})$ such that $f(\varphi)$ is a bounded
operator of finite rank (hence $f (\varphi)$ is of trace class and 
commutes with $\varphi$) with the following properties:

\begin{itemize}
\item[(i)] $F_{\varphi + f(\varphi)} (\lambda) = 0$ for $\lambda < a$ (in
particular, $\varphi + f(\varphi)$ is 1-1);
\item[(ii)] $F_{\varphi + f(\varphi)} (\lambda) \ge \varepsilon$
for $\lambda \ge a;$
\item[(iii)] $F_{\varphi + f(\varphi)} (\lambda) = F_{\varphi} (\lambda)$ for
$\lambda \ge b$.
\end{itemize}
\end{lemma}

\n {\bf Remark} 
This lemma will be applied in the case when $\varphi$ is a pseudodifferential
operator. Then $\varphi$ can
be chosen so that $f(\varphi)$ is a smoothing operator.

\spf

\n {\bf Proof} Let $g : [0, \infty) \rightarrow {\bbr}$ be a smooth,
increasing function defined on $[0, \infty)$ with $g(\lambda) = a$
for $\lambda \le a$, $g(\lambda )\le\lambda$ for $a\le\lambda\le b$ and 
$g(\lambda) = \lambda$ for $\lambda \ge b.$ Let
$f(\lambda) := g(\lambda) -\lambda$ and define

\[f(\varphi) := \int^{\infty}_0 f(\lambda) dP_{\varphi} (\lambda).\ \
\ \mbox{\carre}\]

\spd

\begin{prop} \hspace{8cm} 

\n (A) If $\varphi : {\cal W} \rightarrow {\cal W}$ is of sF type
and $u : {\cal W} \rightarrow {\cal W}$ is a bounded operator of sF
type, then $\varphi + u$ is of sF type. 

\n (B) If, in addition, $\varphi$ and $\varphi + u$ are selfadjoint,
nonnegative operators and $u$ is of trace class, then there exists
$\varepsilon > 0$ so that the difference $\zeta^{{\rm I}}_{\varphi} (s) -
\zeta^{{\rm I}}_{\varphi + u} (s),$ defined for $Res \gg 0,$

\[\zeta^{{\rm I}}_{\varphi} (s) - \zeta^{{\rm I}}_{\varphi+u} (s) =
\frac{1}{\Gamma (s)} \int^1_0 t^{s-1} \Big(\theta_{\varphi} (t) -
\theta_{\varphi+u} 
(t) \Big) dt \]

\n has an analytic continuation to $Res > - \varepsilon$ and its value at $s =
0$ is equal to $dim (Ker\varphi) - dim(Ker(\varphi + u))$. 
\end{prop}

\spf

\n {\bf Remark} If $v$ is a bounded operator and $\varphi$ and $u$ are
operators as in Proposition 1.5(B), one can, instead of $\theta_{\varphi}(t)$ and
$\theta_{\varphi+v} (t)
,$ consider the function $\theta_{\varphi,v} (t) :=
Trve^{-t\varphi} \Big(Id - P_{\varphi} (\{0\})\Big)$ and similarly,
\n $\theta_{\varphi+u,v} (t) := Tr v e^{-t(\varphi + u)} \Big(Id - P_{\varphi + u}
\linebreak
(\{0\})\Big).$ Notice that for $v = Id, \ \theta_{\varphi, Id} (t) = 
\theta_{\varphi} (t) \
(t > 0).$ Both expressions, $\zeta^{{\rm I}}_{\varphi, v} (s) :=
\frac{1}{\Gamma (s)} \int^1_0 t^{s-1} \theta_{\varphi, v} (t) dt$ 
and $\zeta^{{\rm
I}}_{\varphi + u,v} (s) := \frac{1}{\Gamma(s)} \int^1_0 t^{s-1} 
\theta_{\varphi +
u, v} (t) dt,$ define holomor\-phic functions for $Res \gg 0$ and by the same
arguments as in the proof of Proposition 1.5(B) (cf below) their difference is
holomorphic for $Res > - \varepsilon$ for some $\varepsilon > 0$ and its value
at $s = 0$ is $Tr(vP_\varphi\{ 0\})-Tr(v P_{\varphi +u}\{ 0\})$.

\spf

\n {\bf Proof of Proposition 1.5} (A) Since $\varphi$ satisfies $Op(1)-Op(3)$
and $u$ is bounded, $\varphi+u$ satisfies $Op(1)$ and $Op(2)$ as well. One verifies
in a straightforward way that the graph of $\varphi+u$ is closed (and thus
$Op(3)$ is valid) and the null space $Ker (\varphi+u)$ has finite (von Neumann)
dimension.
It remains to check that $Op(4)$ and $Op(5)$ hold for $\varphi+u$. To see
it, notice that, again by the variational characterization of the spectral
distribution function, $F_{\varphi + u} (\lambda) \le F_{\varphi} (\lambda +
||u||) < \infty \ \forall \lambda.$ Thus

\[F_{\varphi+u} (\lambda) = 0 (\lambda^{\alpha}) \ \mbox{for} \ \lambda
\nearrow \infty
\]

\n where $\alpha$ is given by $Op(5)$, satisfied by $\varphi$, and

\[F_{\varphi+u} (\lambda) = 0 (1) \ \mbox{for} \ \lambda \searrow 0\]

\n as $\varphi$ satisfies $Op(4)$ and $Op(5)$.

\n (B) First we prove (B) in the case where $u$ is of the form
$u=f(\varphi)$ with $f \in C^{\infty} \Big([0, \infty) ; {\bbr} \Big)$ being
of  compact support and such that $\varphi + f(\varphi)$ is 1-1. We claim that
$h(t) := \theta_{\varphi+u} (t) - \theta_{\varphi} (t)$ is real analytic near $t = 0$
and satisfies $h(0) = \dim (\textrm{Ker}\varphi )$. From Lemma 1.1 one then 
concludes that the
difference $\zeta^{{\rm I}}_{\varphi} (s) - \zeta^{{\rm I}}_{\varphi +
f(\varphi)} (s) = \frac{1}{\Gamma(s)} \int^1_0 t^{s-1} h(t)dt$  is well defined
and holomorphic in $s$ for $Res > \alpha,$ has an analytic continuation to
$Res > -1$ and takes value $\dim (\textrm{Ker}\varphi )$ at $s = 0$.
To prove that $h(t)$ is real analytic in $t$ near $t = 0$ and satisfies 
$h(0) = \dim (\textrm{Ker}\varphi )$ we argue as follows: Since $\varphi$ and 
$f(\varphi)$ commute,
$e^{-t(\varphi + f(\varphi))} -e^{-t \varphi} = e^{-t \varphi}
(e^{-tf(\varphi)} -1).$ By assumption, $f$ has compact support, i.e. $supp f
\subseteq [0,K]$ for some $K > 0.$ Therefore 

\[f(\varphi) = \int^{\infty}_0 f(\lambda) dP_{\varphi} (\lambda) = \int^K_0
f(\lambda) dP_{\varphi_K} (\lambda) = f(\varphi_K)\]

\n where $\varphi_K := \int^K_0 \lambda dP_{\varphi} (\lambda)$ is $\varphi$
 re\-stric\-ted  to 
$L := \overline{{\rm range} P ([0,K])}$.

\n Using that $F_{\varphi_K} (\lambda)$ is constant for
$\lambda > K$, we conclude, integrating by parts,

\[ \begin{array}{lll}
h_1(t) & := & Tr e^{- t \varphi_K} (e^{-tf(\varphi_K)} -1) = \int^K_0 e^{-t
\lambda} (e^{-tf(\lambda)} -1) dF_{\varphi_K} (\lambda) = \\
& = & e^{-t \lambda} (e^{-tf(\lambda)} -1) F(\lambda) \Big|^K_0 \\
&&  -t \int^K_0
e^{-t\lambda} \Big(e^{-tf(\lambda)} -1 + e^{-tf(\lambda)} f'(\lambda)\Big)
F_{\varphi_K} (\lambda) d\lambda
\end{array}\]

\n which is real analytic in $t$ and satisfies $h_1(0) = 0$.  Further, one obtains
\begin{eqnarray*}
h_2(t) &:= &Tr (e^{-t(\varphi_K+f(\varphi_K))} P_{\varphi_K+f(\varphi_K)} (\{ 0\})
-e^{-t\varphi_K} P_{\varphi_K}(\{ 0\}))\\
&= &e^{-t f(0)} \dim (\textrm{Ker} (\varphi +f(\varphi)))-\dim (\textrm{Ker}\varphi 
)=-\dim (\textrm{Ker}\varphi )
\end{eqnarray*}

\n which is real analytic in $t$ as well and satisfies $h_2(0)=-\dim\textrm{Ker}
(\varphi )$.  We conclude that $h(t)=h_1(t)-h_2(t)$ is real analytic in $t$ and 
$h(0)=\dim (\textrm{Ker}\varphi )$. \newline
To prove (B) in the general case it suffices, in view of the first step
and Lemma 1.4, to consider operators $\varphi$ and $u$ so that, in addition,
$\varphi$ and $\varphi + u$ are 1-1 and satisfy $Op(8).$
\n These additional properties allow to represent $\zeta^{{\rm I}}_{\varphi}
(s)$ and $\zeta^{{\rm I}}_{\varphi + u} (s)$ by a contour integral as follows:
Choose $\varepsilon > 0$ so that $F_{\varphi} (\lambda) = F_{\varphi + u}
(\lambda) = 0$ for $0 \le \lambda < \varepsilon$ and consider the contour
$\Gamma_{\varepsilon/2} = \Gamma_- \cup \Gamma_0 \cup \Gamma_+$ defined by 

\spd

\[\Gamma_- := \{ z = re^{-i \pi} | \infty > r \ge 0 \} ; 
\Gamma_0 := \{ z = \frac{\varepsilon}{2} e^{i \mu} | - \pi \le \mu \le \pi \}; \] 
\[\Gamma_+ = \{ z = re^{i \pi} | \varepsilon/2 \le r < \infty \}.\]

\spd

\n For $\lambda \in \Gamma_{\varepsilon/2}, \lambda - \varphi$ and $\lambda -
\varphi - u$ are both invertible and one computes

\[\begin{array}{lll}
R(\lambda) &:= & (\lambda - \varphi - u)^{-1} - (\lambda - \varphi)^{-1} = \\
& = & (\lambda - \varphi)^{-1} u (\lambda - \varphi - u)^{-1}.
\end{array} \]

\n As $u$ is a bounded operator of trace class we conclude by Proposition 1.3 
that, for $\lambda \in \Gamma_{\varepsilon/2}, R(\lambda)$ is a bounded 
operator of trace class as well. By a standard argument one shows that for 
$\lambda \in \Gamma_{\varepsilon/2},$

\spf

\n {\bf (1.12)} \hsps $
||(\lambda - \varphi)^{-1}|| \le \frac{1}{|\lambda-
\varepsilon|}~~~~~(\le \frac{1}{\varepsilon/2})$

\spd

\n {\bf (1.12')} \hsps $||(\lambda - \varphi - u)^{-1}|| \le
\frac{1}{|\lambda - \varepsilon|}~~~~~(\le \frac{1}{\varepsilon/2}).$

\spd

\n These estimates allow to represent the difference $A(s) := \zeta^{{\rm
I}}_{\varphi} (s) - \zeta^{{\rm I}}_{\varphi + u} (s)$ by the following
contour integral

\spf

\n {\bf (1.13)} \hsps $
A(s) = \frac{1}{2 \pi i} \int_{\Gamma_{\varepsilon/2}}
\lambda^{-s} Tr R(\lambda) d\lambda = A_0 (s) + A_+ (s) + A_- (s)
$

\spd

\n where

\spf
\n {\bf (1.14$\mbox{)}_{\pm}$} \hsps $A_{\pm} (s) = \frac{1}{2 \pi i}
\int_{\Gamma_{\pm}} \lambda^{-s} Tr R (\lambda) d \lambda.$

\spf

\n {\bf (1.14$\mbox{)}_{0}$} \hsps $A_0 (s) = \frac{1}{2 \pi i} \int_{\Gamma_0}
\lambda^{-s} Tr R (\lambda) d \lambda.$

\spd

\n Notice that $A_0 (s)$ is holomorphic function for $s \in {\bbc}$ and
$A_{\pm} (s)$ are holomorphic functions in a neighborhood of $s = 0$ due to
the fact that $\int^{\infty}_{\varepsilon/2} \frac{x^{-s}}{(x+\varepsilon)^2}
dx$ is absolutely convergent for $|s|$ sufficiently small. Moreover, $A_0 (0)
= 0$ and $A_+ (0) + A(0) = 0.$ \ \ \ \carre

\spf

\begin{exam} \hspace{8cm}
\end{exam}
\n 1.  $\varphi : {\cal W}_1 \rightarrow {\cal W}_2$  bounded and $dim {\cal
W}_1 < \infty$:  then $\varphi$ satisfies $Op(1)-Op(7)$ and is of trace class.
In fact, $Op(6)$ follows from (Asy) which holds with $m=0$ and $\rho = 1.$ Since
$Ker \varphi \subseteq {\cal W}_1$ is of finite dimension, $\varphi$ is of
sF type.

\spd

\n 2.  $\varphi : {\cal W}_1 \rightarrow {\cal W}_2$ is 
bounded and $dim {\cal
W}_2 < \infty:$ then the adjoint $\varphi^{\ast}$ of $\varphi$ is as in
Example (1). Therefore $\varphi$ satisfies $Op(1)-Op(6)$ and (Asy) with $m=0,
\rho = 1.$ Moreover $\varphi$ is of trace class. However the dimension of the
nullspace $Ker \varphi$ might not be finite. Therefore, it is useful to reduce
$\varphi$ to $\hat{\varphi} : {\cal W}_{1}/Ker \varphi \rightarrow {\cal
W}_2.$ Then $\hat{\varphi}$ is injective and bounded and $dim ({\cal W}_1 /
Ker \varphi) < \infty.$

\spd

\n 3. $\varphi$  of finite rank: then the reduced $map \ \hat{\varphi} :
{\cal W}_1 / Ker \varphi \rightarrow {\cal W}_2$ is bounded and injective and
fits either into (1) or (2).

\spd

\n 4. Suppose (M,g) is a closed Riemann manifold and ${\cal E}
\rightarrow M$ is a smooth bundle of ${\cal A}$-Hilbert modules of finite
type, equipped with a Hermitian structure $\mu$ and a smooth connection
$\nabla$ so that the ${\cal A}$-action is parallel. Using $g, \mu$
and $\nabla$ one can define inner products $\langle \cdot , \cdot \rangle_s$ of
Sobolev type on the space of smooth sections on ${\cal E}, C^{\infty} ({\cal
E}).$ (For $s = 0,$ the connection $\nabla$ is not needed.) The
completion of $C^{\infty} (E)$ with respect to $\langle \cdot , \cdot 
\rangle_{s}$ is a
Hilbert space denoted by $H_s ({\cal E}).$ Different choices of $g, \mu,
\nabla$ lead to the same topological vector spaces $H_s ({\cal E}),$
 with different inner products, but equivalent norms.

\spd

\n A pseudodifferential operator $(\Psi  DO)$ of order $d$ induces operators
$A : H_s ({\cal E}) \rightarrow H_{s'} ({\cal E})$ which are bounded if $s' \le
s-d.$ 
The operator $(\Delta + id)^{-s'/2} A(\Delta + id)^{s/2} :
L_2 ({\cal E}) \rightarrow L_2 ({\cal E})$ has a kernel of
class $C^k$ if $s' < s - d - k - dimM.$ If $A$ is an
elliptic $\Psi D0$ of order $d > 0$ it satisfies $Op(1) - Op(7)$ and (Asy)
(with $m := dimM$ and $\rho = 1/k$) cf [BFKM] section 2. Hence $A$ is a 
$\zeta$-regular operator of sF type.

\spf

\subsection{Complexes}

\n In this subsection we introduce the class of $\zeta$-regular complexes of
sF type and define their algebraic and reduced cohomology. In the case
where a $\zeta$-regular complex of sF type is algebraically acyclic or, more
generally, of determinant class, one can define its torsion. We recall Milnor's
lemma which relates the torsions of a short exact sequence of complexes of 
finite (von Neumann) di\-men\-sion. Further, within our class of complexes, we
discuss the notion of a mapping cone which is a complex induced by a morphism
$f$ between two complexes and show that if $f$ is of trace class and induces
an isomorphism in algebraic cohomology then the mapping cone is algebraically
acyclic. Proposition 1.15
relates, under appropriate conditions, the torsions of the 
map\-ping cones of two morphisms and their composition. The class of 
$\zeta$-regular complexes of sF type has been chosen so that it includes the 
mapping cone of a regularized version of the integration map 
between the deRham complex and the combinatorial complex, which will be 
discussed in section 2. 

\spd

\n A cochain complex ${\cal C} \equiv ({\cal C}_i, d_i)$ of 
${\cal A}$-Hilbert modules,

\[{\cal C}_0 
\stackrel{d_0}{\rightarrow} {\cal C}_1 
\stackrel{d_1}{\rightarrow} \ldots \rightarrow {\cal C}_{N-1} 
\stackrel{d_{N-1}}{\longrightarrow} {\cal C}_N\]

\n consists of a finite sequence of ${\cal A}$-Hilbert modules 
${\cal C}_i (0 \le i \le N)$ and operators 
$d_i :{\cal C}_i \rightarrow {\cal C}_{i + 1},$
satisfying $Op(1) - Op(3),$ with the property that 
${\rm range}(d_i) \subseteq
{\rm domain} (d_{i+1})$ and $d_{i+1} d_i = 0.$

\spd

\n Notice that the null space $Kerd_i$ is always an ${\cal A}$-Hilbert 
submodule.

\spf

\n \begin{definition} For $0 \le i \le N:$ 

\spd

\n algebraic cohomology: $H^i({\cal C}) := 
Kerd_{i}/{\rm range}(d_{i-1})$

\spd

\n reduced cohomology $:\overline{H}^i ({\cal C}) := 
Kerd_{i}/\overline{{\rm range}(d_{i-1})}.$ 
\end{definition}
\spd

\n We point out that the redu\-ced cohomology
$\overline{H}^i ({\cal C})$ is an ${\cal A}$-Hilbert module whereas the 
algebraic cohomology 
$H^i({\cal C})$ is, in general, not an ${\cal A}$-Hilbert module, as
${\rm range}(d_{i-1})$ needs not to be closed. $H^i({\cal C})$ is always an 
${\cal A}$- module.

\spd

\n Given a complex $({\cal C}_i, d_i)$ we can define the adjoint operators 
$d_i^{\ast} : {\cal C}_{i+1} \rightarrow {\cal C}_i$ and, in turn, the
Laplacians 
$\Delta_i,$

\spf

\n {\bf (1.15)} \hsps $
\Delta_i = d^{\ast}_i d_i + d_{i-1} 
d^{\ast}_{i-1}$

\spd

\n as well as the Hodge decomposition 
${\cal C}_i = {\cal H}_i \oplus {\cal C}^+_i \oplus {\cal C}^-_i$ where

\spf

\n {\bf (1.16)} \hsps ${\cal H}_i := Ker(d_i) \cap Ker(d_{i-1}^{\ast}).$

\spf

\n {\bf (1.17)} \hsps ${\cal C}^+_i := \overline{{\rm range} (d_{i-1})}; 
{\cal C}^-_i := 
\overline{{\rm range} (d^{\ast}_i)}.$ 

\spd

\n With respect to this decomposition, the operators $d_i, d^{\ast}_i$ and 
$\Delta_i$ take the form

\spf

\n {\bf (1.18)} \hsps $ d_i = \bl  {ccc}
0 & 0 & 0 \\
0 & 0 & \underline{d}_i \\
0 & 0 & 0 \er ; d^{\ast}_i = 
\bl {ccc}
0 & 0 & 0 \\
0 & 0 & 0 \\
0 & \underline{d}_i^{\ast} & 0 \er ; 
$

\spd

\hspace{1.1cm} $\Delta_i = \bl {ccc}
0 & 0 & 0  \\
0 & \underline{d}_{i-1} \underline{d}^{\ast}_{i-1} & 0 \\
0 & 0 & \underline{d}^\ast_i \underline{d}_i \er .$

\spd

We write $\underline{\Delta}^+_i:= 
\underline{d}_{i-1} \underline{d}^{\ast}_{i-1}$ and 
$\underline{\Delta}^-_i:= \underline{d}^{\ast}_i \underline{d}_i $
\n Note that $\underline{d}_i : {\cal C}^-_i \rightarrow 
{\cal C}^+_{i+1}$ and 
$\underline{d}^{\ast}_i : 
{\cal C}^+_{i+1} \rightarrow {\cal C}^-_i$ are injective operators
with dense image, (hence so are $\underline{\Delta}^+_i$ and  
 $\underline{\Delta}^+_i$,)and ${\cal H}_i$ is isometric to the reduced 
cohomology $\overline{H}^i ({\cal C}).$

\spf

\begin{definition} \hspace{8cm}

\begin{itemize}
\item[(i)] The cochain complex ${\cal C} =
({\cal C}_i, d_i)$ is said to be of sF type if

\item[(CX1)] the operators $\underline{d}_i$ satisfy $Op(1)- Op(6)$ (and hence,
in view of the injectivity of $\underline{d}_i,$ are of sF type);

\item[(CX2)] \hspace{2cm} $dim {\cal H}_i < \infty.$

\item[(ii)] A complex $({\cal C}_i, d_i)$ of sF type is called
\underline{$\zeta$-regular } if, in addition,

\item[(CX3)] \hspace{2cm} the operators $\underline{d}_i$ satisfy $Op(7).$

\item[(iii)] A complex $({\cal C}_i, d_i)$ is said to be of 
finite type (of finite dimension), if the ${\cal A}$-Hilbert
modules ${\cal C}_i$ are of finite type (of finite dimension) and the operators 
$\underline{d}_i$ are bounded. 
\end{itemize}
\end{definition}

\spf

\n {\bf Remark 1} The conditions (CX1) and (CX2) are equivalent to

\begin{itemize}
\item[$(CX_{12})$] \hspace{2cm} {\em the Laplacians 
$\Delta_i$ are of sF type} 
\end{itemize}

\n and the three properties (CX1)-(CX3) are equivalent to

\begin{itemize}
\item[$(CX_{123})$] \hspace{2cm} {\em the Laplacians $\Delta_i$ are 
$\zeta$-regular
operators of sF type.}
\end{itemize}

\spf

\n {\bf Remark 2} Assume that ${\cal H}_{i_0} = 0$ for some $i_0.$ Using the 
closed graph theorem one verifies that the following statements are 
equivalent: 

\[ \mbox{(a)} \  H^{i_0} 
({\cal C}) = 0; \ \mbox{(b)} \ \underline{d}_{i_0} \ \mbox{has a bounded 
inverse.} \]

\spd

\n Given a complex $({\cal C}_i, d_i)$  of sF type, one can define the spectral
distribution functions $F_i (\lambda)$ and $N_i(\lambda),$

\spf

\n {\bf (1.19)} \hsps $F_i (\lambda) := F_{\underline{d}_i} (\lambda) 
\Big( = F_{|\underline{d}_i|} (\lambda)\Big); N_i (\lambda) := 
F_{\Delta_i} (\lambda).$

\spd

\n Then

\spf

\n {\bf (1.20)} \hsps $N_i(\lambda) = F_i(\lambda) + F_{i-1} (\lambda).$

\spd

\n A complex $({\cal C}_i, d_i)$ of sF type is said to be \underline{algebraically
acyclic} if $H^i ({\cal C}) = 0$ for $0 \le i \le N.$

\spd

\n Notice that for an acyclic complex of sF type, the spectrum of 
$\Delta_i$ does not contain $0$. In fact, in view of Remark 2 after
Definition 1.9, there exists $\varepsilon > 0$ so that $N_i(\lambda) = 0$ for
$0 \le \lambda < \varepsilon$ and the same is true for the spectral
distribution function of $\underline{\Delta}^+_i = \underline{d}_{i-1}
\underline{d}_{i-1}^{\ast}$ and $\underline{\Delta}_i^- = 
\underline{d}_i^{\ast} \underline{d}_i.$ This allows to define the torsion
$T({\cal C})$ of an algebraically acyclic complex of sF type,

\spf

\n {\bf (1.21)} $
\begin{array}{lll}
log T({\cal C}) &:= & \frac{1}{2} \sum_q (-1)^{q+1} q \ log \ det 
\Delta_q \\
&& \Big( = \frac{1}{2} \sum_q (-1)^q log \ det \underline{\Delta}^-_q =
\frac{1}{2} \sum_q (-1)^{q+1} log \ det \underline{\Delta}^+_q \Big).
\end{array}$

\spd

\n More generally, the torsion $T({\cal C})$ can be defined if 
$({\cal C},d)$ is of 
determinant class, i.e.

\[\int^1_{0^+} log \lambda dN_i (\lambda) > - \infty \ \ (0 \le i \le N).\]

\spf

\begin{definition} (i) A morphism $f : {\cal C}^1 \rightarrow 
{\cal C}^2$ between the
complexes $({\cal C}^1_i, d_{1,i})$ and $({\cal C}^2_i, d_{2,i})$ consists
of a collection of bounded operators $f_i : {\cal C}^1_i \rightarrow
{\cal C}^2_i$ so that

\begin{itemize}
\item $f_i ({\rm domain} (d_{1,i})) 
\subseteq {\rm domain}(d_{2,i}); \ d_{2,i} f_i =
f_{i+1} d_{1,i} \ (on \ {\rm domain} (d_{1,i})).$
\item[(ii)] $f$ is said to be of trace class if, in addition,
\item the operators $f_i, d_{2,i}^{\ast} f_{i+1}$ and $f_i 
d^{\ast}_{1,i}$ are bounded operators of trace class.
\end{itemize}
\end{definition}

\spd

\n A morphism $f: {\cal C}^1 \rightarrow {\cal C}^2$ induces ${\cal A}$-linear
maps in {\em algebraic} co\-ho\-mo\-lo\-gy, $H(f); : H^i({\cal C}^1) \rightarrow 
H^i({\cal C}^2)$ and bounded ${\cal A}$-linear maps on the {\em reduced} 
cohomology, $\overline{H} (f)_i : \overline{H}^i ({\cal C}^1) \rightarrow 
\overline{H}^i ({\cal C}^2).$ Further, with respect to the Hodge decomposition,
$f_i$ takes the form,

\[ f_i = \bl {lll}
f_{i,11} & 0 & f_{i,13} \\
f_{i,21} & f_{i,22} & f_{i,23} \\
0 & 0 & f_{i,33} \er \]

\n where $f_{i,11} = \overline{H} (f)_i.$

\spd

\n For the convenience of the reader, we recall Milnor's lemma which was 
extended for complexes of ${\cal A}$-Hilbert modules in [BFK2]. Assume that

\[0 \rightarrow {\cal C}^1 \stackrel{f}{\rightarrow} {\cal C}^2 
\stackrel{g}{\rightarrow} {\cal C}^3 \rightarrow 0\]

\n is a short exact sequence of complexes of finite dimension. It induces a
long weakly exact sequence ${\cal H}$ in reduced cohomology 

\[ \ldots \rightarrow \overline{H}^i ({\cal C}^1) 
\stackrel{\overline{H}(f_i)}{\longrightarrow} 
\overline{H}^i ({\cal C}^2)  \stackrel{\overline{H}(g_i)}{\longrightarrow} 
\overline{H}^i ({\cal C}^3) \stackrel{\overline{H} (\delta_i)}{\longrightarrow} 
\overline{H}^{i+1} ({\cal C}^1) \rightarrow \ldots \]

\spf

\begin{prop} (cf [BFK2, Theorem 1.5]) If three out of the four cochain 
complexes ${\cal C}^1, {\cal C}^2, {\cal C}^3, {\cal H}$ are of determinant
class, then so is the fourth and one has the following equality

\[ \begin{array}{lll}
log T({\cal C}^2) & = & log T({\cal C}^1) + logT({\cal C}^3) + log T ({\cal H}) - \\
& - & \sum_i (-1)^i logT (0 \rightarrow {\cal C}^1_i \rightarrow {\cal C}^2_i
\rightarrow {\cal C}^3_i \rightarrow 0). \end{array} \]
\end{prop}

\spf

\n {\bf Remark} In [BFK2], Lemma B2.8 was only stated for complexes of
finite type. By the same proof, the result remains true for complexes of 
finite di\-men\-sion.

\spd 

\n Given a $\zeta$-regular complex of sF type, $0 \rightarrow {\cal C}_0 
\stackrel{d_0}{\rightarrow} \ldots \stackrel{d_{2N}}{\rightarrow} {\cal C}_{2N+1}
\rightarrow 0,$ we can consider its dual. If ${\cal C}^{\sharp}_j := 
{\cal C}_{2N + 1 - j},$ and $d^{\sharp}_j := d_{2N-j}^{\ast}$ (adjoint of 
$d_{2N-j}$), then

\[0 \rightarrow
{\cal C}^{\sharp}_0 \stackrel{d^{\sharp}_0}{\rightarrow} \ldots
\stackrel{d^{\sharp}_{2N}}{\rightarrow} {\cal C}^{\sharp}_{2N+1} \rightarrow
0\]

\n is a $\zeta$-regular complex of sF type. Notice that 
$(\underline{d}^{\sharp}_j)^{\ast}
\underline{d}^{\sharp}_i = \underline{d}_{2N-j} \ \underline{d}_{2N-j}^{\ast}.$
Thus in case ${\cal C}$ is of determinant class, so is ${\cal C}^{\sharp}$ and
one obtains 

\spf

\n {\bf (1.22)} \hsps $T({\cal C}^{\sharp}) = T({\cal C}).$

\spd

\n We end this subsection with a discussion of the {\em mapping cone} 
${\cal C}(f) = ({\cal C}(f)_i, 
\linebreak
d(f)_i)$ associated to a morphism $f : 
{\cal C}^1 \rightarrow {\cal C}^2$ between $\zeta$-regular complexes 
$({\cal C}^1, d_1)$ and $({\cal C}^2, d_2)$ of sF type. This is a complex 
given by 

\begin{itemize}
\item[(MC1)] \ ${\cal C}(f)_i := {\cal C}^2_{i-1}  \oplus
{\cal C}^1_i;$
\item[(MC2)] \ $d(f)_i = \bl {cc}
-d_{2,i-1} & f_i \\
0 & d_{1,i} \er .$
\end{itemize}

\n The Laplacians $\Delta (f)_i$ of ${\cal C}(f)$ can be computed,

\spf

\n {\bf (1.23)} \ $\Delta (f)_i = \bl {cc}
\Delta_{2,i-1} + f_{i-1} f_{i-1}^{\ast} & 
-d^{\ast}_{2, i-1} \cdot f_i + f_{i-1} \cdot d^{\ast}_{1, i-1} \\
-f^{\ast}_i \cdot d_{2,i-1} + d_{1,i-1} \cdot f^{\ast}_{i-1} &
\Delta_{1,i} + f^{\ast}_i f_i \er
$

\spd

\n The following observation,  will  be used in order 
to in\-tro\-du\-ce the relative $L_2$-torsion without making the assumption that 
the manifold is of determinant class:

\spf

\begin{lemma} (i) Let $f : {\cal C}^1 \rightarrow {\cal C}^2$ be a morphism
of trace class between com\-ple\-xes ${\cal C}^1$ and ${\cal C}^2$ of sF type. 
Then the mapping cone ${\cal C}(f)$ is of sF type.

\begin{itemize}
\item[(ii)] If, in addition, ${\cal C}^1$ and ${\cal C}^2$ are $\zeta$-regular,
then ${\cal C}(f)$ is $\zeta$-regular.
\item[(iii)] If $\varphi : {\cal C}^1 \rightarrow {\cal C}^2$ is a morphism of
trace class between $\zeta$-regular complexes of sF type which induces an 
isomorphism in algebraic cohomology, then the mapping cone ${\cal C}(f)$ is
algebraically acyclic and thus has a well defined torsion $T({\cal C}(f))$
(cf (1.21)).
\end{itemize}
\end{lemma}

\spf

\n {\bf Proof} (i) follows from Proposition 1.5 (A) by choosing
$\varphi := \bl {cc}
\Delta_{2,i-1} & 0 \\
0 & \Delta_{1,i} \er$ and

\[u_i := \Delta (f)_i - \varphi_i = \bl {cc}
f_{i-1} f_{i-1}^{\ast} & -d^{\ast}_{2,i-1} f_i + f_{i-1} d^{\ast}_{1,i-1} \\
-f_i^{\ast} d_{2,i-1} + d_{1,i-1}f_{i-1}^{\ast} & f_i^{\ast} f_i \er \]

\n By Remark 1 after Definition 1.8, the $\varphi_i$ are operators of 
sF type and, by assumption, the $u_i$ are bounded operators of 
trace class. By Proposition 1.5, $\Delta (f)_i$ is then an operator
of sF type. Apply the same Remark 1 once more to conclude that (CX1) holds.
Property (CX2) of Definition 1.8 is easily verified. 

\spd

\n (ii) follows from Proposition 1.5 (B) and Remark 1 after Definition 1.8.

\spd

\n (iii) In view of (1.21) and the statements (i) and (ii) it remains to 
verify that ${\cal C} (f)$ is algebraically acyclic, i.e. that 
${\rm range}(d(f)_{i-1}) = Ker(d(f)_i).$ Let $(v_{i-1}, u_i) \in Ker(d(f)_i) 
\subseteq {\cal C}^2_{i-1} \oplus {\cal C}^1_i$. Then

\[-d_{2,i-1} v_{i-1} + f_i u_i = 0 \ ;  \ \ d_{1,i} u_i = 0.\]

\spd

\n In particular, $f_i u_i \in {\rm range} (d_{2,i-1})$ and, as $f$ induces an 
isomorphism in algebraic cohomology, $u_i \in {\rm range}(d_{1, i -
1}),$ i.e. there
exists $u_{i-1} \in {\cal C}^1_{i-1}$ with $u_i = d_{1, i-1} 
\linebreak
u_{i-1}.$ As
$f_i d_{1, i-1} u_{i-1} = d_{2, i-1} f_{i-1} u_{i-1},$ we conclude that
$-v_{i-1} + f_{i-1} u_{i-1} \in Kerd_{2,i-1}.$ Using once again the assumption
that $f$ induces an isomorphism in algebraic cohomology, one sees that there
exists $v_{i-2} \in {\cal C}^2_{i-2}$ and $\tilde{u}_{i-1} \in Kerd_{1,i-1}$
with $-v_{i-1} + f_{i-1} u_{i-1} = f_{i-1} \tilde{u}_{i-1} + d_{2, i-2} v_{i-2}.$

\spd

\n Thus

\[v_{i-1} = -d_{2,i-2} v_{i-2} + f_{i-1} (u_{i-1} - \tilde{u}_{i-1})\]

\n and, as $\tilde{u}_{i-1} \in Ker d_{1,i-1}$

\[u_i = d_{1, i-1} u_{i-1} = d_{1,i-1} (u_{i-1} - \tilde{u}_{i-1})\]

\n i.e. we have shown that $(v_{i-1}, v_i) \in {\rm range} d(f)_{i-1},$ and
therefore that $Kerd(f)_i
\linebreak
= {\rm range} d(f)_{i-1}.$ \ \ \ \carre

\spd

\n If $f : ({\cal C}, d_1) \rightarrow ({\cal C}, d_2)$ is an isomorphism
of complexes of finite dimension we can use Lemma 1.11 together with Milnor's
lemma to compute the torsion of the mapping cone ${\cal C}(f)$. Recall that
the {\em suspension} $\Sigma{\cal C} = (\Sigma{\cal C}_{i,\Sigma}d_i)$ of a
complex ${\cal C} = ({\cal C}_i, d_i)$ is the complex given by

\[\begin{array}{l}
\Sigma{\cal C}_i \equiv (\Sigma{\cal C})_i := {\cal C}_{i-1} (i \ge 1); 
\ \Sigma{\cal C}_0 := 0; \\
_{\Sigma}d_i : =  -d_{i-1} \ (i \ge 1); \ d_0 = 0.
\end{array}\]
 
\spf

\begin{prop} Assume that ${\cal C}^i \equiv ({\cal C}, d_j)(j=1,2)$ are
$\zeta$-regular complexes of sF type and of finite dimension. If 
$f: {\cal C}^1 \rightarrow {\cal C}^2,$ is an isomorphism, then the mapping
cone ${\cal C}(f)$ is algebraically acyclic and

\[logT({\cal C}(f)) = \sum(-1)^j log \ vol(f_j)\]
\end{prop}

\spf

\n {\bf Proof} By Lemma 1.11, ${\cal C}(f)$ is an algebraically 
$\zeta$-regular complex of sF type.

\spd

\n First we treat the case where, in addition, ${\cal C}^1$ (and then 
${\cal C}^2$ as well) is al\-ge\-brai\-cal\-ly acyclic. Consider the following short
exact sequence of cochain complexes

\[0 \rightarrow \Sigma({\cal C}, d_2) \stackrel{J}{\rightarrow} {\cal C}(f)
\stackrel{P}{\rightarrow} ({\cal C}, d_1) \rightarrow 0\]

\n where $J$ (respectively P) is the canonical inclusion (canonical 
projection).

\spd

\n By Proposition 1.10, the corresponding long weakly exact sequence 
${\cal H}$ in reduced cohomology is of determinant class and

\[log T({\cal C}(f)) = log T ( \Sigma{\cal C}^2) + log T ({\cal C}^1) + logT
({\cal H}).\]

\n Notice that ${\cal H}$ is given by

\[\ldots \rightarrow \overline{H}^i (\Sigma{\cal C}) \rightarrow 
\overline{H}^i ({\cal C}(f)) \rightarrow \overline{H}^i ({\cal C}) 
\stackrel{\overline{H}(\delta_i)}{\longrightarrow} \overline{H}^{i+1} 
(\Sigma{\cal C}) \rightarrow \ldots\]

\spd

\n where the connecting homomorphism $\overline{H}(\delta_i)$ is given by 
the restriction $\overline{H} (f_i)$ of $f_i$ to $\overline{H}^{i+1} (\Sigma
{\cal C}) = \overline{H}^i ({\cal C})$ into $\overline{H}^i({\cal C}).$
Therefore 

\[log T({\cal H}) = \sum(-1)^i log \ vol(\overline{H} (f_i)).\]

\n By Lemma 1.13 below, $log T(\Sigma{\cal C}^2) = - log T({\cal C}^2).$ As
$f_{i+1} d_{1,i} = d_{2,i} f_i$ we con\-clu\-de that $\underline{d}_{1,i} =
(f^+_{i+1})^{-1} \underline{d}_{2,i} f^-_i$ where $f^{\pm}_i : 
{\cal C}^{1, \pm}_i \rightarrow {\cal C}^{2, \pm}_i$ are the restrictions of 
$f_i$ to ${\cal C}^{1, \pm}_i$ and are isomorphisms as well.

\spd

\n Thus

\[\begin{array}{lll}
log T({\cal C}^1) & = & \frac{1}{2} \sum (-1)^i log \ det 
(\underline{d}_{1,i}^{\ast} \underline{d}_{1,i}) \\
& = & \frac{1}{2} \sum (-1)^i log \ det ((f^-_i)^{\ast} 
\underline{d}_{2,i}^{\ast} ((f^+_{i+i})^{-1})^{\ast}(f^+_{i+1})^{-1} 
\underline{d}_{2,i} f_i^-)\\
& = & \frac{1}{2} \sum (-1)^i log \ det (f^+_{i+1} f^{+ \ast}_{i+1})^{-1} \\
& + & \frac{1}{2} \sum (-1)^i log \ det (\underline{d}^{\ast}_{2,i}
\underline{d}_{2,i}) \\
& + & \frac{1}{2} \sum (-1)^i log \ det (f_i^- f_i^{-{\ast}}) \\
& = & log T ({\cal C}^2) + \sum_i (-1)^i (log \ vol (f_i^+) + log \ vol (f_i^-)). 
\end{array} \]

\n Combining the equality above yields

\[ \begin{array}{lll}
log T ({\cal C}(f)) & = & \sum(-1)^i (log \ vol (f_i^+) + log \ vol (f_i^-)) \\
& + & \sum (-1)^i log \ vol (\overline{H} (f_i)) \\
& = & \sum (-1)^i log \ vol (f_i).
\end{array} \]

\spd

\n To prove the result in general we consider a deformation $({\cal C}, d_1(
\varepsilon))$ of the complex $({\cal C}, d_1),$ depending smoothly on the
parameter $\varepsilon,$ so that, for $\varepsilon \neq 0, ({\cal C}, d_1 
(\varepsilon))$ is algebraically acyclic. The operator $d_1 (\varepsilon)$
is constructed as follows: with respect to the Hodge de\-com\-po\-si\-tion 
${\cal C}_i
= {\cal H}_i \oplus {\cal C}^+_i \oplus {\cal C}^-_i, d_{1,i}$ takes the form
$d_{1,i} = \bl {lll} 0 & 0 & 0 \\ 0 & 0 & \underline{d}_{1,i} \\
0 & 0 & 0 \er$ where $\underline{d}_{1,i} : {\cal C}^-_{1,i} 
\rightarrow {\cal C}^+_{1, i+1}.$ Consider the polar de\-com\-po\-si\-tion of
$\underline{d}_{1,i}, \underline{d}_{1,i} = \xi_i, \eta_i$ where 
$\xi_i := \underline{d}_{1,i} (\underline{d}^{\ast}_{1,i} 
\underline{d}_{1,i})^{-1/2}
: {\cal C}^-_i \rightarrow {\cal C}^+_{i+1}$ is an isometry and
$\eta_i$ is given by $\eta_i := (\underline{d}^{\ast}_{1,i} 
\underline{d}_{1,i})^{1/2}$.
The operator $\eta_i$ admits a spectral decomposition; $\eta_i = 
\int^{\infty}_{0^+}
\lambda dP_i (\lambda).$ Define $\eta_i (\varepsilon) := \int^{\infty}_{0^+}
(\lambda + \varepsilon) dP_i (\lambda)$ and let $\underline{d}_{1,i} 
(\varepsilon)
:= \xi_i \eta_i (\varepsilon).$ Then $\underline{d}_{1,i} (\varepsilon)$ is an 
isomorphism and thus of determinant class. Now consider ${\cal C}^1_{\varepsilon}
\equiv ({\cal C}, \underline{d}_1, (\varepsilon))$ where the operator $d_{1,i}
(\varepsilon),$ with respect to the Hodge decomposition of $({\cal C}, d_1),$
is given by $d_{1,i} (\varepsilon) = \bl {lll} 0 & 0 & 0 \\ 0 & 0 & 
\underline{d}_{1,i} (\varepsilon) \\
0 & 0 & 0 \er.$ Define $({\cal C}, d_2 (\varepsilon))$ where $\underline{d}_{2,i}
(\varepsilon) : = f^+_{i+1} \ \underline{d}_{1, \varepsilon} (f_i^-)^{-1}.$ 
Then,
for any $\varepsilon, \ {f:} ({\cal C}, d_1 (\varepsilon)) \rightarrow ({\cal C},
d_2 (\varepsilon))$ is a morphism and, by the argument above, the torsion 
$T({\cal C}_{\varepsilon} (f))$ of the mapping cone ${\cal C}_{\varepsilon} (f)$
induced by $f: ({\cal C}, d_1 (\varepsilon)) \rightarrow ({\cal C}, d_2 
(\varepsilon))$ is given by, for $\varepsilon \neq 0,$

\[log T({\cal C}_{\varepsilon} (f)) = \sum (-1)^i log \ vol (f_i)\]

\n which is independent of $\varepsilon.$ As ${\cal C}_{\varepsilon} (f)$ is
algebraically acyclic for all values of $\varepsilon$ and $log T 
({\cal C}_{\varepsilon} (f))$ depends continuously on $\varepsilon,$ we 
conclude that

\[log T ({\cal C} (f)) = \sum_j (-1)^j log \ vol (f_j). \ \ \ \mbox{\carre} \]

\spf

\begin{lemma} (i) Let ${\cal C}$ be a $\zeta$-regular complex of sF type. Then
the mapping cone ${\cal C}(id)$ is an algebraically acyclic $\zeta$-regular 
complex of sF type and satisfies $log T ({\cal C} (id)) = 0$.

\begin{itemize}
\item[(ii)] If ${\cal C}$ is an algebraically acyclic, $\zeta$-regular complex
of sF type, then so is $\Sigma {\cal C}$ and satisfies

\[log T (\Sigma {\cal C}) = - log T ({\cal C}).\]
\end{itemize}
\end{lemma}

\spd

\n {\bf Proof} (i) By Lemma 1.11, ${\cal C} (id)$ is an algebraically acyclic
$\zeta$-regular complex of sF type. Use Lemma 2.4 to conclude that 

\[log T ({\cal C}(id)) = \frac{1}{2} \sum (-1)^{q+1} log \ det (\Delta_q 
+ id) = 0 .\]

\n (ii) It follows from Definition 1.9 and the assumptions that $\Sigma {\cal C}$
is an al\-ge\-brai\-cal\-ly acyclic, $\zeta$-regular complex of sF type.
Its torsion is therefore well defined, and by (1.21)

\[ \begin{array}{lll}
log T (\Sigma {\cal C}) & = & \frac{1}{2} \sum (-1)^{q+1} log \ det 
(\underline{\Delta (\Sigma {\cal C})}^-_q) \\
&= & - \frac{1}{2} \sum (-1)^{q+1} log \ det (\underline{\Delta}^-_q) = - 
log T({\cal C}). \ \ \ \mbox{\carre}
\end{array} \]

\spd

\n The following Lemma 1.14 is due to Carey, Mathai and Mishchenko [CMM]. For
the convenience of the reader we include its statement. As the proof in the
preliminary preprint [CMM]\footnote{The first author thanks Mishchenko for 
kindly making the preliminary preprint available to him.} is somewhat incomplete,
we refer the reader to Appendix A for a detailed proof. Let

\spf

\n {\bf (1.24)} \hsps $0 \rightarrow {\cal C}^1 \stackrel{I}{\rightarrow}
{\cal C} \stackrel{P}{\rightarrow} {\cal C}^2 \rightarrow 0 $

\spd

\n be an exact sequence of cochain complexes of sF type where the complex
${\cal C} = ({\cal C}_i, d_i)$ is given by ${\cal C}_i = {\cal C}^1_i \oplus
{\cal C}^2_i$ and 

\spf

\n {\bf (1.25)} \hsps $d_i = \bl {ll} d_{1,i} & f_i \\
0 & d_{2,i} \er $

\spd

\n with $f_i : {\cal C}^2_i \rightarrow {\cal C}^1_{i+1}$ satisfying

\spf

\n {\bf (1.26)} \hsps $f_i ({\rm domain} (d_{2,i})) \subseteq {\rm domain}
(d_{1, i+1})$ \ and

\spf

\n {\bf (1.27)} \hsps $f_{i+1} d_{2,i} + d_{1,i+1} f_i = 0 \ \ ({\rm on  \ domain}
(d_{2,i})) $ 

\spd

\n so that $d_{i+1} d_i = 0.$ The morphism I [resp. P] in (1.24) denotes the
canonical inclusion [resp. canonical projection].

\spf

\begin{lemma} ([CMM]) Assume ${\cal C}^1$ and ${\cal C}^2$ are $\zeta$-regular
complexes of sF type which are algebraically acyclic and $f_i : {\cal C}^2_i 
\rightarrow {\cal C}^1_{i+1}$ are bounded maps of trace class so that (1.26)
and (1.27) are satisfied and $d_{1,i}^{\ast} f_i$ as well as $f_i d^{\ast}_{2,i}$ are
bounded operators of trace class. Then the complex ${\cal C},$ given by (1.24) -
(1.25), is an  algebraically acyclic, $\zeta$-regular complex of sF type and

\spf

\n {\bf (1.28)} \hsps $log T ({\cal C}) = log T({\cal C}^1) + 
log T ({\cal C}^2).$
\end{lemma}

\spd

\n  {\bf Proof:} cf [CMM] or Appendix A.

\spf

\begin{prop} Suppose that ${\cal C}^1, {\cal C}^2$ and ${\cal C}^3$ are 
$\zeta$-regular complexes of sF type and $f_1 : {\cal C}^1 \rightarrow 
{\cal C}^2$
and $f_2 : {\cal C}^2 \rightarrow {\cal C}^3$ are morphisms of trace class
which induce isomorphism in algebraic cohomology. Then:

\spd

\begin{itemize}
\item[(i)] $f := f_2 \circ f_1$ is a morphism of trace class;

\item[(ii)] the mapping cones ${\cal C} (f_1), {\cal C}(f_2)$ and ${\cal C}(f)$ are
$\zeta$-regular complexes of sF type;

\item[(iii)] ${\cal C} (f_1), {\cal C}(f_2)$ and ${\cal C} (f)$ are algebraically
acyclic (and thus their torsions are well defined).

\item[(iv)] If, in addition, either ${\cal C}^3$ (case 1) or ${\cal C}^1$ 
(case 2) is of finite type, then $log T ({\cal C}(f)) = log T ({\cal C} (f_1))
+ log T ({\cal C} (f_2)).$ 
\end{itemize}
\end{prop}

\spf

\n {\bf Proof} \ The proof is an application of Lemma 1.13 and Lemma 1.14.

\spf

\n {\bf Case 1} \ (${\cal C}^3$ is of finite type): Consider the diagram

\spf

\n {\bf (1.29)} \hsps $\begin{array}{ccccccc}
&& \Sigma {\cal C} (id_3) &&&& \\
&&&&&& \\
&& \downarrow I_1 &&&& \\
&&&&&& \\
\Sigma {\cal C} (-f_2) & \stackrel{I_2}{\rightarrow} & {\cal C}(h) & 
\stackrel{P_2}{\rightarrow} & {\cal C}(f) & \rightarrow & 0 \\
&&&&&& \\
&& \downarrow P_1 &&&& \\
&&&&&& \\
&& {\cal C}(f_1) &&&& \\
&&&&&& \\
&& \downarrow &&&& \\
&&&&&& \\
&& 0 &&&& \end{array} $

\n where ${\cal C} (f_1), {\cal C}(-f_2), {\cal C}(f)$ and ${\cal C}(id_3)$ 
are mapping cones, hence ${\cal C}(f_1); = {\cal C}^2_{i-1} \oplus {\cal C}^1_i,
{\cal C} (-f_2)_i = {\cal C}^3_{i-1} \oplus {\cal C}^2_i, {\cal C}(f)_i =
{\cal C}^3_{i-1} \oplus {\cal C}^1_i$ and ${\cal C}(id_3)_i = {\cal C}^3_{i-1}
\oplus {\cal C}^3_i.$ The morphism $h : {\cal C}(f_1) \rightarrow {\cal C}
(id_3),$ is given by the operators $h_i : {\cal C}^2_{i-1} \oplus {\cal C}^1_i
\rightarrow {\cal C}^3_{i-1} \oplus {\cal C}^3_i,$

\spd

\n {\bf (1.30)} \hsps $h_i := \bl {cc}
f_{2,i-1} & 0 \\
0 & f_i \er \ , $

\spd

\n One verifies that $h_i$ are bounded maps of trace class. Denote by ${\cal C}
(h)$ the mapping cone of $h,$ hence ${\cal C} (h)_i = ({\cal C}^3_{i-2} 
\oplus {\cal C}^3_{i-1}) \oplus ({\cal C}^2_{i-1} \oplus {\cal C}^1_i).$ 
Further $I_1,$ respectively $I_2$, are the canonical inclusions,

\[I_1 : {\cal C}^3_{i-2} \oplus {\cal C}^3_{i-1} \rightarrow ({\cal C}^3_{i-2}
\oplus {\cal C}^3_{i-1}) \oplus ({\cal C}^2_{i-1} \oplus {\cal C}^1_i)\]

\n and

\[I_2 : {\cal C}^3_{i-2} \oplus {\cal C}^2_{i-1} \rightarrow ({\cal C}^3_{i-2}
\oplus {\cal C}^3_{i-1}) \oplus ({\cal C}^2_{i-1} \oplus {\cal C}^1_i)\]

\n whereas $P_1$ and $P_2$ denote the canonical projections on the complement of the
images of $I_1$ resp. $I_2$. One verifies (cf (1.29)) that

\spf

\n {\bf (1.31)} \hsps $0 \rightarrow \Sigma{\cal C} (id_3) \stackrel{I_1}{\rightarrow}
{\cal C}(h) \stackrel{P_1}{\rightarrow} {\cal C} (f_1) \rightarrow 0$

\spd

\n and

\spf

\n {\bf (1.32)} \hsps $0 \rightarrow \Sigma {\cal C} (-f_2) 
\stackrel{I_2}{\rightarrow} {\cal C}(h) \stackrel{P_2}{\rightarrow}
 {\cal C}(f) \rightarrow 0$

\spd

\n are short exact sequences of cochain complexes. First notice that by 
Proposition 1.3, $f= f_2 \cdot f_1$ is a morphism of trace class which proves
(i).

\spd

\n By Lemma 1.11, ${\cal C}(f_1), {\cal C}(f_2)$ and ${\cal C}(f)$ are 
al\-ge\-brai\-cal\-ly acyclic, $\zeta$-regular com\-ple\-xes of sF type 
and, by Lemma 1.13, 
${\cal C} (id_3)$ and $\Sigma{\cal C}(f_i), \Sigma{\cal C}(f_2), 
\Sigma{\cal C}(id_3)$ are algebraically acyclic, $\zeta$-regular complex of sF
type. We would like to apply Lemma 1.14 to (1.31) and (1.32). For this 
purpose,
observe that, given $h_i : {\cal C}(f_1)_i = {\cal C}^2_{i-1} \oplus {\cal C}_i
\rightarrow {\cal C}(id_3) = {\cal C}^3_{i-1} \oplus {\cal C}^3_i, d(h)_i$
admits the following de\-com\-po\-si\-tion

\spf

\n {\bf (1.33)} \ $
\begin{array}{lll}
d(h)_i & =  & \bl {c|c}
d(\Sigma {\cal C}(id_3))_i & h_i \\
\cline{1-2}
0 & d({\cal C}(f_1)_i \er,  \\
&&\\
& = & \bl{cc|cc}
d_{3,i-2} & -id_{3,i-1} & f_{2,i-1} & 0 \\
0 & -d_{3,i-1} & 0 & f_i \\
\cline{1-4}
0 & 0 & d_{2,i-1} & f_{1,i} \\
0 & 0 & 0 & d_{1,i} \er
\end{array}
$

\spd

\n Using this decomposition, all assumptions in Lemma 1.15 which have not
 yet been proved can be verified in a straight forward way.

\spd

\n Therefore, we can apply Lemma 1.14 to (1.31) to conclude that ${\cal C}(h)$
is an algebraically acyclic, $\zeta$-regular complex of sF type and 

\spf

\n {\bf (1.34)} \hsps $
\begin{array}{lll}
log T({\cal C}(h)) & = & logT(\Sigma{\cal C} (id_3)) + logT({\cal C}(f_1)) \\
&=& log T ({\cal C}(f_1))
\end{array}$

\spd

\n where we used that, by Lemma 1.13,

\[log T(\Sigma{\cal C}(id_3)) = - log T({\cal C}(id_3)) = 0\]

\n In order to apply Lemma 1.14 to (1.32), notice that

\spf

\n {\bf (1.35)} \  $
\begin{array}{lll}
\bl {c|c}
d(\Sigma{\cal C}(-f_2))_i & \tilde{h}_i \\
\cline{1-2}
0 & d({\cal C}(f))_i \er & = &   \\
&&\\
\bl{cc|cc}
d_{3,i-2} & +f_{2,i-1} & \tilde{h}_i &  \\
0 & -d_{2,i-1} & &  \\
\cline{1-4}
0 & 0 & -d_{3,i-1} & f_{i} \\
0 & 0 & 0 & d_{1,i} \er & = & d(h)_i
\end{array}
$

\spd

\n where, in view of (1.33), $\tilde{h}_i : {\cal C}(f)_i = {\cal C}^3_{i-1} 
\oplus {\cal C}^1_i \rightarrow {\cal C}^3_{i-1} \oplus {\cal C}^2_i = \Sigma
{\cal C}(-f_2)_{i+1}$ given by $\tilde{h}_i = \bl {cc} 
-id_{3,i-1} & 0 \\
0 & f_{1,i} \er.$ As ${\cal C}^3$ is of finite type, $id_{3,i} : {\cal C}^3_i
\rightarrow {\cal C}^3_i$ and therefore $\tilde{h}_i$ are bounded maps of trace
class. Again one verifies that all other 
assumption in Lemma 1.14 are satisfied. Thus we can apply Lemma 1.14 to (1.32),
to conclude that

\spf

\n {\bf (1.36)} \ $\begin{array}{lll}
log T ({\cal C}(h)) & = & logT (\Sigma {\cal C}(-f_2)) + logT({\cal C}(f)) \\
& = & - logT ({\cal C}(-f_2)) + log T({\cal C}(f))
\end{array}$

\spd

\n where for the last identity we again used Lemma 1.13. Notice that the
morphism $\Phi : {\cal C}(f_2) \rightarrow {\cal C}(-f_2),$ given by $\Phi_i
= \bl {cc}
Id & 0 \\ 0 & -Id \er : {\cal C}^3_{i-1} \oplus {\cal C}^2_i \rightarrow
{\cal C}^3_{i-1} \oplus {\cal C}^2_i,$ is an isometry and therefore $log T 
({\cal C}(-f_2)) = log T({\cal C}(f_2)).$
Combining this with (1.34) and (1.36) we conclude that statement (iv) in case 1
is proved.

\spf

\n {\bf Case 2} \ (${\cal C}^1$ of finite type): In view of (1.22) it suffices to
consider the dual complexes, to which we can apply case 1. \ \ \ \carre

\spf

\begin{prop} Suppose that ${\cal C}^1, {\cal C}^2$ and ${\cal C}^3$ are 
$\zeta$-regular complexes of sF type and $f_1 : {\cal C}^1 \rightarrow {\cal C}^2,
f_2 : {\cal C}^2 \rightarrow {\cal C}^3$ are morphisms which induce isomorphisms
in algebraic cohomology and denote by $f$ the composition $f:= f_2 \circ f_1$. 
Then the following statements hold:

\begin{itemize}
\item[(i)] If $f_1$ is an isometry and $f_2$ is of trace class, then the 
mapping cones ${\cal C}(f_2)$ and ${\cal C}(f)$ are algebraically acyclic
$\zeta$-regular complexes of sF type. Moreover

\[log T ({\cal C}(f_2)) = log T ({\cal C}(f)).\]

\item[(ii)] If $f_2$ is an isometry and $f_1$ is of trace class, then 
${\cal C} (f_1)$ and ${\cal C}(f)$ are algebraically acyclic $\zeta$-regular
complexes of sF type. Moreover
\end{itemize}

\[log T ({\cal C}(f_1)) = log T ({\cal C}(f)).\]
\end{prop}

\spf

\n {\bf Proof} The assumption imply that in case (i), ${\cal C}(f_2)$ and 
${\cal C}(f)$ are isometric whereas in case (ii), ${\cal C}(f_1)$ and
${\cal C}(f)$ are isometric. From this observation and Lemma 1.11 the claimed
results follow. \ \ \ \carre

\section{Relative torsion and its Witten deformation}

\n In subsection 2.1 we introduce the relative torsion (cf [CMM]). As already
mentioned in the introduction, the relative torsion cannot be defined as the
torsion of the mapping cone associated to the integration map, denoted by Int, as 
Int cannot be extended to a closed morphism on $L_2 (\Lambda(M;{\cal E})).$ 
To circumvent this difficulty  we consider
$s > \frac{n}{2} + 1,$ and the composition

\[L_2 (\Lambda(M; {\cal E})) \equiv H_0(\Lambda(M; {\cal E})) 
\stackrel{(\Delta + Id)^{-s/2}}{\longrightarrow} H_s (\Lambda (M; {\cal E})) 
\stackrel{Int_s}{\rightarrow} {\cal C}\]

\n where $H_s (\Lambda^k (M; {\cal E}))$ denotes the completion of the space
$\Omega^k (M; {\cal E})$ of smooth $k$-forms with respect to the $s$-Sobolev 
norm,
$\Delta_k$ denotes the $k$-Laplacian, $Int_s$ is the extension of $Int$ to 
$H_s (\Lambda (M; {\cal E}))$ and  ${\cal C}$ denotes the combinatorial complex
associated to ${\cal E} \rightarrow M, \mu$ and a generalized triangulation 
$\tau = (h, g').$ It turns out that the composition 
$Int_s \cdot (\Delta + Id)^{-s/2}$ is a morphism of trace class which induces 
an isomorphism in algebraic cohomology and, therefore, the corresponding 
mapping cone ${\cal C} (Int_s (\Delta + Id)^{-s/2})$ is algebraically 
acyclic. We show that it admits a torsion, called the relative torsion 
${\cal R}$ associated to ${\cal E} \rightarrow M, g, \mu, \tau,$ and that this
torsion is independent of $s > \frac{n}{2} + 1.$ In the case where ${\cal E}
\rightarrow M$ is of determinant class (cf [BFKM]), one can show that, 
$\log {\cal R} =  \log T_{an} - \log T_{Reid}$ where $T_{an} [T_{Reid}]$ is the 
analytic 
[Reidemeister] torsion.

\spd

\n To prove that $\log {\cal R}$ is a local quantity (in Section 4), we will use 
the Witten de\-for\-ma\-tion of the deRham 
complex, given by the differentials $d_k (t) : = e^{-th} d_k e^{th}.$
There are dif\-fe\-rent possibilities for defining the corresponding deformation 
${\cal R} (t)$ of the relative tor\-sion. We chose one for which the 
variation $\frac{d}{dt} \log {\cal R}(t)$ can be  computed. By 
definition, ${\cal R} (t)$ is the torsion of the mapping cone associated to

\[ (L_2 (\Lambda (M; {\cal E})), d(t)) \stackrel{e^{th}}{\rightarrow} 
(L_2 (\Lambda (M; {\cal E})), d) 
\stackrel{Int_s \cdot (\Delta + Id)^{-s/2}}{\longrightarrow}
{\cal C}.\]

\spd

\n We will show, Theorem 2.1 subsection 2.2 that  the variation 
$\frac{d}{dt} \log {\cal R} (t)$ is given by 

\[\frac{d}{dt} \log {\cal R} (t) = \dim {\cal E}_x \int_M h e (M, g)\]

\n where $e(M,g)$ is the Euler form of the tangent bundle equipped with
the Levi-Civit\`a connection induced from $g.$

\spd

\subsection{Relative torsion}

\n Let $M$ be a closed smooth manifold, 
${\cal E} \stackrel{\pi}{\rightarrow} M$
a smooth bundle of ${\cal A}$-Hilbert modules of finite type equipped with
a smooth flat connection $\nabla,$ which makes the fiberwise multiplication
with elements of ${\cal A}$ \ $\nabla$-parallel. Unlike in [BFKM], or [BFK1]
we do not restrict to the case where the Hermitian structure $\mu$ provided
by the scalar product $\mu_x$ of the fiber ${\cal E}_x$ of 
${\cal E} \ (x \in M)$ is
$\nabla$-parallel.

\spd

\n We denote by ${\cal E}^{\sharp} \rightarrow M$ the dual bundle of 
${\cal E} \rightarrow
M.$ This is a bundle of ${\cal A}$-Hilbert modules of finite type with fiber
${\cal E}^{\sharp}_x,$ the Hilbert space dual to ${\cal E}_x.$
%
%

%
%
The flat connection $\nabla = \nabla_{{\cal E}}$ induces
a flat connection $\nabla^{\sharp} := \nabla_{{\cal E}^{\sharp}}.$ 
${\cal E} \rightarrow M$ and ${\cal E}^{\sharp} \rightarrow M$
are canonically isometric as bundles by an ${\cal A}$-linear isometry,
but this isometry does not intertwine $\nabla$ and $\nabla^{\sharp}$
 unless the 
Hermitian structure $\mu$ is 
$\nabla_{{\cal E}}$-parallel.
If ${\cal E}\to M,\  \nabla $ is induced by the representation $\rho$ of $\Gamma = \pi_1(M, x_0)$
then  ${\cal E}^{\sharp}\to M,\  \nabla ^{\sharp}$ is induced by $\rho ^{\sharp}.$
Each of these connections induces a covariant 
differentiation,

\[d_k \equiv d_{k, \nabla} : \Omega^k (M; {\cal E}) \rightarrow \Omega^{k+1} 
(M; {\cal E})\]

\spd

\n and

\[d^{\sharp}_k \equiv d_{k, \nabla^{\sharp}} : \Omega^k (M; {\cal E}^{\sharp}) 
\rightarrow \Omega^{k+1} (M; {\cal E}^{\sharp})\]

\n where $\Omega^k (M; {\cal E}) := C^{\infty} (\Lambda^k (T^{\ast}M) 
\otimes {\cal E})$
and
where $\Omega^k (M; {\cal E}^{\sharp})$ is defined similarly. The isomorphism
$\theta$ induces a canonical isomorphism, again denoted by 
$\theta,$
between $\Omega^k (M; {\cal E})$ and $\Omega^k (M; {\cal E}^{\sharp}).$  
To simplify
notation, we write $\Lambda^k(M; {\cal E})$ for the ${\cal A}$-Hilbert module of 
finite type, $\Lambda^k (M; {\cal E}) := \Lambda^k (T^{\ast}M) 
\otimes {\cal E}.$

\n Given a Riemannian metric $g$ on $M$, let ${\cal J}_{k, {\cal E}} :
\Lambda^k (M; {\cal E}) \rightarrow \Lambda^{n-k} (M; {\cal E}^{\sharp})$ (with
$n = \dim M$) be the morphism of ${\cal A}$-Hilbert modules defined by 
${\cal J}_{k, {\cal E}} := {\cal J}_k \otimes \theta$ where
${\cal J}_k : \Lambda^k (T^{\ast} M) \rightarrow \Lambda^{n-k} (T^{\ast} M)$
is the Hodge-star operator induced by $g.$

\spd

\n Denote by $d^{\ast}_k : \Omega^{k+1} (M; {\cal E}) \rightarrow 
\Omega^k (M; {\cal E})$
the operator defined by the com\-po\-si\-tion

\spf

\n {\bf (2.1)} \hsps $d^{\ast}_k \equiv d^{\ast}_{k, \nabla} := 
(-1)^{nk-1} {\cal J}_{n-k, {\cal E}^{\sharp}}
\circ d_{n-k-1; \nabla^{\sharp}} \circ {\cal J}_{k+1, {\cal E}}$.

\spd

\n Further we define a scalar product  $\ll \cdot , \cdot \gg : \Omega^k (M; {\cal E}) \times 
\Omega^k (M; {\cal E}) \rightarrow
\bbc$   given by

\spd

\n {\bf (2.2)} \hsps $\ll w_1, w_2 \gg := \int_M w_1 \wedge_{{\cal E}} \ast w_2
\ d \ vol_g$

\spd

\n where $\ast w_2 \equiv {\cal J}_{k, {\cal E}} w_2$ and $\wedge_{{\cal E}}$ is
defined by

\[\Omega^k (M; {\cal E}) \times \Omega^k(M; {\cal E}) \stackrel{Id \times 
{\cal J}_{k, {\cal E}}}{\rightarrow} \Omega^k (M; {\cal E}) \times \Omega^{n-k} 
(M; {\cal E}^{\sharp})\]

\[ \stackrel{\wedge}{\rightarrow} \Omega^n (M, {\cal E} \otimes
{\cal E}^{\sharp}) \stackrel{ev}{\rightarrow} \Omega^n (M; \bbc),\]

\[(w_1, w_2) \mapsto (w_1, \ast w_2) \mapsto w_1 \wedge \ast
w_2 \mapsto w_1 \wedge_{{\cal E}} \ast w_2\]

\n and $ev$ is the map induced by the evaluation map ${\cal E}_x \otimes 
{\cal E}^{\sharp}_x \rightarrow \bbc$. Notice that with respect to the inner
product $\ll \cdot , \cdot \gg, \ d^{\ast}_k$ is the 
formal
adjoint of $d_k,$ i.e. for
$w_1 \in \Omega^k (M, {\cal E}), w_2 \in \Omega^{k+1} (M; {\cal E})$

\[ \ll w_1, d^{\ast}_k w_2 \gg = \ll d_k w_1, w_2 \gg .\]

\spd

\n Introduce the Laplacians

\spd

\n {\bf (2.3)} \hsps $\Delta_k = d_k^{\ast} d_k + d_{k-1} d^{\ast}_{k-1} :
\Omega^k (M; {\cal E}) \rightarrow \Omega^k (M; {\cal E})$

\spd

\n which are second order, elliptic essentially selfadjoint 
nonnegative ${\cal A}$-linear ope\-ra\-tors. In particular, as $M$ is closed,
$(Id + \Delta_k) : \Omega^k (M; {\cal E}) \rightarrow \Omega^k (M; {\cal E})$ is an 
isomorphism of Fr\'echet spaces. Using functional calculus, one can define
the powers $(Id + \Delta_k)^s$ ($s \in \bbr$ arbitrary). They can be used to
define a family of scalar products on 
$\Omega^k (M, {\cal E})$ by setting

\spf

\n {\bf (2.4)} \hsps $\ll w_1, w_2 \gg_s := \ll (Id + \Delta_k)^{s/2} w_1, 
(Id + \Delta_k)^{s/2} w_2 \gg .$

\spd

\n Clearly, $\ll \cdot , \cdot \gg_s$ depends on the Hermitian structure 
$\mu$ and the Riemannian metric $g$. As $d_k \Delta_k = \Delta_{k+1} d_k$
and $\Delta_k d^{\ast}_k = d^{\ast}_k \Delta_{k+1}, d^{\ast}_k$ is also the
adjoint of $d_k$ with respect to the inner product $\ll \cdot , \cdot \gg_s.$
Denote by $H_s (\Lambda^k (M; {\cal E}))$ the completion of 
$\Omega^k (M; {\cal E})$
with respect to the norm induced by the scalar product 
$\ll \cdot , \cdot \gg_s.$
As $M$ is supposed to be closed, one obtains a family of complexes

\[0 \rightarrow \ldots \rightarrow H_s (\Lambda^k (M; {\cal E})) 
\stackrel{d_k}{\rightarrow} H_s (\Lambda^{k+1} (M; {\cal E})) \rightarrow
\ldots \rightarrow 0\]

\n which we denote by $H_s (\Lambda (M; {\cal E})) \equiv (H_s (\Lambda (M; 
{\cal E})), d).$ Here $d_k : H_s (\Lambda^k (M; {\cal E})) 
\linebreak
\rightarrow 
H_s(\Lambda^{k+1}
(M; {\cal E}))$ is the maximal closed extension of 

\[d_k : \Omega^k 
(M; {\cal E}) 
\rightarrow \Omega^{k+1} 
(M; {\cal E}), \]

\n   its domain is $H_{s+1} (\Lambda^k (M; {\cal E})).$ 

\spd

\n As mentioned in Examples 1.6 (4), $H_s (\Lambda (M, {\cal E}))$ are 
$\zeta$-regular complexes of sF type and the morphisms 

\[(1+ \Delta_k)^{t/2} : H_s (\Lambda^k (M; {\cal E})) 
\rightarrow H_{s-t} (\Lambda^k
(M; {\cal E}))\]

\n establish an isometry between the complexes 
$H_s (\Lambda (M; {\cal E}))$ and
$H_{s-t} (\Lambda (M; {\cal E})).$ The adjoint of $d_k$ with respect to the 
$L_2$-inner product, $d_k^{\ast} : H_s (\Lambda^{k+1} (M; {\cal E})) 
\rightarrow
H_s (\Lambda^k (M; {\cal E})),$ is given by the maximal closed extension
of $d^{\ast}_k : \Omega^{k+1} (M; {\cal E}) 
\rightarrow \Omega^k (M; {\cal E})$ (which
depends on $s$). 

\spd

\n Consider a generalized triangulation 
$\tau = (h, g')$ of $M$ (cf introduction of [BFK], also [BFKM])
where $h : M \rightarrow \bbr$ is a Morse function, $g'$ is a compatible 
Riemannian metric on $M$ and ${\cal O}_h$ is a collection of orientations 
of 
the unstable manifolds defined by   $-\rm{grad}_{\tilde{g^{\prime}}}
 \tilde{h},$ compatible with the deck transformations,
where $\tilde{h}$ and $\tilde{g}'$ are lifts of $h$ respectively $g'$ to the
universal cover of $M$  and denote by ${\cal C} \equiv ({\cal C} (M, \tau, 
{\cal T}), \delta)$ the finite  dimensional cochain complex of 
${\cal A}$-Hilbert modules associated with $\tau$ and ${\cal T} : =
({\cal E}, \nabla, \mu).$ The ${\cal A}$-linear map
$Int_k : \Omega^k (M; {\cal E}) \rightarrow {\cal C}^k (M, \tau, {\cal T})$
provided by integration (cf. [BFKM]) are continuous with respect to the
Fr\'echet topology on $\Omega^k (M; {\cal E})$ and extend, for $s > n/2$, to 
bounded maps

\[Int_{k; s} : H_s (\Lambda^k (M; {\cal E})) \rightarrow {\cal C}^k
(M, \tau, {\cal T}).\]

\n Since the integration map intertwines $d$ and $\delta$, they induce 
morphisms

\spf

\n {\bf (2.6)} \hsps $Int_{k;s} : H_s (\Lambda^k (M; {\cal E})) \rightarrow 
{\cal C}^k (M, \tau, {\cal T}).$

\spd

\n As ${\cal C}$ is a complex of finite type Hilbert modules and
$Int_{k;s}$ is bounded for $s > \frac{n}{2}, \ Int_{k;s}$ is of trace class
(cf Examples 1.6 (2)). For $s' \le s,$ denote by $In_{\ast; s,s'}$ the
embedding

\spd

\n {\bf (2.7)} \hsps $In_{k; s,s'} : H_s (\Lambda^k (M; {\cal E})) \hookrightarrow
H_{s'} (\Lambda^k (M; {\cal E}))$

\spd

\n which is a morphism of cochain complexes. For $s-s' > n, \ In_{k; s,s'}$
is of trace class.

\spd

\begin{prop} (cf [CMM, section 4]) 

\n The following families of morphisms induce isomorphisms in al\-ge\-bra\-ic 
co\-ho\-mo\-lo\-gy:

\begin{itemize}
\item[(i)] $(Id + \Delta)^{s/2} : H_r (\Lambda (M; {\cal E})) \rightarrow H_{r-s}
(\Lambda (M; {\cal E})) \ \ (\forall r, s);$

\item[(ii)] $In_{s,s'} : H_s (\Lambda (M; {\cal E})) \rightarrow H_{s'} (\Lambda
(M; {\cal E})) \ \ (s \ge s');$

\item[(iii)] $Int_s : H_s (\Lambda (M; {\cal E})) \rightarrow {\cal C} (M, \tau,
{\cal T}) \ \ (s > n/2).$
\end{itemize}
\end{prop}

\spd

\n {\bf Proof} (i) is obvious since $(Id + \Delta)^{s/2}$ is an isometric
morphism.

\spd

\n (ii) For ease of notation, let $H_{k;s} := H_s (\Lambda^k (M; {\cal E})).$
Denote by $Q_k : H_{k;0} \rightarrow H_{k;0}$ the $L_2$-orthogonal 
spectral projector corresponding to the interval $[0, r],  r>0$ and consider its 
restriction to $H_{k;s}, \ s>0.$ 
 Due to the ellipticity of $\Delta_k, \ Q_k \omega
\in \Omega^k (M; {\cal E})$ and $Q_k (H_{k;s})$ is isomorphic to $Q_k (\Omega^k
(M; {\cal E})).$ As $H_{k,s} = Q_k (H_{k;s}) \oplus (Id - Q_k) (H_{k;s}), 
H_s (\Lambda (M; {\cal E}))$ is the sum of two subcomplexes, 
${\cal C}^1 \oplus {\cal C}^2$ where ${\cal C}^1_k := Q_k (H_{k;s})$ and 
${\cal C}^2_k := (Id - Q_k)(H_{k,s}).$ Notice that ${\cal C}^2$ is algebraically
acyclic
by remark after (1.20).
${\cal C}^1$ has the same algebraic co\-ho\-mo\-lo\-gy as $H_s (\Lambda (M,
{\cal E}))$ and that
$In_{k;s,s'} = \bl {l|l} Id & 0 \\
\cline{1-2} 
0 & In_{k;s, s'|_{{\cal C}^2_k}} \er.$ Therefore, $In_{s,s'}$, induces an 
isomorphism
in algebraic cohomology.

\spd

\n (iii) We use Witten's deformation of the deRham complex (cf [BFK2, section 3], also [BFKM, 
section 5]). For $t \in \bbr,$ let $d_k (t) := e^{-th} d_k e^{th}$  be the
deformed exterior differential 
where $h : M \rightarrow \bbr$ is the 
Morse function of a generalized triangulation $\tau = (h, g')$ ( h not necessarily selfindexing).
For $t$ 
sufficiently large, the spectrum of the deformed Laplacian $\Delta_k (t)$
splits into a small part contained in the interval [0,1] and a large part,
contained in an interval of the form  $[Ct, \infty]$ with $C > 0.$ Denote by
$(\Omega_{sm} (M; {\cal E}) (t),$ 
$d(t))$ the subcomplex associated to the small
part of the spectrum. In particular, 
$\Omega^k_{sm} (t) \equiv \Omega^k_{sm} (M; {\cal E})(t) = 
Q_k (t) (H_{k;s})$
where $Q_k (t)$ denotes the or\-tho\-go\-nal spectral projector of 
$\Delta_k (t)$ 
corresponding to the interval [0,1]. As $\Delta_k(t)$ 
\linebreak
is elliptic, one 
concludes that
$\Omega^k_{sm} (t) \equiv \Omega^k_{sm} (M; {\cal E}) (t) \subseteq 
\Omega^k (M; {\cal E})$ 
and, by [BFK2, section 3] (cf also [BFKM, section 5]), is an ${\cal A}$-Hilbert module of finite type.
Similarly as in the proof 
of (ii), $H_s (\Lambda (M; {\cal E}))$ is
the sum of two
subcomplexes, $\Omega_{sm} (t) \oplus H_{s, la} (t),$ 
where $H_{k;s,la} :=
(Id -Q_k (t)) (H_{k;s}).$
Introduce

\spf

\n {\bf (2.8)} \hsps $f_k (t) : \Omega^k_{sm} (M; {\cal E})(t) \rightarrow
{\cal C}^k (M, \tau, {\cal T})$

\spd

\n given by $f_k (t) := Int_k e^{th}$. 
One verifies that $f(t) : (\Omega_{sm} (t), d (t)) \rightarrow 
({\cal C}, \delta)$ is a morphism of cochain com\-ple\-xes.
By the Helffer-Sj\"ostrand theory 
(cf [BFK2, section 3], also [BFKM, section 5]), $f(t)$ is an 
isomorphism for $t$ sufficiently large. In particular, it induces an
isomorphism in algebraic cohomology. For $t$ sufficiently large, 
$H_s (\Lambda 
\linebreak
(M; {\cal E}))$ can be decomposed into two subcomplexes, 
${\cal C}^1 \oplus {\cal C}^2$ where ${\cal C}^1 = (e^{th} \Omega_{sm}
\linebreak
(t), d)$ and ${\cal C}^2 = (e^{th} H_{s, la} (t), d).$ Notice that ${\cal C}^2$
is algebraically acyclic and that ${\cal C}^1$ has the same algebraic cohomology
as $H_s (\Lambda (M; {\cal E})).$ For $s > n/2,$ 

\[Int_s : H_s 
(\Lambda 
(M; {\cal E})) =
{\cal C}^1 \oplus {\cal C}^2 \rightarrow {\cal C}\]

\n is well defined and takes
the form $(f(t) e^{-th}, Int_s |_{{\cal C}^2}).$ Therefore $Int_s$ 
induces an isomorphism in algebraic cohomology. \ \ \ \carre

\spf

\begin{rema} For $s-s' > n + 1$ the inclusion  $In_{s,s'}$ is a morphism of
trace class (cf Definition 1.9). \end{rema}

\spd

\n As $H_s (\Lambda (M; {\cal E}))$ is a $\zeta$-regular complex of sF type (for
any $s \in \bbr),$ we conclude from Proposition 2.1 and Lemma 1.12 that the
mapping cone ${\cal C} (In_{s,s'})$ is al\-ge\-brai\-cal\-ly acyclic and thus has a 
well defined torsion $T({\cal C} (In_{s,s'} )) (s-s' > n +1)$.

\spf

\begin{prop} For $s-s' > n + 1, \ \log T ({\cal C} (In_{s,s'})) = 0.$
\end{prop}

\spf

\n {\bf Proof} Recall that the mapping cone ${\cal C} (In_{s,s'})$ is 
given by 

\[{\cal C} (In_{s,s'})_k := H_{s'} 
(\Lambda^{k-1} (M; {\cal E})) \oplus
H_s
(\Lambda^k (M; {\cal E}));\]

\[d(In_{s,s'} )_k = \bl {cc}
-d_{k-1} & In_{k; s,s'} \\
0 & d_k \er.\]

\n Therefore, by (1.23) and (2.4)
, 

\spf

\n {\bf (2.9)} \hsps $\Delta (In_{s,s'})_k = \bl {cc}
\Delta_{k-1} + (Id+\Delta_{k-1})^{s'-s} & 0 \\
0 & \Delta_k + (Id+\Delta_k)^{s'-s} \er$

\spd

\n Proposition 2.3 then follows from Lemma 2.4 below. \ \ \ \carre

\spf

\begin{lemma} Assume that $({\cal C}, \delta)$ is a $\zeta$-regular cochain 
complex of sF type. Then 
\begin{enumerate}
\item[(i)] $\sum_k (-1)^k \log \ \det (\Delta_k + Id) = 0$;
\item[(ii)] $\sum_k(-1)^k \log \ \det (\Delta_k+(Id+\Delta_k)^{s'-s})=0$.
\end{enumerate}
\end{lemma}

\spd

\n {\bf Proof} As the two statements can be proved in the same way, we consider 
only the second one.  With respect to the Hodge decomposition, $\Delta_k + Id$ 
takes the form (cf (1.18))

\spf

\n {\bf (2.10)} \hsps $\Delta_k + (Id+\Delta_k)^{s'-s} = $
\[= diag (Id, 
\underline{\Delta}^+_k + (Id+\underline{\Delta}_k^+)^{s'-s},
\underline{\Delta}^-_k + (Id+\underline{\Delta}_k^-)^{s'-s}).\]

\spd

\n where $\underline{\Delta}^+_k := \underline{d}_{k-1} 
\underline{d}^{\ast}_{k-1}$ and $\underline{\Delta}^-_k := 
\underline{d}^{\ast}_k \underline{d}_k$ are $\zeta$-regular operators of sF
type and 
 
\[\log \ \det (\Delta_k + (Id+\Delta_k)^{s'-s}) = \log \ \det (\underline{\Delta}^+_k 
+ (Id+\underline{\Delta}_k^+)^{s'-s}) +\]
\[+\log \ \det (\underline{\Delta}^-_k + 
+(Id+\underline{\Delta}_k^-)^{s'-s}).\]

\n As $\underline{\Delta}^+_k$ and $\underline{\Delta}^-_{k-1}$ have the same 
spectral distribution function we conclude that

\spf

\n {\bf (2.11)} \hsps $\log \ \det (\underline{\Delta}^+_k + (Id+\underline{\Delta}_k^+
)^{s'-s}) = \log \ \det
(\underline{\Delta}_{k-1}^- + (Id+\underline{\Delta}_{k-1}^-)^{s'-s}).$ 

\spd

\n This proves the claimed statement. \ \ \ \carre

\spf

\begin{prop} Let $s > \frac{n}{2} +1, \tilde{s} > s + n + 1$ and $p \in \bbr_{>0}$.
Then

\begin{itemize}
\item[(i)] $Int_s : H_s (\Lambda(M; {\cal E})) \rightarrow
 {\cal C} (M, \tau, {\cal T})$
\end{itemize}

\n is a morphism of trace class which induces an isomorphism in algebraic 
co\-ho\-mo\-lo\-gy. Moreover

\spf

\n {\bf (2.12)} \hsps $\log T ({\cal C} (Int_{\tilde{s}})) = \log T ({\cal C} 
(Int_s))$

\spd

\n and

\spf

\n {\bf (2.13)} \hsps $\log T ({\cal C} (Int_s)) = \log T ({\cal C} (Int_s 
(\Delta + Id)^{p/2} )),$  

\spd

with  $(\Delta + Id)^{p/2}:  H_{s+p} (\Lambda(M; {\cal E})) \rightarrow
 H_s (\Lambda(M; {\cal E}))$ \end{prop}

\spd

\n {\bf Proof} By Proposition 2.1, $Int_s$ is a morphism of cochain complexes 
which induces an iso\-mor\-phism in algebraic cohomology. We claim that $Int_s$ is
a mor\-phism of trace class (cf Definition 1.10).
As we have already observed, $Int_{k;s}$ is an operator of trace class for
$s > n/2$ as it is bounded and ${\cal C}^k$ is a Hilbert module of finite
type.
By the same reason, $\delta^{\ast}_k Int_{k + 1; s}$ is of trace class for 
$s > n/2.$ Finally we have to verify that $Int_{k;s} d^{\ast}_k$ is of trace
class. Use that $Int_{k;s} = Int_{k; s-1} In_{k;s, s-1}$ and that 
$In_{k; s, s-1} d_k^{\ast}$ is a bounded operator to deduce from 
Proposition 1.3 that $Int_{k;s} d^{\ast}_k$ is of trace class for 
$s > \frac{n}{2} + 1.$
Further notice that $Int_{\tilde{s}} = Int_s \cdot In_{\tilde{s}, s}$ for
arbitrary $\tilde{s} > s + n + 1$ where $In_{\tilde{s},s}$ is a
morphism of trace class (Remark 2.2) which induces an isomorphism in
algebraic cohomology (Proposition 2.1). From Proposition 1.15 and
Proposition 2.3 we then conclude that

\[\begin{array}{lll}
\log T ({\cal C} (Int_{\tilde{s}} )) & = & \log T ({\cal C} (In_{\tilde{s},s}))
+ \log T ({\cal C} (Int_s)) \\
& = & \log T ({\cal C} (Int_s)).
\end{array}\]

\n It follows from Proposition 1.3 and the fact that $(\Delta + Id)^{p/2}$ 
is an isometry of cochain complexes that $Int_s (\Delta+Id)^{p/2}$ is a
morphism of trace class which induces an isomorphism in algebraic cohomology.
Applying Lemma 1.11 and Proposition 1.16 yields formula (2.13). \ \ \ \carre

\spd

\n We are now ready to define the relative torsion. Our definition 
differs slightly from the one first introduced by [CMM].

\spd

\n Given $M, {\cal T}, g, \tau$, define ${\cal R}_s = {\cal R}_s (M, {\cal T},
g, \tau),$

\spf

\n {\bf (2.14)} \hsps $\log {\cal R}_s := \log T ({\cal C} (Int_s 
(\Delta + Id)^{-s/2})) \hsps (s > n/2 +1).$

\spd

\n Notice that, by Proposition 2.5, $\log {\cal R}_s = \log T ({\cal C} (Int_s))$ 
and ${\cal R}_s$ is independent of $s > \frac{n}{2} + 1.$

\spf

\begin{definition} The \underline{relative torsion} 
is defined by 
${\cal R} := {\cal R}_s (M, {\cal T}, g, \tau).$ 
\end{definition}

\spd

\n The name of relative torsion is justified by the observation, due to
[CMM], that in the case where ${\cal E} \rightarrow M$ is of determinant class - 
and therefore, the analytic torsion $T_{an}$ and the Reidemeister torsion 
$T_{Reid}$ are well defined (cf [BFKM]) - ${\cal R}$ is given by 

\spf

\n {\bf (2.15)} \hsps $\log {\cal R} = \log T_{an} - \log T_{Reid}$

\spd

\n To see that (2.15) hold we would like to apply Proposition 1.10 (Milnor's
Lemma). However $H_s (\Lambda (M; \xi))$ is not a complex of finite dimension, 
we therefore decompose the mapping cone ${\cal C} (Int_s)$ as follows:
As shown in the proof of Proposition 2.1 (iii), 
$H_s (\Lambda (M; {\cal E})) = {\cal C}^1 \oplus {\cal C}^2$ and $f : Int_s
|_{{\cal C}^1} : {\cal C}^1 \rightarrow {\cal C}^2$ is a bijective morphism
where ${\cal C}^2$ is algebraically acyclic and ${\cal C}^1$ is given by
${\cal C}^1 := (e^{th} \Omega_{sm} (t), d).$ Notice that ${\cal C}^1$ has the
same algebraic cohomology as $H_s (\Lambda (M; {\cal E})).$ Denote by ${\cal C}^3$
the mapping cone ${\cal C} (f).$ As $f$ is a bijective morphism, ${\cal C}^3 =
{\cal C} (f)$ is algebraically acyclic. As $f$ satisfies the assumption of 
Lemma 1.14 one concludes that

\[\log {\cal R} = \log T ({\cal C} (Int_s)) = \log T({\cal C}^3) + \log T 
({\cal C}^2).\]

\n The term $\log T ({\cal C}^3)$ has to be analyzed further.
By assumption, $(M, \rho)$ is of determinant class. 
Therefore ${\cal C}^1$ is of determinant class as well. We then obtain the
following short exact sequence of cochain complexes, each of them being of 
determinant class,

\[0 \rightarrow \Sigma {\cal C} \rightarrow {\cal C}^3 \rightarrow {\cal C}^1
\rightarrow 0.\]

\spd

\n As $\Sigma {\cal C}, {\cal C}^3$ and ${\cal C}^1$ are of finite dimension
we can apply Proposition 1.10 to conclude that

\[\log T ({\cal C}^3) = \log T(\Sigma {\cal C}) + \log T({\cal C}^1) + \log T 
({\cal H})\]

\n where ${\cal H}$ denotes the long weakly exact sequence in reduced 
cohomology and is of determinant class. By Lemma 1.14 and the 
de\-fi\-ni\-tion of
the combinatorial torsion $T_{comb}$ (cf [BFKM]), $\log T(\Sigma {\cal C}) = - 
\log T({\cal C}) = - \log T_{comb}.$ By the de\-fi\-ni\-tion of the analytic torsion
(cf [BFKM]), $\log T({\cal C}^1) + \log T({\cal C}^2) = \log T_{an}.$ Further
${\cal H}$ is given by

\[\ldots \rightarrow \bar{H}^i(\Sigma {\cal C}) \rightarrow 
\bar{H}^i ({\cal C}^3)
\rightarrow \bar{H}^i ({\cal C}^1) \stackrel{\bar{H} (\delta_i)}{\rightarrow}
\bar{H}^{i+1} (\Sigma {\cal C}) \rightarrow \ldots \]

\n Recall that ${\cal C}^3$ is algebraically acyclic and 
$\bar{H}(\delta_i) : \bar{H}^i
({\cal C}^1) \rightarrow \bar{H}^i({\cal C})$ is given by the restriction
of $Int_{i;s}$ to ${\cal H}_i = ker (\Delta_i)$ and, by  the definition of the
metric part $T_{met}$ of the torsion

\[\log T ({\cal H}) =  - \log T_{met}.\]

\n Combining the various equalities above leads to (2.15).

\spd

\subsection{Witten deformation of the relative torsion}

\n To obtain a local formula for $\log {\cal R}$ we consider the 
deformed relative torsion ${\cal R}(t)$ obtained by considering the Witten
deformation of the deRham complex. 

\spd

\n Using the generalized triangulation $\tau = (h, g', {\cal O}_h)$ we 
consider for $s > n/2 +1$ the following composition $ g_s (t)$ of
morphisms.

\spd

\n {\bf (2.16)} \  $(L_2 (\Lambda (M; {\cal E})), d (t)) 
\stackrel{e^{th}}{\longrightarrow} (L_2 (\Lambda (M; {\cal E})), d) 
\stackrel{(\Delta + Id)^{-s/2}}{\rightarrow} (H_s (\Lambda (M, {\cal E})), d) $

\[\stackrel{Int_s}{\rightarrow} ({\cal C}, \delta).\]

\spd

\n where we recall that $d_k (t) = e^{-th} d_k e^{th}$. For $s > n/2 + 1,$
$Int_s$ is a morphism of trace class (cf Proposition 2.5) and thus, by
Proposition 1.3, the composition $g_s (t)$ is a morphism of trace class. Due to
Proposition 2.1, $g_s (t)$ induces an isomorphism in algebraic cohomology, hence
$\log T({\cal C} (g_s(t)))$ is well defined.

\spf

\begin{prop} For $s + p > s > 2n+2$ and $t \in \bbr$,

\[\log T ({\cal C} (g_s (t))) = \log T({\cal C}(g_{s+p} (t))).\]
\end{prop}

\spf

\n {\bf Proof} Using that $Int_s = Int_{s/2} In_{s, s/2}, \ g_s(t)$ is given by
$g_s(t) = Int_{s/2} \alpha_s(t)$ where $\alpha_s(t) := In_{s, s/2} (\Delta + Id)^{-s/2}
e^{th}$.
As $s > 2n+2, \ \alpha_s(t)$ (cf. Remark 2.2) and $Int_{s/2}$ (cf. Proposition
2.5) are morphisms of trace class. By Proposition 2.1, both $Int_{s/2}$ and
$\alpha_s(t)$ induce isomorphisms in algebraic cohomology. Therefore we can apply
Proposition 1.16 to conclude that

\spf

\n {\bf (2.17)} \hsps $\log T ({\cal C}(g_s(t))) = \log T ({\cal C} (Int_{s/2} ))
+ \log T ({\cal C} (\alpha_s(t))).$

\spd

\n Observe that, $In_{s+p, \frac{s}{2} + p} (\Delta + Id)^{-\frac{s+p}{2}} =
(\Delta + Id)^{-p/2} In_{s, s/2} (\Delta + Id)^{-s/2}$
to deduce that

\[g_{s+p}(t) = Int_{p+ \frac{s}{2}} (\Delta + Id)^{-p/2} \alpha_s(t).\]

\n Using the same arguments as above we can apply Proposition 1.16 to obtain

\spf

\n {\bf (2.18)} \hsps $\begin{array}{lll}
\log T ({\cal C} (g_{s+p})) & = & \log T ({\cal C} (Int_{p+s/2} \cdot (\Delta +
Id)^{-p/2})) + \\
&& + \log T ({\cal C} (\alpha_s)).
\end{array}$

\spd

\n By Proposition 2.5,

\spf

\n {\bf (2.19)} \hsps $\log T({\cal C}(Int_{\frac{s}{2} +p})) = 
\log T ({\cal C} 
(Int_{\frac{s}{2} + p} (\Delta + Id)^{-p/2}))$

\spd

\n and (cf (2.14))

\spf

\n {\bf (2.20)} \hsps $\log T ({\cal C} (Int_{s/2 + p} )) = \log T ({\cal C} 
(Int_{s/2})).$

\spd

\n Combining (2.17)-(2.20) leads to the claimed result. \ \ \ \carre

\spd

\n In view of Proposition 2.7, 
$\log T ({\cal C}(g_s (t)))$
is independent of
$s$.

\spf

\begin{definition} The Witten deformation ${\cal R}(t)$ of the relative torsion
is defined by

\[\log {\cal R}(t) := \log T ({\cal C}(g_s (t))) \ (s > n +2).\]
\end{definition}

\spd

\n Using the same arguments as for the verification of (2.15) one sees that,
for ${\cal E} \rightarrow M$ of determinant class,

\spf

\n {\bf (2.21)} \hsps $\log {\cal R}(t) = \log T_{an} (t) - \log T_{comb} - \log
T_{met} (t)$

\spd

\n where $T_{an} (t)$ denotes the analytic torsion of the deformed deRham 
complex (cf [BFKM]), $T_{comb}$ denotes again the combinatorial torsion and
$\log T_{met} (t) = -\log T ({\cal H} (t))$ is the torsion of the long weakly
exact sequence in cohomology obtained in a similar fashion as in the case of 
(2.15).

\spd

\n Our next objective is to calculate the variation for $\log {\cal R}(t):$

\spf

\begin{theor}

\[\frac{d}{dt} \log {\cal R}(t) = \dim \ {\cal E}_x \int_M h \cdot e(M,g)\]

\n where $e(M,g)$ is the Euler form of the tangent bundle equipped with the
Levi-Civit\`a connection induced from $g.$ In particular, $\frac{d}{dt} \log
{\cal R}(t)$ is a local quantity, independent of $t$ and the Hermitian
structure. Moreover if $n = \dim M$ is odd, then

\[\frac{d}{dt} \log {\cal R}(t) \equiv 0.\]
\end{theor}

\spf

\n {\bf Proof} For $s > n+2,$ introduce the morphism $g_s := Int_s (\Delta +
Id)^{-s/2}.$
Then $g_s(t) = g_s e^{th}.$ Observe that for all $t$, the mapping cone 
${\cal G} (t) := {\cal C} (g_s (t))$ is a cochain complex whose $k' th$ component
${\cal G}_k : = {\cal C}(g_s(t))_k$ is independent of $t$, given by

\spf

\n {\bf (2.22)} \hsps ${\cal G}_k = {\cal C}^{k-1} (M, \tau, {\cal T}) 
\oplus L_2 (\Lambda^k
(M; {\cal E})).$

\spd

\n With respect to the decomposition (2.22),
$d_k (t) \equiv d^{{\cal G}}_k (t) : 
{\cal G}_k \rightarrow {\cal G}_{k+1}$ takes the form

\spf

\n {\bf (2.23)} \hsps $d_k (t) = E_{k+1} (-t) D_k E_k (t); \ d_k (t)^{\ast} =
E_k (t) D_k^{\ast} E_{k+1} (-t).$

\spd

\n where

\[ E_k (t) := \bl{cc} Id & 0 \\
0 & e^{th} \er \ , \ D_k := \bl {cc}
- \delta_{k-1} & g_{k;s} \\
0 & d_k \er \]

\n According to Proposition 2.5 and 1.12, the mapping cone ${\cal C} (g_s (t))$
is a $\zeta$-regular, algebraically acyclic cochain complex of sF type. Denote
by 

\[\Delta_k (t) \equiv \Delta^{{\cal G}}_k (t) \ \mbox{the} \ k\mbox{-Laplacian
associated to} \ {\cal G} \equiv {\cal C} (g_s (t)),\]

\[\Delta_k(t) := \Delta^+_k (t) + \Delta^-_k(t)\]

\n where

\[\Delta^+_k (t) := d_{k-1} (t) d^{\ast}_{k-1} (t) \ , \ \Delta^-_k (t) :=
d^{\ast}_k (t) d_k(t).\]

\n With respect to the Hodge decomposition ${\cal G}_k = {\cal G}^+_k (t) 
\oplus {\cal G}^-_k (t),$

\spf

\n {\bf (2.24)} \hsps ${\cal G}^+_k (t) := 
\overline{d_{k-1} (t) ({\cal G}_{k-1})} ; 
\ {\cal G}_k^- (t)
: = \overline{d^{\ast}_k (t) ({\cal G}_{k+1})},$

\spd

\n $d_k (t)$ is of the form $d_k (t) = \bl {cc}
0 & 0 \\
0 & \underline{d}_k (t) \er$ and thus $\Delta^+_k (t) = \bl {cc}
\underline{\Delta}^+_k (t) & 0 \\
0 & 0 \er$ and $\Delta^-_k (t) = \bl {cc}
0 & 0 \\
0 & \underline{\Delta}^-_k (t) \er$ where $\underline{\Delta}^+_k (t) := 
\underline{d}_{k-1} (t) \underline{d}^{\ast}_{k-1} (t)$ and 
$\underline{\Delta}^-_k (t) := \underline{d}^{\ast}_k (t) \underline{d}_k (t).$
As ${\cal C} (g_s (t))$ is a $\zeta$-regular cochain complex, the operators
$\underline{d}_k (t),$ according to Definition 1.9, satisfy Op(1)-Op(7). As
${\cal C} (g_s (t))$ is algebraically acyclic, they also satisfy Op(8). (In
fact, the operators $\underline{d}_k(t)$ satisfy $Op (7)$,
because the corresponding
heat traces admit an asymptotic expansion of the form (Asy) (cf Lemma 1.1).)

\spd

\n To compute $\frac{d}{dt} \log {\cal R}(t)$ we proceed similarly as in
[BFKM] where we compute $\frac{d}{dt} \log T_{an} (t).$ According to 
Definition 2.8, the relative torsion ${\cal R}(t) = T({\cal C}(g_s (t)))$ is
given by 

\[\log {\cal R} (t) = \frac{1}{2} \sum_k (-1)^{k+1} \log \ \det 
\underline{\Delta}^+_k (t)\]

\n and

\spf

\n {\bf (2.25)} \hsps $\frac{d}{dt} \log \det (\underline{\Delta}^+_k (t)) = 
F.p._{z=0} Tr (\frac{d}{dt} (\underline{\Delta}^+_k (t)) (\underline{\Delta}^+_k
(t))^{-z-1})$

\[= F.p._{z=0} Tr (\frac{d}{dt} (E_k (-t)  D_{k-1} E_{k-1} (2t) 
D^{\ast}_{k-1} E_k (-t)) (\Delta^+_k (t))^{-z-1}).\]

\n (By $F.p._{z=0}$  we denote the constant term in the Laurent expansion  
at $z=0$ of the
expression which follows it.) To simplify writing we suppress momentarily the
subscript $k.$ Then

\spf

\n {\bf (2.26)} \ $\begin{array}{l}
\frac{d}{dt} (E (-t) D E (2t) D^{\ast} E (-t)) = \\
\ = \bl {cc} 0 & 0 \\ 0 & -h \er E(-t) D E (2t) D^{\ast} E (-t) \\
\ + E (-t) D \bl {cc} 0 & 0 \\ 0 & 2h \er  E (2t) D^{\ast} E (-t) \\
\ + E (-t) D E(2t) D^{\ast} E (-t) \bl {cc}
0 & 0  \\
0 & -h \er . \end{array}$

\spd

\n Substituting (2.26) into (2.25) we obtain

\spf

\n {\bf (2.27)} \ $\begin{array}{l}
\frac{d}{dt} \log \ \det (\underline{\Delta}^+ (t)) = \\
\ = F.p._{z=0} Tr \Big( \bl {cc} 0 & 0 \\ 0 & -h \er
d (t) d(t)^{\ast} \Delta^+ (t)^{-z-1}\Big) \\
\ + F.p._{z=0} Tr \Big(d(t) \bl {cc} 0 & 0 \\
0 & 2h \er d(t)^{\ast} \Delta^+(t)^{-z-1}\Big) \\
\ + F.p._{z=0} Tr \Big(d(t) d(t)^{\ast} \bl {cc} 0 & 0 \\
0 & -h \er \Delta^+ (t)^{-z-1} \Big).
\end{array}$

\spd

\n Using that $d (t) d(t)^{\ast} \Delta^+ (t)^{-z-1} = \Delta^+ (t)^{-z}$
together with the commutativity of the trace we conclude

\spf

\n {\bf (2.28)} \hsps $\begin{array}{l}
\frac{d}{dt} \log \det (\underline{\Delta}^+_k (t)) = \\
\ = F.p._{z=0} Tr\Big(\bl {cc}
0 & 0 \\
0 & -2h \er \Delta^+_k (t)^{-z}\Big) \\
\ + F.p._{z=0} Tr \Big( \bl {cc} 0 & 0 \\
0 & 2h \er d^{\ast}_{k-1} (t) \Delta^+_k (t)^{-z-1} d_{k-1} (t)\Big).
\end{array}$

\spd

\n From $d^{\ast}_{k-1} (t) \Delta^+_k (t) = d^{\ast}_{k-1} (t) d_{k-1} (t)
d_{k-1} (t)^{\ast} = \Delta^-_{k-1} (t) d_{k-1} (t)^{\ast}$
we conclude that

\[d^{\ast}_{k-1} \Delta^+_k (t)^{-z-1} d_{k-1} (t) = (\Delta^-_{k-1} (t))^{-z-1}
\Delta^-_{k-1} (t) = \Delta^-_{k-1} (t)^{-z}\]

\n and therefore, after substituting into (2.28)

\spf

\n {\bf (2.29)} \hsps $\frac{d}{dt} \log \ \det (\underline{\Delta}^+_k (t)) =$

\[ = F.p._{z=0} Tr \Big(\bl {cc}
0 & 0 \\ 0 & -2h \er \
(\Delta^+_k (t)^{-z} - \Delta^-_{k-1} (t)^{-z})\Big).\]

\spd

\n This leads to

\spf

\n {\bf (2.30)} \hsps $\frac{d}{dt} \log T ({\cal C} (g_s (t))) = \sum_k (-1)^k 
F.p._{z=0} Tr \Big( \bl {cc} 0 & 0 \\ 0 & h \er \Delta_k (t)^{-z} \Big).$

\spd

\n As ${\cal C} (g_s (t))$ is algebraically acyclic, $\Delta_k (t)$ has no 
spectrum near $0$ and therefore $\zeta^{\rm II} (z) := \frac{1}{\Gamma(z)} 
\int^{\infty}_{1} \lambda^{z-1} Tr \Big( \bl {cc} 0 & 0 \\ 0 & h \er e^{-\lambda
\Delta_k (t)} \Big) d \lambda$ is an entire function, with
$\zeta^{\rm II} (0) = 0.$

\spd

\n In view of the fact that the heat trace $Tr \Big( \bl {cc} 0 & 0 \\
0 & h \er e^{-\lambda \Delta_k(t)}\Big)$ admits an asymptotic expansion at
$\lambda = 0$ of the form (Asy) one obtains (cf Lemma 1.1 and Op(7))

\spf

\n {\bf  (2.31)} \ $\begin{array}{ll}
F.p._{z=0} & Tr \Big( \bl{cc} 
0 & 0 \\
0 & h \er
\Delta_k (t)^{-z} \Big) = \\
& = F.p._{z=0} \Big( \frac{1}{\Gamma (z)} \int^{\infty}_0 \lambda^{z-1} Tr
( \bl {cc} 0 & 0 \\
0 & h \er e^{-\lambda \Delta_k(t)}) d \lambda\Big) \\
& =  F.p._{z=0} \Big( \frac{1}{\Gamma(z)} \int^1_0 \lambda^{z-1} Tr (\bl {cc}
0 & 0 \\
0 & h \er e^{-\lambda \Delta_k (t)} ) d \lambda\Big).
\end{array}$

\spd

\n With respect to the decomposition (2.22) of the mapping cone ${\cal C}
(g_s (t)),$ the Laplacian $\Delta_k (t) \equiv \Delta^{{\cal G}}_k (t)$ 
takes the form

\[ \Delta^{{\cal G}}_k (t) = \bl {cc} 0 & 0 \\
0 & \Delta_k (t) \er + \Phi_k (t)\]

\n where $\Phi_k (t) : {\cal G}_k \rightarrow {\cal G}_k$ is given by

{\arraycolsep0.3mm
\[\Phi_k (t) : = \bl {l|l}
\Delta^{comb}_{k-1} + g_{k-1; s} (t) g_{k-1; s} (t)^{\ast} & 
- \delta^{\ast}_{k-1} g_{k;s} (t) + g_{k-1;s} (t) d^{\ast}_{k-1} (t) \\
\cline{1-2}
-g_{k;s} (t)^{\ast} \delta_{k-1} + d_{k-1} (t) g_{k-1;s} (t)^{\ast} & 
g_{k;s} (t)^{\ast} g_{k;s} (t) 
\er.\]}

\n As $g_{k;s} (t)$ and $g_{k;s} (t)^{\ast}$ are bounded operators of 
trace
class
for $s> \frac{n}{2} +1$ (cf Proposition 2.5) the maps $g_{k;s} (t) d^{\ast}_k (t)$
and $d_k (t) g_{k;s} (t)^{\ast}$ are bounded operators of trace class for 
$s > n/2 + 1$ and we conclude that $\Phi_k (t)$ is a bounded operator of trace 
class.
By Proposition 1.5 and its remark, we conclude that

\spd

\n {\bf (2.32)} \ $F.p._{z=0} \Big( \frac{1}{\Gamma(z)} \int^1_0 \lambda^{z-1}
Tr (\bl {cc} 0 & 0 \\
0 & h \er e^{-\lambda \Delta^{{\cal G}}_k (t)}) d\lambda\Big) = $

\[= F.p._{z=0} \Big(\frac{1}{\Gamma(z)} \int^1_0 \lambda^{z-1} Tr (h e^{- \lambda
\Delta_k (t)} (Id - P_{\Delta_k (t)} (\{ 0 \} )) \Big) d\lambda\]

\n where $P_{\Delta_k(t)}(\{ 0\})$ denotes the orthogonal projection onto 
$Ker \Delta_k (t).$

\spd

\n One verifies (cf [BZ, p 81]) that $e^{th} \Delta_k (t) e^{-th}$
are equal to the Laplacians $\Delta^{\prime}_k$
associated
to the undeformed deRham complex, the Riemannian metric $g$ and the
Hermitian structure $e^{-2th} \mu.$

\spd

\n Combining this together with (2.30)-(2.32), we obtain

\spf

\n {\bf (2.33)} $
\begin{array}{l}
\frac{d}{dt} \log {\cal R} (t) = \\
= \sum_k (-1)^k F.p._{z=0} \Big( \frac{1}{\Gamma(z)} \int^1_0 \lambda^{z-1} Tr
(he^{-\lambda \Delta_k (t)} (Id - P_{\Delta_k (t)} (\{ 0 \} ))) d \lambda \Big)
\\
= \sum_k (-1)^k F.p._{z=0} \Big( \frac{1}{\Gamma(z)} \int^1_0 \lambda^{z-1} Tr
(h e^{-\lambda \Delta^{\prime}_k} (Id - P_{\Delta^{\prime}_k} (\{ 0 \} ))) d
\lambda\Big)
\end{array}
$

\spd

\n where we used that, as $\Delta_k (t) = e^{-th} \Delta^{\prime}_k e^{th},
e^{- \lambda \Delta_k (t)} = e^{-th} e^{- \lambda \Delta^{\prime}_k} e^{th}$
and $e^{th} (Id - P_{\Delta_k(t)} (\{ 0 \})) e^{-th} = Id - P_{\Delta^{\prime}_k}
(\{ 0 \}) .$

\spd

\n It is known (cf e.g. [BFKM]) that the restriction $K_k (x; \lambda, t)$ of 
the Schwartz kernel of $e^{- \lambda \Delta^{\prime}_k} 
(Id - P_{\Delta^{\prime}_k}
(\{ 0 \} ))$ to the diagonal has an asymptotic expansion with respect to 
$\lambda$ for $\lambda \searrow 0$ of the form
$\sum_{j \ge 0} A_{k; n-j} (t) \lambda^{- \frac{n-j}{2}}$ where $A_{k; n-j} (t)$
are smooth densities on $M$ with values in $End (\Lambda^k (T^{\ast} M) \otimes
{\cal E}),$ depending on the parameter $t.$ Substituting into (2.33), and 
integrating with respect to $\lambda$ leads to 

\spf

\n {\bf (2.34)} \ $
\begin{array}{l}
\frac{d}{dt} \log {\cal R} (t) = \sum_k (-1)^k Tr (h A_{k;0} 
(t)) \\
= Tr (h (\sum_k (-1)^k A_{k;0} (t))) \\
= \dim \ {\cal E}_x Tr (h e(M,g)) \end{array}$

\spd

\n where $e (M,g)$ is the Euler form, which can be proven to be equal to 
$\frac{1}{\dim {\cal E}_x} \sum_k
\linebreak
(-1)^k A_{k;0} (t).$ If $\dim M$ is odd, 
$e(M,g) \equiv 0.$ \ \ \ \carre

\section{Proof of  Proposition 0.1 in the case $ g = g'$}

\subsection{Outline of the Proof}

\n In this section we prove Proposition 0.1 in the case where the Riemannian metric $g$ is 
the same as the metric $g^{\prime}$ of the generalized triangulation $\tau = 
(h, g^{\prime}).$ Throughtout this section we assume that the Hermitian structure $\mu$ is $\tau$-admissible, i.e. 
that there exists a neighborhood $U_h$ of the critical points of $h$ so that on 
$U_h$, $\mu$ is parallel with respect to the canonical flat connection on 
${\cal E}\to M$.  

\spd

From Theorem
2.9 we know that

\spd

\n {\em (A) $\log {\cal R} (t) = \log {\cal R} + \int^t_0 \frac{d}{dt} \log 
{\cal R} (t) dt$ is an affine function of $t,$}

\spd

\n hence $\log {\cal R} (t)$  has an asymptotic expansion with $\log {\cal R}$ 
as the free term. 
 
\spd

Recall from the 
proof of Proposition 2.1 (iii) that, for $t$ sufficiently large, 
\newline  $(H_0 (\Lambda (M; {\cal E})), d(t))$ is the direct sum of two 
sub\-com\-ple\-xes, 

\[(H_0 (\Lambda (M; {\cal E})), d(t)) = (\Omega_{sm} (t), d(t)) 
\oplus (H_{0, la} (t), d(t))\]

\n where  \[H_{k;0}:= H_0 (\Lambda^k (M; {\cal E})), \] 
\[\Omega^k_{sm} (t) = Q_k (t) (H_{k;0}),\]
$H_{k;0,la} (t) = (Id - Q_k 
(t)) (H_{k;0})$ and $Q_k(t)$ denotes the orthogonal spectral projector $Q_k(t)
: H_{k;0} \rightarrow H_{k;0}$ of the k-Laplacian $\Delta_k (t)$ corresponding
to the interval [0,1]. Notice that 
$\Omega^k_{sm} (t)$ consists of smooth forms  and that
$(H_{0,la} (t), d(t))$ is an algebraically
acyclic $\zeta$-regular complex of sF type. Notice also that
\newline $(\Omega_{sm} (t), d(t))$ is a
$\zeta$-regular complex of finite type for $t$ large enough.  The finite type 
property 
follows from [BFKM] Theorem 5.5. 
Denote by $g_{s,sm}(t)$ the
restriction of $g_s (t)$ to $\Omega_{sm}(t).$ This is a morphism of trace
class for $s > \frac{n}{2} + 1$ since  $g_s(t)$ is 
of trace class by Proposition 2.5. We have the 
following short exact sequence of algebraically acyclic, $\zeta$-regular 
complexes of sF type,

\spf

\n {\bf (3.1)} \hsps $0 \rightarrow {\cal C}(g_{s,sm}(t)) 
\stackrel{I}{\rightarrow} {\cal C} (g_s (t)) 
\stackrel{P}{\rightarrow} H_{0,la} (t) \rightarrow 0$

\spd

\n where $I$, resp. $P$, is the obvious inclusion, resp. 
projection. One verifies in a straightforward way that (3.1) satisfies the
assumptions in Lemma 1.14 and therefore, with 

\spf

\n {\bf (3.2)} \hsps ${\cal R}_{sm} (t) := T({\cal C} (g_{s, sm} (t))),$

\spd

\n we obtain the decomposition

\spd

\n {\em (B) $\log {\cal R} (t) = \log {\cal R}_{sm} (t) + \log T (
H_{0, la} (t)).$}

\spd

\n For consistency with the notation in [BFK2] and [BFKM] we write
\[\log T_{la} (t):=  \log T (H_{0, la} (t)).\]
\n The two terms, $\log {\cal R}_{sm} (t)$ and $\log T_{la} (t),$ will 
each be treated separately. To obtain an asymptotic expansion of 
$\log {\cal R}_{sm} (t)$ as $t \rightarrow \infty,$ we proceed in two steps. 
Recall that $f(t)$ is given by the composition 

\[f(t) : (\Omega_{sm} (t), d(t)) \stackrel{e^{th}}{\rightarrow} (\Omega, d)
\stackrel{Int}{\rightarrow} ({\cal C}, \delta).\]

\n We show in subsection 3.2, Proposition 3.1 (i), that

\[{\cal R}_{sm} (t) = T({\cal C} (f(t))).\]

\n The morphism $f(t),$ unlike $g_{s,sm} (t),$ can be studied using 
the Witten-Helffer-Sj\"ostrand theory (cf [BFK2, section 3], also [BFKM, section 5]). 
We show in subsection 3.2, Proposition  3.1 (ii) that $\log T({\cal C} (f(t)))$
has an asymptotic expansion which we calculate  up to  terms of order 
$0 ( \frac{1}{t}).$
Proposition 3.1   implies that:

\spd

\n {\em (C) $\log {\cal R}_{sm} (t)$ admits an asymptotic expansion for $t 
\rightarrow \infty$ of the form 

\[\begin{array}{lll}
\log {\cal R}_{sm} (t) & = & \dim {\cal E}_x \cdot (\sum_j (-1)^j \sum_{x \in Cr_{j}(h)} h(x))t \\
& - & (\sum_j (-1)^j (\frac{n}{4} - \frac{j}{2}) m_j \cdot \dim 
{\cal E}_x) \log(\frac{\pi}{t})
+ 0 ( \frac{1}{t}). \end{array}\]}




\n Concerning $\log T_{la}(t),$  we will establish in subsection 3.3:

\spd

\n {\em (D) (i) $\log T_{la} (t)$ admits an asymptotic expansion as $t \rightarrow
\infty$ and

\[FT_{t=\infty} (\log T_{la} (t)) = - FT_{t= \infty} (\log {\cal R}_{sm} (t)) +
\int_M  \alpha \]

\spd

\n where $\alpha$ is a local density on $M$  which vanishes in a neighborhood 
of the critical points of $h.$

\spd

(ii) If $\mu$ is parallel with respect to the canonical flat connection, 
$\int_M \alpha 
= 0.$ }

\spd

\n The results (A)-(D) prove Proposition 0.1 in the case where $g = g^{\prime}.$

\spd

\subsection{Asymptotic expansion of $\log {\cal R}_{sm} (t)$}

\spf

\n In this subsection we prove the statement (C) of subsection 3.1. By the
Witten-Helffer-Sj\"ostrand theory, $f(t) := 
Int  \cdot e^{th}$ is an isomorphism of
cochain complexes for $t$ sufficiently large. Therefore the mapping cone
${\cal C} (f(t))$ is an algebraically acyclic, $\zeta$-regular complex of
sF type and thus admits a torsion, $T({\cal C} (f(t))).$

\spd

\begin{prop} For $t$ sufficiently large, 

\begin{itemize}
\item[(i)] $\log {\cal R}_{sm} (t) = \log T ({\cal C} (f(t)));$

\item[(ii)] $\log T ({\cal C} (f(t)))$ admits an asymptotic expansion for 
$t \rightarrow \infty$ of the form
\end{itemize}

\[\begin{array}{lll}
\log T ({\cal C} (f(t))) & = & \dim {\cal E}_x \cdot (\sum_j (-1)^j 
\sum_{x \in Cr_{j}(h)} h(x)) t\\
& - & ( \sum_j (-1)^j (\frac{n}{4} - \frac{j}{2}) m_j \dim {\cal E}_x) \log
(\frac{\pi}{t}) + 0 (\frac{1}{t}).
\end{array}\]

\end{prop}

\n {\bf Proof} (i) Notice that the morphism $f(t)$ is equal to the composition
$(s > n+2)$

\spf

\n {\bf (3.3)} \hsps $(\Omega_{sm} (t), d(t)) \stackrel{e^{th}_s}{\rightarrow}
H_s (\Lambda (M; {\cal E})) \stackrel{Int_s}{\rightarrow} {\cal C} 
(M, \tau, {\cal T})$

\spd

\n whereas $g_{s, sm} (t)$ is given by

\spd

\n {\bf (3.4)} 
\hsps
$(\Omega_{sm} (t), d(t)) \stackrel{e^{th}}{\rightarrow} (H_0 
(\Lambda (M; {\cal E})), d) \stackrel{(\Delta + Id)^{-s/2}}{\rightarrow} 
(H_s (\Lambda
(M; {\cal E})), d)$

\[\stackrel{Int_s}{\rightarrow} {\cal C} (M, \tau, {\cal T}).\]

\n By Proposition 2.5, $Int_s$ is a morphism of trace class, multiplication by
$e^{th}$ is a morphism of trace class as $\Omega_{sm} (t)$ is a cochain
complex of finite rank (cf Examples 1.6) and $(\Delta + Id)^{-s/2}$ is an
isometry. Each of them induces an isomorphism in algebraic cohomology (cf 
Proposition 2.1). Therefore we can apply Proposition 1.15 to (3.3) and
(3.4) to obtain 

\spd

\n {\bf (3.5)} \hsps $\log T ({\cal C} (f (t))) = \log T ({\cal C} (e^{th}_s)) +
\log T ({\cal C} (Int_s))$

\spd

\n and

\spd

\n {\bf (3.6)} \ $\log T ({\cal C} (g_{s,sm} (t))) = \log T ({\cal C} (e^{th}))
 + \log T ({\cal C} (Int_s (\Delta + Id)^{-s/2} )).$

\spd

\n In view of Proposition 1.16 and the fact that $(\Delta + Id)^{-s/2}$ is
an isometry we conclude that

\spd

\n {\bf (3.7)} \hsps $\log T ({\cal C} (Int_s)) = \log T ({\cal C} 
(Int_s (\Delta + Id)^{-s/2})).$

\spd

\n Further, $e^{th} = In_{s,0} \ e^{th}_s$ is a morphism

\[(\Omega_{sm} (t), d(t)) \stackrel{e_s^{th}}{\longrightarrow} H_s (\Lambda
(M; {\cal E})) \stackrel{In_{s,0}}{\longrightarrow} H_0 (\Lambda (M; {\cal E}))\]

\n with $In_{s,0}$ and $e^{th}_s$ being both morphisms of trace class and 
inducing iso\-mor\-phisms in algebraic cohomology (cf Proposition 2.1). Thus, 
applying Proposition 1.15 once more, we obtain, in view of Proposition 2.3,

\spd

\n {\bf (3.8)} \hsps $\log T({\cal C} (e^{th})) = \log T 
({\cal C} (e_s^{th})).$

\spd

\n Combining (3.5) - (3.8) leads to the proof of statement (i).

\spd

\n (ii) Recall from [BFK2, section 3.3] the following commutative diagram (for
$t$ sufficiently large)

\spd

\n {\bf(3.9)} \hsps $\begin{array}{ccl}
(\Omega_{sm} (t), d(t)) & \stackrel{f(t)}{\longrightarrow}
& ({\cal C}, \delta) \\
&& \\
\Phi+ 0 (\frac{1}{t}) \searrow && \ \downarrow S(t) \\
&& \\
&& ({\cal C}, \tilde \delta) \end{array}$

\spd

\n where $\Phi$ is a morphism which is an isometry,
 $\tilde \delta = \tilde \delta (t)$ is chosen so that the diagram above is commutative,
$\tilde \delta_k (t) := S_{k+1}(t)^{-1} \cdot \delta_k \cdot S_k(t),$
and $S(t)$ is the
morphism given by 
$$
S_k(t) = (\bigoplus_{x \in Cr_{k}(h)}  S_x(t)) :
{\cal C}^k  \longrightarrow  {\cal C}^k
$$
where ${\cal C}^k := \bigoplus_{x \in Cr_{k}(h)} {\cal E}_x $ 
and $S_x(t):{\cal E}_x\to {\cal E}_x$ is given by 

\n $S_x (t) := e^{-th(x)} (\frac{\pi}{t})^{\frac{n}{4} - \frac{k}{2}}Id.$

\spd

\n For $t$ sufficiently large, $S(t), f(t), \Phi + 0 (\frac{1}{t})$ are 
iso\-mor\-phisms of cochain com\-ple\-xes. Thus, by Lemma 1.11, the 
mapping cones
${\cal C} (\Phi + 0 (\frac{1}{t})), {\cal C}(S(t))$ and ${\cal C} (f(t))$ are
algebraically acyclic. By Proposition 1.15,

\[\log T ({\cal C} (\Phi + 0 (\frac{1}{t}))) =  \log T({\cal C} (f(t))) +
\log T ({\cal C} (S(t))). \]

\spd

\n Using Proposition 1.12, one verifies that $\log T ({\cal C} (Id + 
0 (\frac{1}{t} )))
= 0 (\frac{1}{t}).$

\spd

\n Therefore 

\spd

\n {\bf (3.10)} \hsps $\log T ({\cal C} (f (t))) = - \log T ({\cal C} (S (t))) +
0 (\frac{1}{t}) .$

\spd

\n By Proposition 1.12,

\[\begin{array}{lll}
\log T ({\cal C}(S(t))) 
& = & \sum (-1)^j  \sum_{ x \in Cr_{j} (h)}
\log (e^{-th(x)} (\frac{\pi}{t})^{\frac{n}{4} - \frac{j}{2}})
\dim {\cal E}_x.
\end{array}\]

\spd

 \n In view of (3.10) we obtain

\[\begin{array}{lll}
\log {\cal C}(f(t)) & = & (\sum (-1)^{j} \sum_{ x \in Cr_{j} (h)} h(x)
\dim {\cal E}_x) t \\
&& - (\sum (-1)^j (\frac{n}{4} - \frac{j}{2}) m_j \dim {\cal E}_x) 
\log \frac{\pi}{t}
+ 0 (\frac{1}{t}). \ \ \ \mbox{\carre} \end{array}\]

\spd

\subsection{Asymptotic expansion of $\log T_{la} (t)$}

\spf

\n In this section we prove the statement (D) in subsection 3.1. As 
$\log T_{la}(t)= \log {\cal R}(t) -\log {\cal R}_{sm}(t)$
and in view of Theorem 2.9 and Proposition 3.1, $\log T_{la}(t)$
admits an asymptotic expansion for $t\to \infty.$

The arguments
to prove (D) follow closely  the ones in [BFKM, section 6.2]: we will derive 
statement (D)
from a relative version of this statement , Proposition 3.2, and from the validity of  
(D) in some simple cases. 

\spf

We consider systems 
${\cal F} := (M, {\cal E}, \nabla, \mu,  g,\tau )$, 
where $(\cal E,\nabla)$ is a bundle of $\cal A$ -Hilbert modules equipped 
with a flat connection, $\mu$ a Hermitian structure, $g$ a Riemannian metric 
and $\tau=(h,g')$ a generalized triangulation with the  metric $g'= g.$
We suppose that $\mu$ is admissible with respect to $\tau,$
i.e there exists an open neighborhood $U_h$ of the critical points
so that on $U_h, \ \mu$ is parallel with respect to the connection $\nabla.$ 
Introduce  
${V}_{\cal F}(t,\epsilon ):= 
\frac{1}{2} \sum_q (-1)^{q+1} q \log\det ({\Delta}_q (t) + \epsilon).$ 

\spd

\begin{prop} For any $\epsilon>0 $ there exists a smooth 
density $\alpha_{\cal F}(\epsilon) $
on  
\newline $M\setminus Cr (h),$ which is polynomial in $\epsilon$ and a local quantity
(cf subsection 0.2 and [BFKM] section 2) so that the following statement hold:

\spd

\n (i) If ${\cal F'}:= (M, {\cal E}, \nabla, \mu,  g,\tau_{\cal D} )$
 denotes the system obtained from
${\cal F}$  by replacing the triangulation $\tau = (h, g)$
by its dual $\tau_D = (n-h, g)$  ($n= \dim M$), then 

\spd

\n {\bf (3.11)} \hsps
$\alpha_{{\cal F}}(\epsilon) + (-1)^{n+1} \alpha_{{\cal F}'}(\epsilon) = 0.$

\spd

\n(ii) If ${\cal F}$ and $\tilde{\cal F} $ are two systems as above and
${\cal F}\otimes\tilde{\cal F}$ denotes the system defined by 
\[{\cal F}\otimes\tilde{\cal F}:= (M\times {\tilde{M}}, {\cal E}\otimes 
\tilde{\cal E},
\overline{\nabla}, \mu\otimes\tilde{\mu},  g\times \tilde{g}, 
\tau \times \tilde{\tau})\]
with $\overline{\nabla}=\nabla\otimes id + id\otimes\tilde{\nabla},$
then 

\spd

\n {\bf (3.12)} \hsps
$\alpha_{{\cal F}\otimes\tilde{\cal F}} (\epsilon)=( 
{\alpha}_{\cal F}(\epsilon)\otimes e(\tilde{M},\tilde{g}) + 
e(M,g)\otimes \tilde{\alpha}_{\tilde{\cal F}}(\epsilon))\dim {\cal W},$

\spd

\n where $e(M,g)$ denotes the Euler form of $(M, g).$ 

\spd

\n (iii) The density $\alpha_{\cal F}(\epsilon)$ vanishes on $U_h\setminus Cr(h).$ 

\spd

\n (iv) Assume ${\cal F}$ and $\tilde{\cal F}$ are two systems as above
so that  $\sharp {\cal C}r_q (h) = \sharp Cr_q (
\tilde{h}) \ (0 \le q \le n)$,
and the bundles $\cal E$ and $\tilde{\cal E}$ have isomorphic 
$\cal A-$Hilbert modules as fibers.
Then $V_{\cal F} (t,\epsilon) - V_{\tilde{\cal F}} (t,\epsilon)$
has an asymptotic expansion with 

\spd

\n (a) $FT_{t=\infty} (V_{\cal F}(t,\epsilon) - V_{\tilde{\cal F}}(t,\epsilon))=
\int_{M \setminus Cr(h)} \alpha_{\cal F}(\epsilon) - \int_{\tilde{M}\setminus 
Cr(\tilde{h})} \alpha_{\tilde{\cal F}}(\epsilon);$

\spd

\n (b) $ FT_{t=\infty} (\log T_{la} (t)) - FT_{t=\infty} (\log \tilde{T}_{la} (t))
=\int_{M \setminus Cr(h)} \alpha_{\cal F}(0) - \int_{\tilde{M} \setminus 
Cr(\tilde{h})} \alpha_{\tilde{\cal F}}(0).$
\end{prop}

\spd

\n {\bf Proof:} We begin with the construction of the density 
$\alpha_{\cal F}(\epsilon)$. Away from the critical  points of $h$
we choose a coordinate system for ${\cal E} \to M.$ In these coordinates 
we calculate inductively the 
quantities $r^q_{-2-j}(h,\epsilon; x,\xi ,t,\mu ),$  by the formulae 6.59 in 
[BFKM],
from the symbol of the operator
with parameter $\Delta_q(t)+\epsilon$, which is elliptic with parameter $t$
away from $Cr(h)$
(cf [BFKM] (3.2)). We are  interested in the quantity $r^q_{-2-n} $
only. We use  formula (3.4) in [BFKM]  to obtain from $r^q_{-2-n}$ the 
density $\alpha_{\cal F}^q(h,\epsilon)$  on $M\setminus Cr(h)$. We define 

\spd

\n {\bf (3.13)} \hsps $\alpha_{\cal F}(h,\epsilon):=
1/2  \sum_q (-1)^{q+1}q \alpha_{\cal F}^q(h,\epsilon).$

\spd

\n Since by construction  $r^q_{-2-n}$ is a polynomial in $\epsilon$ 
of degree smaller 
than $n$  so is $\alpha_{\cal F}(h,\epsilon).$

\spd

(i) follows from the following homogeneity property (cf [BFKM] (6.60)):
\[
r_{-2-n}(h,x,-\xi ,-t,\mu )=(-1)^n r_{-2-n}(h,x,\xi ,t,\mu )\
\]
and by a straightforward verification 
\[
r_{-2-n}(n-h,x,\xi ,t,\mu )= r_{-2-n}(h,x,\xi ,-t,\mu )\,.
\]

\spd

To prove (ii) we  consider 
systems  ${\cal G} = (M, {\cal E}, \nabla, \mu,  g,\omega )$,
where $M, {\cal E}, \nabla, \mu, g $ are as above and $\omega$  is a 
closed 1-form on $M.$
A system ${\cal F}= (M, {\cal E}, \nabla, \mu, g, \tau=(g, h ))$ gives 
rise to a 
system ${\cal G}$ by  
taking $\omega= dh.$
For any such  ${\cal G}$, define the Witten deformation by 
taking $\omega$
instead of $dh$ and then the Witten Laplacians
$\Delta_q(t)$ for any  real number $t.$   If the 1-form $\omega$ 
has no zeros,
$\Delta_q(t) $ is elliptic with parameter and the general theory in 
[BFKM], section 2
implies that for any fixed $\epsilon >0,\
{V}_{\cal G}(t,\epsilon ):= 
\frac{1}{2} \sum_q (-1)^{q+1} q \log ({\Delta}_q (t) + \epsilon)$
has an asymptotic expansion for $t\to \infty,$ whose free term  is
$\int_M \alpha_{\cal G}$, with $\alpha_{\cal G}$ a local density 
on $M.$
Following [BMFK] section 2 this density is  calculated in  the same 
way as the density $\alpha_{\cal F}$ (using $\omega$ instead of
$dh).$

\n Since ${\cal F}|_{M\setminus Cr(h)}$ is locally isomorphic to the 
restriction to 
an open set of a system 
${\cal G}= (\tilde {M}, \tilde{\cal E},  \tilde{\nabla},  \tilde{\mu},  
\tilde {g}, 
\tilde{\omega} )$ such that $ \tilde M$ is closed and $\tilde{\omega}$  
has no zeros,
it suffices to check
(3.11) for $\alpha_{\cal G}$ and  $\alpha_{\tilde{\cal G}}$ instead of
$\alpha_{\cal F}$ and  $\alpha_{\tilde{\cal F}}.$

For a system ${\cal G}$ and $u>0, t>0$
consider the operator
$e^{-u\Delta_q(t)}$ which is a smoothing operator and denote by 
$\lambda_q(u,t)$
the pointwise (von Neumann)trace
of the restriction of its Schwartz kernel to the diagonal. This 
provides a two parameter 
family (in $u$ and $t$) of densities on $M.$ 
Denote by \[\lambda_{\cal G} (u,t ):=1/2\sum_q (-1)^{q+1}q\lambda_q(u,t).\]  

\n Define the smooth three parameter family, $\eta(s,t, \epsilon),$ 
 of densities on $M$  for ${\Re}e s$ sufficiently large,  

\spd

\n {\bf (3.14)} \hsps $  \eta(s,t, \epsilon) = 
\frac{1}{\Gamma(s)} \int_0^{\infty} u^{s-1}e^{-u\epsilon}
\lambda_{\cal G}(u,t)du. $

\spd

\n which, by analytic continuation, is  holomorphic in $s$ near $0\in \bbc$. 
Denote by $\theta_{\cal G}(t,\epsilon)$ the density
valued function in $t$ and $\epsilon,$

\[\theta _{\cal G}(t,\epsilon):=\frac{d}{ds} \vert _{s=0}\eta(s,t,\epsilon).\]

\n For the system ${\cal G}\otimes \tilde{\cal G},$ we can prove 
that 

\spd

\n {\bf (3.15)} \ $  \theta_{{\cal G}\otimes \tilde{\cal G}}(t,\epsilon)=
\theta_{\cal G}(t,\epsilon)\otimes e(\tilde M,\tilde g)\dim {\cal W}+ 
 e(M,g)\otimes \theta_{\tilde{\cal G}}(t,\epsilon) \dim{\cal W}.$

\spd

\n This can be derived as in the proof of Proposition 2.4 in [BFKM] from 
the following facts:

\n (i)  $\Delta^{{\cal G}\otimes \tilde{\cal G}}_q(t)=
\Delta^{\cal G}_q\otimes Id + Id\otimes \Delta^{\tilde{\cal G}}_q(t),$

\n whose proof is similar to  the proof of Proposition 2.4 in [BFKM] and for 
any $t,$

\n (ii)   $\eta(0,t,\epsilon=0) = e(M,g) \dim {\cal W},$

\n whose proof follows from  the local index theorem 
 of Bismut and Zhang. This takes care of (3.12). 

\spd

\n For (iii), we use the locality of the density 
$\alpha_{\cal F}(\epsilon)$ and the explicit formula of 
$\Delta_q(t).$  In admissible coordinates in a 
neighborhood of the critical points
$\Delta_q(t),$ is the same as $\Delta_q(t),$ for a product of systems ${\cal G}$ 
with underlying
manifolds of dimension one. Using  (3.12) and the vanishing of the Euler 
form on $M= \bbr ,$  the result follows. 

\spd

\n (iv) (a) is actually Proposition 6.6 in [BFKM] and (b) follows from Lemma 6.7 
of [BFKM].
\ \ \ \carre

\spd

\n {\bf Proof of statement }(D).  We want to apply Proposition 3.2. Set 
$\tilde{M} := M, \ \tilde{\tau} := \tau, g= \tilde g,$  $\tilde{\rho}$ the 
trivial
representation of $\Gamma \equiv \pi_1 (M)$ on ${\cal W}.$ Then the 
associated bundle 
$\tilde{{\cal E}} \rightarrow M$ is trivial ; we choose  $\tilde{\mu}$ the 
trivial Hermitian 
structure.  Then by (2.15) and the equality of analytic and Reidemeister 
torsion for the trivial
representation we have  $\tilde{\cal R}=1.$ Since $\log \tilde{\cal R}(t)$
admits an asymptotic expansion by $(B)$ so does $\log \tilde{T}_{la} (t)$.
Since the free term of 
the ex\-pan\-sion  of $\log \tilde{\cal R}$ is zero, one obtains   

\spd

\n {\bf (3.16)} \hsps $FT_{t= \infty} \log \tilde{T}_{la} (t) = 
\Big(\sum_j (-1)^j (\frac{n}{4} - 
\frac{j}{2}) m_j
\dim {\cal E}_x \Big) \log \pi.$

\spd

\n Statement (i) in (D) follows then from  Proposition 3.2  (iv) . To prove 
statement (ii), 
observe that since
$\mu$ is parallel with respect to the canonical flat  connection $\nabla$, the
Hodge $\ast$ operator
induces an isometry between $L_2 (\Lambda^k 
(M; {\cal E}))$ and $L_2 
(\Lambda^{n-k} (M; {\cal E}))$
and conjugates $\Delta_k (t)$ with $\Delta_{n-k} (-t).$
Thus $\Delta_k (t)$ and $\Delta_{n-k} (t)$ are isospectral and then by 
the definition of 
$V_{\cal F}(t,\epsilon )$ we have :

\[V_{\cal F}(t,\epsilon )= (-1)^{n+1} V_{\cal F'}(t,\epsilon )\]
and 
\[V_{\tilde{\cal F}}(t,\epsilon )= (-1)^{n+1} V_{\tilde{\cal F}'}(t,\epsilon ).\]

\n Therefore 

\[V_{\cal F}(t,\epsilon ) - V_{\tilde{\cal F}}(t,\epsilon )=
(-1)^{n+1} V_{\cal F'} (t,\epsilon )- (-1)^{n+1} V_{\tilde{\cal F}'}(t,\epsilon ),\]
and thus by Proposition 3.2 (iv) we have  

\[FT_{t=\infty}\log T_{la}({\cal F}, t) -
FT_{t=\infty}\log T_{la}(\tilde{\cal F}, t)=\]
\[=(-1)^{n+1} FT_{t=\infty}\log T_{la}({\cal F'}, t) - (-1)^{n+1} 
FT_{t=\infty}\log T_{la}(\tilde{\cal F}', t).\]
Again by Proposition 3.2 (iv) 

\[\int_{M\setminus Cr(h)}\alpha_{\cal F}(0)-\int_{M\setminus Cr(h)}
\alpha_{\tilde{\cal F}}(0)=\]
\[(-1)^{n+1}(\int_{M\setminus Cr(h)}\alpha_{\cal F'}(0)-
\int_{M\setminus Cr(h)}\alpha_{\tilde{\cal F}'}(0)).\]
In view of (3.12), (D) (i) and 
\[\int_{M\setminus Cr(h)}\alpha_{\tilde{\cal F}}(0)= 
\int_{M\setminus Cr(h)}\alpha_{\tilde{\cal F}'}(0)=0, \]
because $\mu$ is parallel, one obtains  
\[\int_{M\setminus Cr(h)}\alpha_{\cal F}(0)=(-1)^{n+1}\int_{M\setminus Cr(h)}
\alpha_{\cal F'}(0) .\]
Statement (D)(ii)  follows by combining this last formula with Proposition  3.2 (iv)
 \ \ \ \carre

\spd

\section{ Anomalies of the relative torsion}

\subsection{ Metric and Hermitian anomalies}

\n In this subsection, we investigate how the relative torsion ${\cal R}=
{\cal R}(M,\rho , \mu, g,\tau )$ varies with respect to the Riemannian 
metric $g$  and the Hermitian structure $\mu$. We denote by ${\cal E}\to M, \nabla$ 
the flat bundle associated to $\rho.$
\n In order to formulate the results we  introduce a number of 
additional quantities:

\spd

\n {\it The closed 1-form $\theta = \theta (\rho, \mu) \in \Omega^1 (M)$}: 
Choose a finite covering $M$ by open sets $(U_j)_{j\in J}$, which are 
simply connected, and points $x_j\in U_j$.  Define $v_j\colon U_j\to\bbr$ 
by

    \[ v_j(x):=\log{vol} (T_{x,x_j})={1 \over 2}\log\det 
(T^\ast_{x,x_j}T_{x,x_j})\]

\n where $T_{x,x_j}\colon ({\cal E}_x,\mu_x)
\to ({\cal E}_{x_j},\mu_{x_j})$ 
denotes the parallel transport from ${\cal E}_x$ to 
${\cal E}_{x_j}$ along 
any curve joining $x$ and $x_j$ inside $U_j$.  (As the connection 
$\nabla$ 
of ${\cal E}\to M$ is flat, the map $T_{x,x_j}$ does not depend on the 
choice of the curve.)  Denote by $\theta =\theta (\rho, \mu )$ the 
smooth $1$-form 
on $M$ defined by $\theta_j :=dv_j.$   Notice that if the 
Hermitian structure 
$\mu$ is parallel with respect to the canonical flat connection 
$\nabla,$ then  
$\theta(\mu) =0.$ We mention that in [BZ] the quantity 
$\theta$ is defined to be $2\theta(\rho,\mu)$.

If the representation $\rho$ is unimodular,
(i.e. $\log vol(\rho (g))= 0$ for any $g \in \Gamma$) then one can find 
$\mu$ so that $\theta (\rho, \mu)= 0.$ To see this  let 
$\det_{\bbr}{\cal E}\to M$
denote the 1-dimensional real vector bundle obtained by assigning 
to each point $x\in M$
the vector space $\det_{\bbr}{\cal E}_x,$  as described in [CFM], 
and let $\det \nabla$ denote the induced flat connection in 
$\det_{\bbr}{\cal E}\to M$.
This is the flat bundle associated with the representation $\det_{\bbr} \rho$.  
The unimodularity 
of $\rho$ implies the existence of a parallel section $s_0$ in 
$\det_{\bbr}{\cal E}\to M.$
Take a  Hermitian structure $\mu$ in ${\cal E}\to M,$ and denote by 
$s$ the tautological 
section in $\det_{\bbr}{\cal E}\to M$ induced by $\mu.$
Clearly there exists a smooth nonzero function $f:M\to \bbr_+$ so 
that $s_0= f\cdot s.$
If $\mu$ is parallel above the open set $U,$ then $f |_U=1.$ 
It is immediate from definition of $\theta$ that $\theta (\rho, f\cdot\mu)= 0.$ 

If $(\rho_i,\ \mu_i)$, $i=1,2$ are two pairs consisting of a representation and 
a Hermitian structure on 
the induced flat bundle, then one can verify in a straightforward manner
that

\spd

\n {\bf (4.1)}\hsps $\theta (\rho_1\otimes \rho_2, \mu_1\otimes \mu_2)= 
\theta (\rho_1, \mu_1)\otimes 1 +
1\otimes \theta (\rho_2, \mu_2).$

\spd

\n {\it The function $V=V(\rho, \mu_1,\mu_2) \in \Omega^0(M)$}: 
If $\mu_j,  j=1,2$ are two
Hermitian structures of ${\cal E}\to M$, define the smooth function 

\spd

\[V(x):= V(\rho, \mu_1,\mu_2)(x)=
\log {vol} ({Id}_x: ({\cal E}_x,\mu_1(x))\to ({\cal E}_x,\mu_2(x))).\]

\spd

Note that

\spd

\n {\bf (4.1')}\hsps $\theta (\rho,\mu_1) - \theta (\rho,\mu_2)= 
d V(\rho, \mu_1, \mu_2),$

\spd

\n {\bf (4.1'')}\hsps $V(\rho,\mu_1, \mu_2)= -V(\rho,\mu_2, \mu_1),$

\spd

\n {\bf (4.1''')}\hsps $V(\rho,\mu_1, \mu_3)= V(\rho,\mu_1, \mu_2)+ 
V(\rho,\mu_2, \mu_3).$

\spd

\n {\it The Euler form  $e(M,g) \in \Omega^n(M; {\cal O}_M$)}: 
Denote by $ e(M,g)$ the Euler 
form associated to the 
Levi-Civit\`a connection on the tangent bundle $TM.$ This is a form with
coefficients in the orientation bundle of $M,$ a flat, 1-dimensional real vector bundle.

\spd

\n {\it  The Chern-Simon element $[e_{CS}(M,g_1,g_2)]\in\Omega^{n-1}
(M; {\cal O}_M)/d\Omega^{n-2}(M; {\cal O}_M)$}: For two Riemannian metrics $g_1$ and $g_2$ on $M$
denote by $[e_{CS}=e_{CS}(M,g_1,g_2)]$
the Chern-Simon class , cf [BZ], p 46.
Recall that 
\[d(e_{CS}(M,g_1,g_2))=e(M,g_2)-e(M,g_1)\]
and that there is 
a canonical construction, due to Chern-Simon, for a representative $e_{CS}$ 
of $[e_{CS}]$ so that for a smooth $1$-parameter family $g_2 (u)$ of 
Riemannian 
metrics on $M$, $e_{CS}(g_1,g_2(u))$ is a smooth $1$-parameter 
family of $(n-1)$ forms.

\spd

\n The following object appears in Theorem 0.1 and 
is related to the 
previous quantities.

\spd

\n {\it  The form $\Psi (TM,g)\in \Omega^{n-1}(TM\setminus M; {\cal O}_{TM}):$}  Let 
$\pi: TM\to M$ be 
the tangent bundle of $M$. Clearly, ${\cal O}_{TM} = \pi^\ast {\cal O}_M.$
 Following Mathai-Quillen [MQ] Theorem 6.4, or Bismut Zhang [BZ] 
 Theorem 3.4,
 one can construct a smooth form $\Psi (TM,g) \in 
 \Omega^{n-1}(TM \setminus M; {\cal O}_{TM})$
 with the following properties:

\spd

\n {\bf (4.2')}\hsps $ d \Psi (TM,g) =\pi^\ast (e(M,g)) $

\spd 

If $g_1$ and $g_2$ are two Riemannian metrics on $M$ then 

\spd

\n {\bf (4.2'')}\hsps $ \Psi (TM,g_2)-  \Psi (TM,g_1) = \pi^\ast e_{CS}(M, g_1,g_2)$ (modulo exact forms)

\spd

If $\lambda: TM\setminus M\to TM\setminus M$ denotes multiplication by 
the real number $\lambda\ne 0$ then 

\spd

\n{\bf (4.2''')}  \hsps $\lambda^{\ast} (\Psi (TM,g))= 
(\lambda/{|\lambda|})^{n} \Psi (TM,g). $

\spd

\n {$\bf (4.2^{iv}$)} For $x \in M$ fixed, denote by $\omega$ the volume form on $T_x M$ and by $Y$ the
Euler vectorfield on $T_x M$ ( in polar coordinates, $Y = r \frac{\partial}{\partial r},$  $r$ the distance to the origin).
Then 

\spd

$\Psi (TM,g) \vert _{T_{x}M} = - \frac {\Gamma(n/2)}{2\pi^{n/2}} \frac {\iota_{Y}\omega}{|Y|^n}.
$

\spd

It is easy to check (cf [BZ] ) that for 
two Riemannian manifolds $(M_i, g_i), \ i=1,2$ we have  

\spf

\n {$\bf (4.2^{v}$)} \ $ \begin{array}{l} \Psi (T(M_1\times M_2), g_1\times g_2)= 
\Psi (T(M_1),g_1)
\otimes e(M_2,g_2) +\\
+ e(M_1,g_1)\otimes \Psi (T(M_2),g_2).
\end{array}$

\spd

\begin{prop} (Metric anomaly of the relative torsion)
\newline
(i)  $\log{\cal R} (g_2)-\log{\cal R}(g_1)=-\int_M \theta\wedge e_{CS} 
(g_1,g_2)$.\newline
(ii) If $\dim M$ is odd or $\mu$ is a Hermitian structure parallel with respect 
to the canonical flat connection, the relative torsion ${\cal R}$ is 
independent of $g$.
\end{prop}

\begin{prop}(Hermitian anomaly of the relative torsion)
\[\log{\cal R}(\mu_2) - \log{\cal R}(\mu_1)=
-\int_M V(\rho, \mu_1,\mu_2)e(M,g) +\]
\[+\sum_{x\in Cr (h)} (-1)^{\rm{ind}\ x}V(\rho,\mu_1, \mu_2)(x).\]
\end{prop}

\n In the case ${\cal A}= \bbc$ these results were  established in
 by Bismut Zhang, cf [BZ].  

\n Propositions 4.1 and 4.2  will  be proved at the 
 same time. 
 Their  proof is reduced 
  to some  local index  type results established in [BZ].

 \spd

 \n {\bf Remark} Theorem 2.1. can be also derived as a special case 
 of  a  (slightly more general) version of  Proposition 4.2(i). 

\spd

\n {\bf Proof of Propositions 4.1, 4.2:}  (ii) follows from (i) by noting that 
$e(M,g) = 0$  and  $e_{CS}(M,g_1,g_2)=0 $ if ${dim}M$ is 
odd and  that $\theta(\mu)=0$ if  the Hermitian structure $\mu $ is parallel with
respect to the canonical flat connection.

\spd

\n To prove (i), consider a smooth $1$-parameter family $g(u)$ of 
Riemannian metrics and a smooth $1$-parameter family $\mu(u)$ of
Hermitian structures, ($-1<u<+1$).  We want to compute ${d \over du}
\log{\cal R}(g(u),\mu_0)$
and ${d \over du}\log{\cal R}(g_0,\mu(u)).$
We begin by analyzing  ${d \over du}\log{\cal R}(g(u),\mu(u))$.  

\spd

\n Denote by 
$\langle\langle \cdot ,\cdot\rangle\rangle (u)$ the scalar product defined 
by $g(u)$ and $\mu(u)$ on $\Omega^k(M;{\cal E}))$ (cf. 2.2) and by $\langle 
\langle \cdot , \cdot \rangle\rangle_s(u)$ the one given by (2.4).
Clearly 
$\langle\langle \cdot ,\cdot \rangle\rangle (u)=\langle 
\langle \cdot, \cdot \rangle\rangle_0(u).$

Denote  by  $\Delta_k(u)\colon\Omega^k(M;{\cal E})\to\Omega^k(M;{\cal E})$
the Laplacian on $\Omega^k(M;{\cal E})$ induced by $g(u)$ and $\mu(u)$
and by  $H_{s}(u) (\Lambda (M;{\cal E}))$  the completion of 
$\Omega (M;{\cal E})$ with respect to $\langle\langle\cdot,\cdot\rangle
\rangle_{s}(u).$  

\spd

Consider the following commutative diagram

\[
\begin{array}{ccc}
&\gamma_{s+s'}(u) &\\
H_{0}(u) (\Lambda (M;{\cal E})) &\longrightarrow &{\cal C}(\tau)(u)\\
&&\\
\downarrow\varphi_{s,s'}(u) &&\downarrow \iota(u)\\
&\gamma_s(0)&\\
H_{0}(0) (\Lambda (M;{\cal E}))&\longrightarrow &{\cal C}(\tau)(0)\end{array}
\]

where 

\n $\gamma_s(u):= \textrm{Int}_s(\textrm{Id}+\Delta (u))^{-s/2},\  \varphi_{s,
s'}(u):=(\textrm{Id} +\Delta (0))^{+{s/2}}(\textrm{Id}+\Delta (u))^{-{s+
s' \over 2}},$

\n and $\iota(u)$ is the identity map on ${\cal C}(\tau)$ regarded as an 
isomorphism of ${\cal C}(\tau)$ with the scalar product induced from 
$\mu(u)$ and ${\cal C}(\tau)$ with the scalar product induced by 
$\mu(0).$

\n We recall from subsection 2.1 that, for $s$ and $s'$ sufficiently large, 
$\gamma_s(u)$ and $\varphi_{s,s'}(u)$ are morphisms of trace class.  Since 
$\gamma_{s+s'}(u)$ and $\gamma_s(0)$ induce isomorphisms in algebraic cohomology 
(cf. Proposition 2.5), so does $\varphi_{s,s'}(u)$.

\n We thus can apply Proposition 1.15 to obtain

\spd

\n {\bf (4.3)} \hsps $ \log T({\cal C}(\iota (u))+ \log{\cal R}(g(u))=\log{\cal R}(g(0))+ \log T({\cal C} 
(\varphi_{s,s'}(u))$

\spd

\n where ${\cal C}(\varphi_{s,s'}(u))$ denotes the mapping cone associated to 
$\varphi_{s,s'}(u)$.

\n Since $(\textrm{Id}+\Delta (u))^{s/2}\colon H_{s}(u)
(\Lambda (M;{\cal E}))\to
H_{0}(u)(\Lambda (M;{\cal E}))$ is an isometry, we conclude from 
Proposition 1.16 that

\spd

\n {\bf (4.4)} \hsps $T({\cal  C}(\varphi_{s,s'}(u))=
T({\cal C}(\textrm{In}_{s+s';s}
(u))$

\spd

\n where
\[
\textrm{In}_{s+s',s}(u)\colon H_{s+s'}(u)(\Lambda (M;{\cal E}))\to H_{s}(0)(
\Lambda (M;{\cal E}))
\]
denotes the canonical inclusion.  
To analyze the torsion of the mapping cone 
\newline $({\cal C}_u,D):= ({\cal C}(\textrm{In}_{s+s';s}(u)),
d(\textrm{In}_{s+s',s}(u)))$, 
note that

\begin{eqnarray*}
{\bf (4.5)}~~~~~~~~~~~~~~{{\cal C}_u}_k &:= &H_{s}(0)
(\Lambda^{k-1}(M;{\cal E}))\oplus 
H_{s+s'}(u)(\Lambda^k(M;{\cal E}))\,;\\
D_k &= &\left(\begin{array}{cc}
-d_{k-1} &\textrm{In}_{s+s',s}(u)\\
0 &d_k
\end{array}\right)
\end{eqnarray*}

\n with inner product $(\omega_1,\omega_1'\in\Omega^{k-1},
\omega_2,\omega_2'\in
\Omega^k)$

\spd

\n {\bf (4.6)}\hsps $\langle\langle (\omega_1,\omega_2),
(\omega_1',\omega_2')
\rangle\rangle :=\langle\langle\omega_1,\omega_1'\rangle\rangle_{s}(0)+
\langle
\langle \omega_2,\omega_2'\rangle\rangle_{s+s'}(u)$.
\spd

\n In order to calculate ${d \over du}\log T({\cal C}(u))$, introduce the zeroth 
order differential operators $A_k(0)\colon \Omega^k(M;{\cal E})\to\Omega^k 
(M;{\cal E})$ defined by 
\[
\langle\langle\omega ,\omega '\rangle\rangle_{0}(u)=\langle\langle A_k(u)
\omega ,\omega '\rangle\rangle_{0}(0)
\]
and consider the zeroth order pseudo-differential operator
\[
B_k(s+s';u)\colon H_{s+s'}(0)(\Lambda^k(M;{\cal E}))\to H_{s+s'}(0)(\Lambda^k 
(M;{\cal E}))
\]
defined by
\[
B_k(u)\equiv B_k(s+s',u):= (\textrm{Id} +\Delta_k (u))^{-(s+s')}A_k(u)^{-1}
{d \over du} (A_k(u)(\textrm{Id}+\Delta_k(u))^{s+s'})
\]
so that the following identity holds
\[
{d \over du} \langle\langle\omega ,\omega '\rangle\rangle_{s+s'}(u) = 
\langle\langle B_k(u)\omega ,\omega '\rangle\rangle_{s+s'}(0)\,.
\]

\n Let ${\cal B}_k(s+s',u)$ be the bounded operator in 
\newline ${\cal L}(H_{s}(0)(
\Lambda^k(M;{\cal E}))\oplus H_{s+s'}(0)(\Lambda^k(M;{\cal E})))$ defined by
\[
{\cal B}_k(u)\equiv {\cal B}_k(s+s',u):=\left(\begin{array}{cc}
0 &0\\
0 &B_k(s+s',u)\end{array}\right)\,.
\]
Then we have, in view of (4.6),
\[
{d \over du}\langle\langle (\omega_1,\omega_2),(\omega_1',\omega_2')\rangle
\rangle = \langle\langle{\cal B}_k(u)(\omega_1,\omega_2),(\omega_1',\omega_2')
\rangle\rangle\,.
\]

\n Proceeding as in [BFKM, Appendix 3], one obtains (cf. [BFKM, A3.10]) in 
view of the fact that the complex ${\cal C}_u$ is algebraically acyclic,

\spd

\n {\bf (4.7)}\hsps ${d \over du}\log T({\cal C}_u) = -Fp_{z=0}{1 \over \Gamma (z)}
\int^1_0 x^{z-1} {1 \over 2}\sum\limits^{n}_{q=0} (-1)^q Tr ({\cal B}_q(u)e^{-x
\tilde{\Delta}_q(u)})dx$

\spd

\n where $\tilde{\Delta}_q(u)$ denotes the $q$-Laplacian of ${\cal C}_u$ with 
respect to the inner product (4.6).  

\spd

\n By formula (1.23),
\[
\tilde{\Delta}_q(u)= \left(\begin{array}{cc}
\Delta_{q-1}(0)+\psi_{q-1}(u) &\eta_{q}(u)\\
\eta_q(u)^\ast &\Delta_q(u)+\psi_q(u)\end{array}\right)
\]
where $\left(\begin{array}{cc}
\psi_{q-1}(u) &\eta_q(u)\\
\eta_q(u)^\ast &\psi_q(u)\end{array}\right)$ is a nonnegative selfadjoint 
operator of trace class ($s,s'$ sufficiently large).  

\spd

By the remark after 
Proposition 1.5 we conclude, in view of definitions (1.5) and (1.7),
\begin{eqnarray*}
{\bf (4.8)}~~~~~ &- &{1 \over 2} Fp_{z=0}{1 \over \Gamma (z)}\int^1_0 
x^{z-1}\sum^n_{q=0} (-1)^q Tr({\cal B}_q(u)e^{-x\tilde{\Delta}_q(u)})dx\\
= &- &{1 \over 2} Fp_{z=0} {1 \over \Gamma (z)}\int^1_0 x^{z-1} \sum^n_{q=0} 
(-1)^q Tr(B_q(u)e^{-x\Delta_q (u)}(\textrm{Id} - P_{q;u}))dx\\
&+ &{1 \over 2}\sum^n_{q=0}(-1)^q Tr(B_q(u)P_{q;u})\\
= &- &{1 \over 2} Fp_{z=0} {1 \over \Gamma (z)} \int^1_0 x^{z-1} \sum^n_{q=0} 
(-1)^q Tr(B_q(u)e^{-x\Delta_q (u)})dx
\end{eqnarray*}

\n where $P_{q;u} \equiv P_{q;u}(\{ 0\})$ is the orthogonal projector onto 
the null space of $\Delta_q(u)$.  

Introduce
\begin{eqnarray*}
{\bf (4.9)}~~~~~ a(x;u) &:= &\sum^n_{q=0} (-1)^q Tr (B_q(s;u)e^{-x\Delta_q 
(u)})\\
&= &\sum^n_{q=0} (-1)^q Tr ((\textrm{Id} +\Delta_q(u))^{-s}\\
& &~~~\qquad~~~A_q(u)^{-1}{d \over 
du} (A_q(u)(\textrm{Id}+\Delta_q(u))^s)e^{-x\Delta_q(u)})\,.
\end{eqnarray*}

\n To investigate the right hand side of (4.8) introduce

\spd

\n {\bf (4.10)}\hsps $b(x;u):=\sum\limits^n_{q=0} (-1)^q 
Tr \left( A_q(u)^{-1}\left({d 
\over du} A_q(u)\right) e^{-x\Delta_q(u)}\right).$

\spd

\n According to [BZ],
\[
Fp_{z=0}{1 \over \Gamma (z)}\int^1_0 x^{z-1} b(x;u)dx
\]
is, in the case where  $\mu(u)= \mu$ is constant, given  by

\spd

\n {\bf (4.11)}\hsps $-\int_M \theta (\rho, \mu)\wedge {d \over du} 
e_{CS}(M,g(0),
g(u))$

\spd

\n and,  when $g(u)=g$  is constant by 

\spd

\n {\bf (4.12)}\hsps $\int_M {d 
\over du} V (\rho, \mu(u), \mu(0)) e(M,g).$

\spd

\n Actually, in [BZ], these formulae are proven only for the case $\cal A=\bbc$,
but the same arguments work "word by word" for $\cal A$ arbitrary. In this paper $\theta(\rho,\mu)$ is defined to be $1/2$ of the corresponding quantity in [BZ].

\spd

\n We will show that 
\[
a(x;u)=b(x;u)\,.
\]

\n From the equality
\[
{d \over du}((\textrm{Id}+\Delta_q(u))^s (\textrm{Id} +\Delta_q(u))^{-s})=0
\]
and the commutativity of the trace, $Tr AB=TrBA$, we obtain
\begin{eqnarray*}
&~ &Tr((\textrm{Id}+\Delta_q(u))^{-s} A_q(u)^{-1}{d \over du}(A_q(u)
(\textrm{Id}+
\Delta_q(u))^s)e^{-x\Delta_q(u)})\\
&~ &= Tr(A_q(u)^{-1}({d \over du} A_q(u)) 
e^{-x\Delta_q(u)})\\
&~ &~~~~~~~~~~+ Tr((\textrm{Id}+\Delta_q(u))^{-s} {d \over ds}
((\textrm{Id}+\Delta_q(u))^{-s})
e^{-x\Delta_q(u)})\\
&~ &= Tr (A_q(u)^{-1}({d \over du} A_q(u)) e^{-x\Delta_q(u)}) +\\
&~ &~~~~~~~~~~+s\,Tr ((\textrm{Id}+\Delta_q(u))^{-1}{d \over du}
(\Delta_q(u)) e^{-x\Delta_q(u)})\\
&~ &= Tr (A_q(u)^{-1}{dA_q(u) \over du} e^{-x\Delta_q(u)})+
s {d \over du} Tr f(\Delta_q
(u))
\end{eqnarray*}
where $f(x,y):= -\int^\infty_y {e^{-xr} \over (1+r)} dr$.

\n Since $\Delta_q(u)$ is selfadjoint and nonnegative one obtains, for $x>0$, 
\[
\Vert f(x,\Delta_q(u))\Vert_{tr} := 
\int^\infty_0 \lambda\vert f(x,\lambda )\vert d
N_{\Delta_q(u)}\le \int^\infty_0\lambda e^{-\lambda x} dN_{\Delta_q(u)}(\lambda ) 
<\infty\,.
\]
Therefore, $f(x,\Delta_q(u))$ is of trace class for $x>0$ and 
\begin{eqnarray*}
\sum^n_{q=0}(-1)^q Tr f(x,\Delta_q(u)) &= &
\sum^n_{q=0} (-1)^q Tr f(x,\Delta^+_q(u)) +\\
&~ &~~\quad~~+\sum^n_{q=0} (-1)^q Tr f(x,\Delta^-_q (u))\\
&= &\sum^n_{q=1} (-1)^q (Tr f(x,\Delta_q^+(u))-Tr 
f(x,\Delta^-_{q-1}))\\
&= &0
\end{eqnarray*}
where we used that $\Delta^+_q(u)$ and $\Delta^-_q(u)$ are 
isospectral and, therefore, 
\[
Tr f(x,\Delta^+_q(u))=Tr f(x,\Delta^-_q(u))\,.
\]
Thus, for any $x\ge 0$,
\begin{eqnarray*}
a(x; u) &= &b(x; u) + s{d \over du}\sum^n_{q=0} (-1)^q Tr f(x,\Delta_q(u))\\
&= &b(x; u)\,.
\end{eqnarray*}

\n Note that by Proposition 1.12 
\[\log T({\cal C}(\iota (u))= \sum_{x\in Cr h} (-1)^{\rm{ind}\ x}V(\rho, \mu(u),\mu(0))(x).\]

Combining these with (4.8)-(4.11), statement (i) in Proposition 4.1
and  
Proposition 4.2
follow. \ \ \ \carre

\spd

\subsection{ Proof of Proposition 0.1}

We now prove Proposition 0.1 in full generality.
Let $g$ be an arbitrary Riemannian metric of $M$ and $g'$ the metric
given by the generalized triangulation $\tau.$
 Denote the corresponding relative torsions by
${\cal R} (g)$ and ${\cal R} (g').$  Proposition 4.1 and 
Proposition 0.1 in the special case considered in section 3
 imply Proposition 0.1: Indeed, denote by $\alpha (g')$ the density given
 in subsection 3.1 (D) and define 
 $ \alpha (g) := \alpha (g') +  \theta\wedge e_{CS} (g,g')$
 Clearly, $\alpha $ is a local quantity and  satisfies (cf subsection 3.1 (B) and (D))
$\log {\cal R}(g)  = \int_{M \setminus Cr(h)} \alpha (g).$
Further $\alpha(g)$ vanishes if $\mu$ is parallel 
(cf subsection 3.1 (D) and Proposition 4.1(ii)).
\ \ \ \carre

\subsection{ Anomaly with respect to triangulations}

\n In  this subsection we investigate how the relative torsion depends 
on the triangulation. 

\n Suppose $\tau_1 = (h_1,g_1)$ and $\tau_2 = (h_2,g_2)$ are 
two generalized 
triangulation which admit a common subdivision $\tau_0 = (h_0,g_0)$.  Recall 
from [BFKM, Subsection 6.3] that $\tau_0$ is called a subdivision of $\tau_1$ if :

\[ (i) \  Cr_q(h_1)\subseteq Cr_q (h_0), \  (0\le q\le n), 
\  W^{-}_x(\tau_0)\subseteq 
W^{-}_x (\tau_1), \  (x\in Cr(h_1) )\] 
and for any $y\in Cr (h_1) $ 
\[ (ii) \  W^{-}_y(\tau_1)=\cup_{x\in Cr(h_0), x\in  W^{-}_y(\tau_1)} W^-_x(\tau_0),\] 
where $Cr_q(h_{i})$ denotes the set of critical 
points of index $q$  of $h_{i} \ (i=0,1,2),$   $W^{-}_x(\tau_1)$ denotes the 
unstable  manifold of $x\in Cr(h_1)$ with 
respect to the gradient flow $-\textrm{grad}_{g_1} (h_1)$ and
$g_0=g_1,\  h_0=h_1$ in a neighborhood of the critical points of $h_1.$

\n Given $x\in Cr(h_0)$, there exist unique elements $x_1\in Cr(h_1)$, $x_2\in 
Cr(h_2)$ so that $x\in W^{-}_{x_j}(\tau_j)$ ($j=1,2$).

\spd

\n Introduce the function $w=
w_{\tau_0} = w(\tau_1,\tau_2 ;\tau_0)\colon Cr(h_0)\to \bbr$ by setting 

\spd

\n {\bf (4.13)}\hsps $w(x):=\log\textrm{vol}(T^{\tau_2}_{x,x_2}
\circ (T^{\tau_1}_{x,x_1})^{-1})$

\spd

\n where $T^{\tau_j}_{x,x_j}\colon {\cal E}_x \to {\cal E}_{x_j}$ ($j=1,2$) 
denotes the 
parallel transport along an arbitrary curve in $W^-_{x_j} (\tau_j)$ joining $x$ and 
$x_j$.  Define $\omega_{\tau_1,\tau_2} = \omega (\tau_1,\tau_2 ;\tau_0)$ by 

\spd

\n {\bf (4.14)}\hsps $\omega_{\tau_1,\tau_2} := \sum\limits_{x\in Cr(h_0)} 
(-1)^{\textrm{index} (x)} w(x)$.

\spd

\n Notice that $\omega (\tau_1,\tau_2;\tau_0)$ is independent of the choice of 
$\tau_0$, i.e., 

\spd

\n {\bf (4.15)}\hsps $\omega (\tau_1,\tau_2 ;\tilde{\tau}_0)= 
\omega (\tau_1,\tau_2;
\tau_0)$

\spd

\n for any generalized triangulation $\tilde{\tau}_0 = (\tilde{h}_0,\tilde{g}_0)$ 
which is a subdivision of $\tau_0$.  To see it, notice that for any 
$y\in Cr (\tilde{h}_0
)$, there exists a unique $x\in Cr(h_0)$ with  $y\in W^-_x(\tau_0)$.  Use that 
$T^{\tau_j}_{y,x_j} = T^{\tau_j}_{x,x_j} T^{\tau_0}_{y,x}$ ($j=1,2$) to 
conclude that 
$w_{\tilde{\tau}_0}(y)= w_{\tau_0}(x)$ and thus
\[
\sum_{y\in Cr (\tilde{h}_0)\cap W^-_x(\tau_0)} (-1)^{\textrm{index}(y)} 
w_{\tilde{\tau}_0} 
(y) = w_{\tau_0}(x)\sum_{y\in Cr(\tilde{h}_0)
\cap W^-_x(\tau_0)}(-1)^{\textrm{index}(y)}
\,.
\]

\n As $\tilde{\tau}_0$ is a subdivision of $\tau_0$
\[
\sum_{\begin{array}{c}
y\neq x\\
y\in Cr(\tilde{h}_0)\cap W^-_x(\tau_0)
\end{array}} (-1)^{\textrm{index}(y)} = 0
\]

\n we conclude that $\omega (\tau_1,\tau_2;\tilde{\tau}_0) = \omega (\tau_1,
\tau_2;\tau_0)$.  Moreover, if $\tau_1$, $\tau_2$ and $\tau_3$ are 
three generalized 
triangulation which admit a common subdivision $\tau_0$, then

\spd

\n{\bf (4.16)}\hsps $\omega_{\tau_1,\tau_2} +\omega_{\tau_2,\tau_3} = 
\omega_{\tau_1,
\tau_3}$.

\begin{prop} Assume that the generalized triangulations $\tau_1$ and $\tau_2$ 
admit a common subdivision $\tau_0$.  Then

\n (i)~~~$\log{\cal R}(\tau_1)-\log{\cal R}(\tau_2) =  \omega_{\tau_1,\tau_2}$

\n (ii)~~If $\mu$ is parallel with respect to the canonical flat connection of 
${\cal E}$, then ${\cal R}(\tau_1 )= {\cal R}(\tau_2 )$ provided $\tau_1$ 
and $\tau_2$
have a common subdivision.
\end{prop}

\n {\bf Proof:}  Statement (ii) follows from (i) by noticing that $w(x)=0$, 
$\forall~x
\in Cr(h_0)$ and thus $\omega_{\tau_1,\tau_2} = 0$.

\n To prove (i), notice that
\[
\log{\cal R}(\tau_2)-\log{\cal R}(\tau_1)= (\log{\cal R}(\tau_2)-
\log{\cal R}(\tau_0))+
(\log{\cal R}(\tau_0)-\log{\cal R}(\tau_1))\,.
\]
In view of (4.16), it suffices to consider the case where $\tau_2=\tau_0$.  
The subdivision 
$\tau_2$ can be obtained from $\tau_1$ by a sequence $\tau_1=\sigma_1,\dots ,
\sigma_N = \tau_2,$ 
\newline where $\sigma_{j+1} = (h_{j+1},g'_{j+1})$ is a subdivision of $\sigma_j =
(h_j,g'_j)$ 
with $Cr(h_{j+1})=Cr(h_j)\cup\{ x_{j+1},y_{j+1}\}$ so that there 
exists $z_j\in Cr(h_j)$ 
with $x_{j+1},y_{j+1}\in W^-_{z_j}$  and 
$\textrm {index}(x_{j+1})=\textrm{index}
(z_j)=\textrm{index}(y_{j+1})+1$.  Recall that $W^-_{x} \equiv  W^-_{x; h_j}$ denotes 
the unstable manifold, associated to $x \in Cr(h_j)$ and the 
gradient flow $-\textrm{grad}_{g_j} h_j.$  In view of (4.15), it suffices 
to consider the case 
where $\tau_2=\sigma_2$.  To ease notation, we write $\tau:=
\tau_1=(h_1, g'_1)$, $\sigma :=
\tau_2=(h_2, g'_2)$ $z:=z_1$, $x':=x_2$, $y':=y_2$.  Then, 
with $q_0=\textrm{index}(z)$
\[
Cr_q(h_2) = \left\{ \begin{array}{ll}
Cr_q(h_1) &q\not= q_0 ,q_0 -1\\
Cr_q(h_1)\cup\{ y'\} &q=q_0 -1\\
Cr_q(h_1)\cup\{ x'\} &q=q_0 \end{array}\right.
\]

\n Consider the following commutative diagram
\[
\begin{array}{ccc}
&\textrm{Int}_s(\tau)\\
H_s(\Lambda (M,{\cal E})) &\longrightarrow &{\cal C} (\tau )\\
&&\\
\textrm{Int}_s (\sigma )\searrow &&\nearrow A\\
&&\\
&{\cal C} (\sigma )\end{array}
\]
where $A_q$ ($0\le q\le n$) is defined as follows:  for a section $s\in 
C^{\infty} ({\cal E}
\mid_{Cr_qh_2} )$ and a critical point $x\in Cr_q(h_1)$, set
\[
A_q (s)(x) := \sum_{y\in Cr_q (h_2)\cap W^-_x} T_{yx} (s(y))
\]
where $T_{yx}\colon{\cal E}_y\to{\cal E}_x$ denotes the parallel 
transport from $y$ 
to $x$ along a curve in $W^-_x = W^-_x (\tau )$.  The map $A$ is a
morphism of cochain 
complexes of ${\cal A}$-Hilbert modules of finite type which induces 
an isomorphism in 
algebraic cohomology.

\n As $\textrm{Int}_s$ and $A$ are of trace class, we then conclude 
from Proposition 1.15 
that 
\[
\log {\cal R}(\tau ) = \log{\cal R}(\sigma )+ \log T({\cal C} (A))\,.
\]
Thus the claimed result follows once we show that

\spd

\n {\bf (4.17)}\hsps $\log T({\cal C} (A))= \omega_{\sigma\tau}$.

\spd

\n Formula (4.17) can be verified, using a localization argument:
 Consider first the case where $\mu$ is parallel. Clearly, $\omega_{\sigma\tau} = 
 0.$ On the other hand,
 $\log{\cal R}(\tau ) = \log{\cal R}(\sigma )= 0$  (Proposition 0.1 (ii)) and thus 
$\log T 
({\cal C} (A)) = \log{\cal R}(\tau ) - \log{\cal R}(\sigma )=0$ .
 Hence (4.17)  holds in this case.

\n Further, if ${\cal E}$ admits a Hermitian structure $\mu_0$ which is 
parallel, then 

\n (4.17) remains valid.  Indeed, consider the following commutative diagram
\[
\begin{array}{ccc}
&A\\
{\cal C} (\sigma ,\mu ) &\longrightarrow &{\cal C} (\tau ,\mu )\\
&&\\
\Big\downarrow \textrm{Id}_{\sigma ,\mu,\mu_0} &&\Big\downarrow 
\textrm{Id}_{\tau ,\mu,
\mu_0}\\
& A\\
{\cal C} (\sigma ,\mu_0) &\longrightarrow &{\cal C} (\tau ,\mu_0 )
\end{array}
\]

\n Then, according to Proposition 1.12,
\[
\log T({\cal C} (\textrm{Id}_{\sigma ,\mu,\mu_0}))+
\log T({\cal C} (A,\mu_0)) = \log T ({\cal C} 
(A,\mu ))+\log T({\cal C} (\textrm{Id}_{\tau ,\mu,\mu_0}))\,.
\]
By Proposition 1.14, ($j=1,2$)
\begin{eqnarray*}
\log T({\cal C} (\textrm{Id}_{\sigma ,\mu,\mu_0})) &- &
\log T({\cal C} (\textrm{Id}_{\tau ,
\mu,\mu_0}))\\
&= &\sum_q (-1)^q \log\textrm{vol}(\textrm{Id}_{q;\sigma ,
\mu,\mu_0})\\
&~ &-\sum_q (-1)^q 
\log\textrm{vol}(\textrm{Id}_{q,\tau ,\mu,\mu_0})\\
&= &\sum_{y\in Cr(h_2)}(-1)^{\textrm{index}(y)}\log
\textrm{vol}(\textrm{Id}_{y,\mu,
\mu_0})\\
&~ &-\sum_{y\in Cr(h_1)}(-1)^{\textrm{index}(y)} \log
\textrm{vol}(\textrm{Id}_{y,
\mu,\mu_0})\\
&= &(-1)^{\textrm{index}(x')}\log\textrm{vol}
(\textrm{Id}_{x',\mu,\mu_0})\\
&~ &+ (-1)^{\textrm{index}(y')}\log\textrm{vol}
(\textrm{Id}_{y',\mu,\mu_0})\,.
\end{eqnarray*}

\n Combining the above equalities with $\log T({\cal C} (A,\mu_0 )) =
0$ we obtain
\[\log T({\cal C} (A,\mu ))= (-1)^{\textrm{index}(x')}\log\textrm{vol}
(\textrm{Id}_{x',\mu,
\mu_0}) +\]
\[+(-1)^{\textrm{index}(y')}\log\textrm{vol}(\textrm{Id}_{y',\mu,\mu_0})\,.
\]
Observe that
\[
\log\textrm{vol}(\textrm{Id}_{x',\mu,\mu_0}) = 
\log\textrm{vol}_\mu (T_{x',z'}) +\log
\textrm{vol}(\textrm{Id}_{z',\mu,\mu_0})
\]
and, similarly,
\[
\log\textrm{vol}(\textrm{Id}_{y',\mu,\mu_0}) = 
\log\textrm{vol}_\mu (T_{y',x_0}) +
\log\textrm{vol}(\textrm{Id}_{z',\mu,\mu_0})\,.
\]

\n As $\textrm{index}(x')=\textrm{index}(y')+1$, this leads to
\begin{eqnarray*}
\log T({\cal C} (A,\mu )) &= &(-1)^{\textrm{index}(x')}
\log\textrm{vol}_\mu (T_{x'z'})\\
&+ &(-1)^{\textrm{index}(y')}\log\textrm{vol}_\mu (T_{y'z'})\\
&= &\omega_{z\tau}
\end{eqnarray*}
which establishes (4.17) in the case where ${\cal E}$ admits 
a Hermitian structure which 
is parallel.

\n Now use a standard localization to conclude that (4.17) is 
true in general. \ \ \ \carre

\section{Proof of Theorem 0.1}

\subsection{Proof of Theorem 0.1}

Consider the system 
 ${\cal{F}} := (M,\rho, \mu,  g, \tau ) 
 \equiv (M,{\cal{E}}, \nabla, \mu, g,  \tau ),$  with $\tau= (h,g')$ and 
denote by  $X$ the vector field $X:=-\textrm{grad}_{g'}h.$
$X$ defines a smooth map $X: M\setminus Cr(h) \to TM\setminus M.$ 
Denote by 
\[\beta_{\cal F} \equiv \beta (M,\rho, \mu, g, \tau) := 
(-1)^{n+1}\theta(\rho,\mu)\wedge 
X^{\ast}(\Psi (TM, g)) \in \Omega^n(M\setminus Cr(h); \cal O_M).\] 
Observe that if  $\mu$ is parallel in a neighborhood of $Cr(h)$ 
both $\alpha \equiv \alpha_{\cal F},$  by Proposition 3.2 (iii), and  $\beta_{\cal F}$ 
vanish in the neighborhood of $Cr(h).$

Proposition 3.2 (ii) and the equalities (4.1) and ($4.2^{v}$) imply that 
for
\newline $\gamma (M,\rho, \mu, g, \tau)= \alpha (M,\rho, 
\mu, g, \tau) \  \rm{or} \    
\beta (M,\rho, \mu, g, \tau)$ we have 
\begin{eqnarray*}
 {\bf (5.1)} \hsps \gamma (M_1\times M_2, \rho_1\otimes 
 \rho_2, g_1 \times g_2, 
\mu_1\otimes \mu_2, \tau_1\times \tau_2)=\\
\gamma (M_1, \rho_1, g_1, \mu_1, \tau_1) e(M_2,g_2) 
\rm{dim}({\cal {W}}_2) + \gamma (M_2, \rho_2, g_2, \mu_2, \tau_2) 
 e(M_1,g_1) \rm{dim}({\cal {W}}_1).
\end{eqnarray*}
\spd

Choose a Hermitian structure $\mu_0$ which is parallel in a neighborhood of $Cr(h).$
Introduce the quantity
\begin{eqnarray*}
{\cal S }(M, \rho, \mu, \mu_0, g, \tau) := \log {\cal R}(M, \rho, \mu, g, \tau)
- \int_{M\setminus Cr(h)} \beta (M, \rho, \mu_0, g, \tau)-\\
- \int_M V(\rho, \mu,\mu_0) e(M, g) +\sum_{x\in Cr(h)}(-1)^{\rm{ind} (x)}V(\rho, \mu, \mu_0)(x)
\end{eqnarray*}

The proof of Theorem 0.1 uses the following four statements:

\spd

\n {\em (A)  ${\cal S }(M, \rho, \mu, \mu_0, g, \tau)$ is independent of 
$ \mu,\mu_0, g, \tau.$}

\spd 

\n Therefore we write $ {\cal S }(M, \rho )$
instead of  ${\cal S }(M, \rho, \mu, \mu_0, g, \tau).$ 

\spd

\n {\em (B)  If $\rho$ is isomorphic to its dual, then ${\cal S}(M, \rho) +
(-1)^{n+1}{\cal S}(M,\rho)=0.$}

\spd

\n {\em (C) {${\cal S} (M, \rho) =(-1)^{n+1} {\cal S} (M,\rho^{\sharp})$  
where $\rho^{\sharp}$ is the dual  representation .}

\spd

\n {\em (D) {$ {\cal S} (M_1\times M_2,\rho_1\otimes \rho_2) = 
{\cal S}(M_1,\rho_1) 
\chi (M_2) \rm{dim} {\cal {W}}_1 +
\chi (M_1) {\cal S} (M_2,\rho_2)\rm{dim} {\cal {W}}_2. $} 
\spd

\spd
\n  (A) follows  from Propositions 4.1-4.3. and the formulas 4.1''-4.1''' and
4.2'-4.2'''. Precisely,
the independence of $g$ follows
from Proposition 4.1. and formula 4.2'', and the independence of $\mu$ from 
Proposition 4.2 and the formulas 4.1'' and 4.1'''.
To check the independence   on $\mu_0$ one uses $g =g'$ and proceeds 
as follows: 
Choose  coordinates in a  neighborhood of the critical points so that (T3) 
in the definition 
of a generalized triangulation is satisfied. Denote by 
$D_x(\epsilon)$ the disc of radius $\epsilon$ and centered at $x\in Cr (h)$ and
put $D(\epsilon):= \bigsqcup_{x\in Cr (h)}D_x(\epsilon).$
Let $\tilde \mu_0$ and $\mu_0$ be two Hermitian structures parallel in the neighborhood of 
the critical points.  A straightforward application of the definition
of ${\cal S}$ combined with $4.1'-4.2^v$ leads to 
\[{\cal S}(M,\rho, \mu, \tilde\mu_0, g', \tau) - {\cal S}(M,\rho, \mu, \mu_0, g', \tau)=\] 
\[=\lim_{\epsilon \to 0} ( (-1)^n \int_{M\setminus D(\epsilon)} d(V(\rho, 
\tilde\mu_0,\mu_0) X^\ast(\Psi(TM,g')) -\]
\[-(-1)^n  \int_{M\setminus D(\epsilon} 
V(\rho, \tilde\mu_0,\mu_0)e(M,g'))+\]
\[+\int_{M} 
V(\rho, \tilde\mu_0,\mu_0)e(M,g') 
-\sum_{x\in Cr (h)}
(-1)^{\rm{ind}\ x}(V(\rho, 
\tilde\mu_0,\mu_0)(x).\] 
Using Stokes theorem and, when $n$ is odd, the vanishing of 
$e(M,g')$  one obtains 
\[{\cal S}(M,\rho, \mu, \tilde\mu_0, g', \tau) - {\cal S}(M,\rho, \mu, \mu_0, g', \tau)=\] 
\[=\lim_{\epsilon \to 0} ((-1)^n \int_{\partial (M\setminus D(\epsilon))} 
V(\rho, \tilde\mu_0,\mu_0) X^\ast(\Psi(TM,g')) -\]
\[-\sum_{x\in Cr (h)}(-1)^{\rm{ind}\ x}(V(\rho, 
\tilde\mu_0,\mu_0)(x)=\]
\[=\lim_{\epsilon \to 0} ( (-1)^{n+1} \int_{\partial (D(\epsilon))} 
V(\rho, \tilde\mu_0,\mu_0) X^\ast(\Psi(TM,g')) -\]
\[-\sum_{x\in Cr (h)}(-1)^{\rm{ind}\ x}V(\rho, 
\tilde\mu_0,\mu_0)(x).\] 
Since $V(\rho, 
\tilde\mu_0,\mu_0)$ is constant on each $D_x(\epsilon)$ one obtains 
\[{\cal S}(M,\rho, \mu, \tilde\mu_0, g', \tau) - {\cal S}(M,\rho, \mu, \mu_0, g', \tau)=\] 
\[ =\sum_{x\in Cr (h)} V(\rho, 
\tilde\mu_0,\mu_0)(x)( (-1)^{n+1}\int_{\partial D_x(\epsilon)} 
X^\ast(\Psi(TM,g')) - (-1)^{\rm{ind}\ x}).\]

\n Since in the coordinates we have chosen 
\[X= \sum_{k=1}^{\rm{ind}\ x}
x_k \frac{\partial}{\partial x_k} -  \sum_{k=\rm{ind}\ x+1}^n
x_k \frac{\partial}{\partial x_k}\] $4.2^{iv}$
implies that 
\[  (-1)^{n+1}\int_{\partial D_x(\epsilon)} 
X^\ast(\Psi(TM,g')) = (-1)^{\rm{ind}\ x}\]
and therefore 
\[{\cal S}(M,\rho, \mu, \tilde\mu_0, g', \tau) = {\cal S}(M,\rho, \mu, \mu_0, g', \tau).\]

\n The independence of $\tau$ can be verified in the following way:
One considers a subdivision $\tau'$ of $\tau$ so that only  those cells  
of $\tau $, lying in a 
contractible open set $U,$
are subdivided. We choose $ \mu_0$ to be parallel on $U$
in addition of being  parallel in a neighborhood of the critical points of $h$
and let $\mu := \mu{_0}.$
Since 
$\omega_{\tau, \tau'}$ as defined in subsection 4.3 is zero, it follows from  
Proposition 4.3 that 
$ {\cal S} (M,..,\tau)= {\cal S} (M, ..,\tau').$

\n If $\tau_1$ and $\tau_2$ are two generalized triangulations, 
one can find a third
triangulation, $\tau_0,$ which is simultaneously a  
subdivision of $\tau_1$ and , up to a composition with a diffeomorphism of $M$ isotopic to
the identity, a subdivision of $\tau_2.$  Since  one can pass from 
$\tau_i, \ i=1,2$ 
to $\tau_0$  by a finite sequence of subdivisions each of them  
involving subdivisions 
of cells lying in a contractible open set,  the result follows. 

\spd

\n To check  (B) we set $\mu: =\mu_0$ and choose $\mu_0$ 
with $ \theta (\rho, \mu_0)= 0.$ In view of the assumption that $\rho$ is 
isomorphic to $\rho^{\sharp}$ 
this is possible (cf subsection 4.1). 
The result follows then from Proposition 3.2 (i). 

\spd

\n To check (C )  we choose again $\mu= \mu_0.$ 
As $\theta(\rho^{\sharp},\mu^\sharp_0)= - \theta(\rho,\mu_0)$ and 
 $X_{\tau_{\cal D}} = -grad_{g'}(n-h)= grad_{g'} h= -X_{\tau}$
 we conclude in view of 4.2''' that 
\[\int_M\theta(\rho,\mu_0)\wedge X^{\ast}_{\tau}\Psi + 
(-1)^n \int _M \theta(\rho^{\sharp},\mu_0^{\sharp})\wedge 
X^{\ast}_{{\cal D}_\tau}\Psi
=0.\] 

\n Statement (C ) follows once we verify that 
\[\int_M\alpha (M,\rho, g, \mu, \tau,\epsilon)+  
(-1)^n \int _M \alpha (M,\rho^{\sharp}, g, \mu^{\sharp}, 
\tau_{\cal D}, \epsilon)=0.\] 

\n With the notation from subsection 3.2 and  in view of formula 3.13 
and the fact that 
the $q$-Laplacian corresponding
to $(\rho, \mu, g, h)$ is conjugated by the Hodge star operator to the 
$(n-q)$-Laplacian corresponding to $(\rho^{\sharp}, \mu^{\sharp}, g, n-h),$
the left hand side of the equality above is $\int _M \delta_{\cal F}
(\epsilon)$ with

\spd

\n {\bf (5.2)} \hsps $\delta_{\cal F}(\epsilon)= \delta_{\cal F}(h,\epsilon):=
n/2  \sum_q (-1)^{q+1} \alpha_{\cal F}^q(h,\epsilon).$

\spd

\n We want to check that $\int _M \delta_{\cal F}
(0)=0 $ and proceed  as in the proof of Proposition 3.2 (iv). Given the system 
${\cal F}$ one chooses 
a system $\tilde{\cal F}$  with the same underlying 
Riemannian manifold $(M,g),$ the same triangulation $\tau,$ but with 
$\tilde{\rho}$ the 
trivial representation 
over the same  Hilbert module and  with $\tilde{\mu}$  a parallel 
Hermitian structure. 

\n In view of 5.2, introduce  
${W}_{\cal F}(t,\epsilon ):= 
\frac{n}{2} \sum_q (-1)^{q}  \log \det ({\Delta}_q (t) + \epsilon).$ 
By the same arguments as in the proof of Proposition 3.2 (iv) (Mayer
Vietoris arguments), we conclude that

$$FT_{t=\infty} (W_{\cal F}(t,\epsilon) - W_{\tilde{\cal F}}(t,\epsilon))=
\int_{M \setminus Cr(h)} \delta_{\cal F}(\epsilon) - \int_{\tilde{M}\setminus 
Cr(\tilde{h})} \delta_{\tilde{\cal F}}(\epsilon).$$

\spd

The left side of the above identity is zero because by a straightforward computation 
\[W_{\cal F}(t,\epsilon)= n\chi (M) \rm{dim}{\cal W} \log \epsilon = W_{\tilde {\cal F}}(t,\epsilon).\]
On the right  side of the identity, the term 

\n $\int_{{M}\setminus Cr(h)}
\delta_{\tilde{\cal F}}(\epsilon) = \int_{M\setminus Cr(h)}(\alpha _{\tilde{\cal F}}+ (-1)^n 
\alpha _{\tilde{\cal F}^{\sharp} })$ 
is zero  because in view of Proposition 0.1 both 
$\int_{M \setminus Cr(h)}\alpha_{\tilde{\cal F}}$ and $\int_{M \setminus Cr(h)} 
\alpha_{\tilde{\cal F}^{\sharp}}
(\epsilon)$ are zero for the trivial representation and a parallel  Hermitian 
structure. 

\n Consequently $\int_{M \setminus Cr(h)} \delta_{\cal F}(0)=0 .$

\spd

\n ( D) Choose  $\mu:= \mu_0.$ Then the statement follows from formula 5.1 
Proposition 0.1 and Proposition 3.2 (ii).  
\spd

\n {\bf Proof of Theorem 0.1} {\it Step 1:} We first show that ${\cal S}(M,\rho)$ is of the form 
${\cal S}(M,\rho ) = F(\rho ; \Gamma) \cdot \chi (M).$
Note that when $\rho$ is  isomorphic to $\rho^{\sharp },$  
it follows, by using (B) for $n$ odd and (C) for $n$ even,  that ${\cal S}(M, \rho) = 0.$
\n If  $\chi (M)= 0,$
choose  an even dimensional manifold $M'$ with nonzero Euler Poincar\'e 
characteristic
and with the same fundamental group as $M$ and choose the representation 
$\rho' :=\rho^{\sharp}.$

\n Since $\rho\otimes\rho^{\sharp}$ is isomorphic to its dual, and in view 
of (D)  
\[0= {\cal S}(M\times M',\rho\otimes \rho^{\sharp})= {\cal S}(M,\rho)
\cdot \chi (M')\dim ({\cal W}).\]
Hence 
the result is true if
$\chi(M)=0.$

\vspace{0.3cm}

Suppose  $\chi(M)\ne  0.$ For any  $M'$ with the same fundamental group as $M$ 
and $\rho':= \rho^{\sharp},$ by the same argument as above for the case $\chi (M) =0,$
one concludes from (D) that 
\[{\cal S}(M',\rho^{\sharp})=- {\cal S}(M,\rho)\cdot 
\frac{\chi (M')}{\chi(M)}.\]
This implies that for any closed manifold with fundamental group $\Gamma$,
${\cal S}(M,\rho)= \chi (M)\cdot F(\rho; \Gamma)$, where  
$F(\rho; \Gamma)$ depends only on the representation $\rho$ up to an isomorphism. 
\n Notice that if $\rho$ is isomorphic to its dual $\rho^{\sharp}$, then by (B) and (C),
$F(\rho; \Gamma) \cdot \chi (M) = 0 .$

\spd

\n {\it Step 2:} 
It remains to show that $F(\rho;\Gamma) = 0 .$ 
Given $\Gamma$, choose a closed manifold $(M, x_0)$ with $\dim M = 2n$ and 
$\Gamma =\pi_1(M, x_0).$
Let $M_1 := M$, $M_2 := M \# \bbc P(n)$, $M_3 := \bbc P(n)$, $M_4 := S^{2n}.$
Observe that $  \pi_1(M_2, x_0) = \pi_1 (M_1, x_0) = \Gamma$ and $\chi (M_2) = \chi (M_1) + (n-1).$
Let $({\cal E}_i \to M_i , \nabla _i)$  be the following flat bundles: for $i=1,2$ the ones induced
by the representation $\rho$ and 
for $i=3,4$ the flat bundles induced  by $\varepsilon,$ the unique representation on ${\cal W}$
of the group with one element.
We show below that one can choose Hermitian structures, Riemannian metrics and generalized 
triangulations so that,
by the locality  of the density given in Proposition 0.1 (cf subsection 0.2, (P1),) 
$$
\n {\bf (5.3)} \hsps
{\cal S}(M \# \bbc  P (n) , \rho )  - {\cal S}(M, \rho)
= {\cal S}(\bbc P (n)  , \varepsilon) - {\cal S}(S^{2n} , \varepsilon) 
$$
On the other hand,  by Step 1,
$$
\n {\bf (5.4)}
{\cal S}(M \# \bbc P (n), \rho) -  {\cal S}(M,\rho)
=
F(\rho) \cdot (n-1) ;\,\,\, 
{\cal S}( \bbc P (n) , \varepsilon) = 0, \,\,\, {\cal S}(S^{2n} , \varepsilon) = 0
$$
and the claimed result follows.

\n To prove (5.3),
 decompose $M_i$ as a union of compact manifolds with boundary, $M_i = M_i(-) \cup M_i(+)$
so that 
\begin{itemize}
\item[(a)] $ M_i(-)  \cup M_i(+) = M_i; \,\,\,\,   N_i := M_i(-)  \cap  M_i(+)  \sim S^{2n-1};$ 
\item[(b1)] $ M_1(-) \sim M_2(-) \sim M \setminus Int(D^{2n});$
\item[(b2)] $ M_1(+) \sim M_4(+) \sim D^{2n};$
\item[(b3)] $M_2(+) \sim M_3(+) \sim \bbc P(n) \setminus Int(D^{2n});$
\item[(b4)] $M_3(-) \sim M_4(-) \sim D^{2n}.$
\end{itemize}

\n Let  $W_i(-) $ and $W_i(+)$ be the bordisms
(cf [BFK2]), 
$$W_i(-) :=  ( M_i(-), \emptyset , \partial(M_i(-)) ); \,\,\,\
W_i(+) := ( M_i(+),  \partial (M_i(+)), \emptyset ).
$$ 
The isomorphisms in (b1) - (b4) then become isomorphisms of bordisms.

\n Denote by $({\cal E}_i(+), \nabla _i)$ respectively $({\cal E}_i(-), \nabla _i)$ 
the restriction  of the flat bundles  ${\cal E}_i$ to $M_i(+),$  
respectively to $M_i(-).$ 
Further we have the isomorphisms of flat bundles
\begin{itemize}
\item[(c1)] $({\cal E}_1(-), \nabla _1) \sim ({\cal E}_2(-), \nabla _2);$
\item[(c2)] $({\cal E}_1(+), \nabla _1) \sim ({\cal E}_4(+), \nabla _4);$
\item[(c3)] $({\cal E}_2(+), \nabla _2) \sim ({\cal E}_3(+), \nabla _3);$
\item[(c4)] $({\cal E}_3(-), \nabla _3) \sim ({\cal E}_4(-), \nabla _4).$
\end{itemize}

Choose Riemannian metrics $g_i$ on $M_i$ which are product like near $ N_i = M_i(-)  \cap  M_i(+) $
and generalized triangulations $\tau _i = (h_i, g_{i}')$ with $g_i'$ product like near $N_i$,
$0$ a regular value for $h_i,$ and $h_i^{-1}(0) = N_i.$ They induce generalized triangulations on $W_i(-)$ and $W_i(+).$
Choose Hermitian structures on ${\cal E}_i,$ parallel near the critical points  
so that (c1)-(c4) are isomorphisms of
Hermitian flat bundles. (This is always possible since the space of Hermitian structures which are parallel on a subset 
is contractible.) Using the locality  of the density given in Proposition 0.1 (cf subsection 0.2, (P1)),
(5.3) follows.
\spd

\n {\bf Remark to Theorem 0.1} Theorem 0.1 answers the question raised by several experts in the field, if the formula of Bismut-Zhang can be extended to $L_2-$torsion. 
Consider the system  $(M, {\cal E}, \nabla, \mu, g, \tau )$ where
$(M,g)$ is a closed Riemannian manifold, $(\cal E,\nabla)$  
a finite dimensional flat bundle
equipped with a Hermitian structure $\mu$ and $\tau$ a 
generalized triangulation of $M$. Suppose
that $M$ is connected and let $x_0\in M.$
Let $\Gamma$ be the fundamental group of $M,$ $\Gamma = \pi_1(M,x_0),$
$\cal N$ the  finite von Neumann algebra associated to $\Gamma,$  
 $\rho$ the holonomy 
representation of $\Gamma$ on ${\cal E}_{x_0}$
and  $\rho_{reg}$ the regular representation 
of $\Gamma$ on the $\cal N-$Hilbert module $l_2(\Gamma).$ 
Denote by $\rho'$ the representation of 
$\Gamma$ 
on the $\cal N-$Hilbert module $l_2(\Gamma)\otimes_{\bbc}{\cal E}_{x_0}$ 
given by 
$\rho':= \rho_{reg}\otimes_{\bbc}\rho,$ 
by $\cal M$ the flat bundle associated with 
the regular representation $\rho_{reg}$ and by $\mu_0$ the canonical parallel 
Hermitian structure on ${\cal M}$ induced from the  scalar product 
on $l_2(\Gamma).$ 
The flat bundle induced by the representation $\rho'$ is 
$\cal M\otimes_{\bbc}\cal E$  and is equipped with the Hermitian structure 
 $\mu'= \mu_0\otimes \mu.$
 Under the hypothesis of 
determinant class  for the pair $(M,\rho')$ (conjecturally
always satisfied,) one can use Lott's approach [L] to define 
the $L_2-$ 
analytic
torsion of  $(M, {\cal E}, \nabla, \mu, g )$ resp. the $L_2-$Reidemeister 
torsion of 
$(M, {\cal E}, \nabla, \mu, g, \tau).$ 
It is straightforward to check that these torsions are the analytic torsion 
of $(M, g, \rho', \mu')$ resp.  
the Reidemeister torsion of 
$(M, g, \rho', \mu', \tau)$  and that the local quantities expressing 
$\log {\cal R},$ cf  Theorem 0.1, are the same for 
$(M, \rho',\mu', g, \tau)$ and 
$(M, \rho,\mu, g, \tau).$ 
Therefore, by (2.15), the quotient of the $L_2-$analytic and 
$L_2-$Reidemeister 
torsions is the same as the quotient of  
analytic and Reidemeister torsions .
 \ \ \ \carre

\spd

 \subsection{An invariant for odd dimensional manifolds.}

For $\Gamma$ a finitely presented group, and denote by 
$Rep( \Gamma , {\cal W})$ the space
of representations of $\Gamma$ on the ${\cal A}$-Hilbert module of finite type
$\cal{W}$.
This space has a structure of a complex analytic space of infinite dimension when $\cal{A}$ is of infinite dimension. 

Let $(M,x_0)$ be a base pointed, closed, smooth odd dimensional 
manifold, with $\pi_1 (M, x_0)= \Gamma,$ and let $E$ be 
an Euler structure in the sense of Turaev cf [Tu]. Recall cf [B2],[Tu] 
that the 
set of Euler structures on a base pointed connected manifold 
$(M,x_0)$ can be defined as the set of connected components of the space 
of continuous vector fields with the only zero at 
$x_0$ equipped with the compact open topology.
Recall also that $H_1(M;\bbz)$ acts on this set freely and transitively 
and because $\chi(M)=0$ there is a canonical identification between 
sets of Euler structures obtained by using different base points.
By using the relative torsion we will 
associate a (real analytic) function 
$F(M,E):Rep( \Gamma , {\cal W})\to \bbr$ which is a 
smooth invariant for the triple $(M,x_0, E).$ We 
will calculate $F(M,E)$ 
for $M= S^1\times N,$ 
$N$ is a closed simply connected even dimensional manifold 
and $E$ the canonical Euler 
structure, defined below.
In this case $\pi_1 (M, x_0)= \bbz$ and  
the space $Rep ( \bbz , {\cal W})$ can be identified with 
$Gl_{\cal A} ({\cal W}).$ Therefore we will denote a 
representation $\rho$  by $A\in Gl_{\cal A} ({\cal W}), \ A=\rho(1).$

\spd
{\bf Definition of F(M,E):}

First note that  an Euler class $E$  can be represented  by the following data: 

\n (1)  a generalized triangulation $\tau= (h,g'),$ 

\n (2)  a spray, i.e a collection of disjoint smooth embeddings 
$\alpha_x:[0,1]\to M$ with $\alpha_x (0)= x_0, 
\ \alpha_x(1)= x, \ x\in Cr (h).$ Let $K:=\cup_{x\in Cr (h)}\alpha_x([0,1])\subset M$ 
and refer to $K$ as the spray defined by the curves $\alpha_x.$ 

The data (1) and (2) provide a  vector field $X$ 
by taking
$X=-grad_{g'} h$ outside of  
a contractible  smooth regular neighborhood $U$ of $K$ and extending $X$ 
inside $U$ with the only zero at $x_0.$ This vector field determines an Euler structure and 
any Euler structure can be obtained in this way (cf [B2]).
    
\n For a representation $\rho\in Rep( \Gamma , {\cal W})$ denote by 
$\cal{E}_{\rho}$ the flat bundle 
induced from $\rho$ and choose a Hermitian structure $\mu$ in 
$\cal{E}_{\rho},$ so that $\mu$  
is parallel in the neighborhood of $x_0$ and $\mu(x_0)$ is exactly the scalar product of 
$\cal{W}$ ( ${\cal{E}}_{\rho} (x_0)= {\cal{W}}).$ 

Define
 $$F(M,E)(\rho):=\log {\cal R}(M,\rho,\mu,\tau, g).$$
Using Propositions 4.1-4.3 one verifies that $F(M,E)(\rho)$ is well defined, 
i.e. is independent of the choice of triangulation, metric and  
Hermitian structure cf [B2] Main Theorem (1).
\spd

\n Let $M= N \times S^1$ with $N$ simply connected
and of even dimension. 

\n Choose the generalized triangulation $\tau_{S^1}$
of $S^1= \bbr /\bbz$
given by the canonical metric $g_0$ and a Morse function with 
two critical points,
a minimum at $t_1$ and a maximum at $t_2,$ and choose as a base point $t_0=1/2.$
Suppose that 
$0< t_1< 1/2<  t_2<  1.$ 
Take $K_{S^1}= [t_1,1/2]\cup [1/2,t_2].$ 

\n Choose a generalized triangulation $\tau_N$ of $N,$ a base point $n_0$  
and a spray $K_{N}.$

\n Consider the Euler structure on $S^1\times N$ represented by 
$\tau_{S^1}\times 
\tau_N$ and spray $K= K_{S^1}\times K_N$ with base point $x_0= t_0\times n_0.$ 
Denote  by $E_0$ the Euler structure defined by these data. Since $N$ is 
simply connected this Euler structure is independent of the spray $K_N$
and the triangulation $\tau_N$ and  will be refered 
as the canonical Euler structure.

\n The flat bundle induced by $A$ is actually the tensor product of a 
flat bundle induced by $A$ on $S^1$ for which we choose a Hermitian 
structure parallel above $[t_1-\epsilon,
t_2+\epsilon],$ with the trivial line bundle on $N,$ for which we choose 
a parallel Hermitian structure.


\begin{prop} For the canonical Euler structure $E_N$ described above
 
$F(N \times S^1,E_0)(A) = -\chi (N)/2 \log det (A^{\ast} A)^{1/2}.$  

\end{prop}

\n {\bf Proof} By the product formula for the relative torsion 
(use $\chi(S^1) = 0$)
\[\log {\cal R} (N \times S^1 , A, \mu, g, \tau ) = \chi(N) \log {\cal R} 
(S^1 , A, \mu_{S^1}, g_0, \tau_{S^1} ).\]

\n We compute $\log {\cal R} (S^1, A, \cdots)$ using Theorem 0.1. 

Note that $\mu_{S^1}$ can be viewed as a smooth family $\mu(t),$ $t\in \bbr$
of scalar products in $\cal{W}$ with $\mu(t+1)={A^{-1}}^\ast\mu(t)A^{-1}$ 
and $\mu(t)= \mu_0$ for $t\in [t_1+k-\epsilon, t_2+k+\epsilon],\ k\in \bbz$ 
and $\mu_0$
the scalar product of $\cal{W}.$
Recall (subsection 4.1) that 
$\theta (\mu ) (t) = d \alpha (t) $  where $ \alpha (t) = \log vol
(id:({\cal{W}},\mu(t))\to ({\cal{W}},\mu_0))$ which in particular says that 
$\alpha(t)=0$ for $t\in [t_1-\epsilon,
t_2+\epsilon].$
Clearly
$
\alpha (1)-\alpha (0) =  \log vol A^{-1}.
$

Further, $\Psi \in \Omega ^0 ( TS^1 \setminus S^1 ) $ (cf property $4.2^{iv}$)  is the function given by
$ \Psi (t, \xi) = 1/2 \,\,(\forall \xi> 0)$ and $ \Psi (x, \xi) = - 1/2 \,\,(\forall \xi < 0).$

The vectorfield
$X : S^1 \setminus Cr(h) \to TS^1 \setminus  S^1$ 
is given by $-grad_{g'} h$ and thus the pullback of $\Psi$ by
$X$ is 
\[
X^{\ast} (\Psi) (t) = 1/2\ \  (\forall t\in [0,1]\setminus[t_1,t_2])
\] 
and 
\[
X^{\ast} (\Psi) (t) = -1/2\ \  (\forall t\in (t_1,t_2))
\]
\n Note that $X^{\ast} (\Psi) (t) =-1/2$ on the points where 
$\theta (\mu)$ is zero.
Combining the above computations, one obtains
$$
\int_{ S^1 \setminus Cr(h)}  \theta (\mu)   X^{\ast}(\Psi) = 
1/2 \int_{ S^1 \setminus Cr(h)}  \theta (\mu) =
\frac{1}{2}  \log vol (A^{-1}).
$$ 
By the formula given in Theorem 0.1, the claimed statement follows.\\ \ \carre

\spd

\n {\large{\bf Appendix A \ Lemma of Carey-Mathai-Mishchenko}}

\vspace{1cm}

\n In this Appendix, we prove Lemma 1.14, using a deformation argument. 

\spd

\n In\-tro\-du\-ce
$f_q (t) := t f_q$
and $d_q (t) := \bl {cc} d_{1,q} & tf_q \\ 0 & d_{2,q} \er$
and obtain in this way a complex $({\cal C}^{\ast},
d_{\ast} (t)).$

\spd

\n By the same arguments as in the proof of Lemma 1.10 one concludes that
$({\cal C}^{\ast}, d_{\ast} (t))$ is an algebraically acyclic, $\zeta$-regular
complex of sF type.

\spd

\n Therefore, $({\cal C}^{\ast}, d_x (t))$ has a well defined torsion 
$T(t) := T({\cal C}^{\ast}, d_{\ast} (t)),$ given by (cf 
\linebreak
(1.21))

\[\log T (t) = \frac{1}{2} \sum (-1)^q \log \det (\underline{d}^{\ast}_q (t)
\underline{d}_q (t)).\] 

\n If $\frac{d}{dt} \log T (t) = 0,$ then

\[\log T ({\cal C}) = \log T(1) = \log T (0) = \log T ({\cal C}^1) + \log T 
({\cal C}^2)\]

\n where for the last equality we have used that $d_q (0) = \bl {cc}
d_{1,q} & 0 \\
0 & d_{2, q} \er$ and therefore $\Delta_q (0) = \bl {cc} \Delta_{1,q} & 0 \\
0 & \Delta_{2,q} \er.$

\spd

\n The remaining part of this Appendix is devoted to the proof of the 
statement 

\[\frac{d}{dt} \log T (t) = 0.\]

\n Introduce the Hodge decomposition of ${\cal C}^j, {\cal C}^{j}_q
= {\cal C}^{j+}_q
\oplus
{\cal C}^{j-}_q (j = 1,2).$ The differential $d_{j,q}$ then have the form 

\[d_{j,q} = \bl {cc}
0 & \underline{d}_{j,q} \\
0 & 0 \er .\]

\n Further, decompose ${\cal C}_q,$

\spf

\n {\bf (A.1)} \hsps ${\cal C}_q = {\cal C}^{1+}_q \oplus {\cal C}^{1-}_q
\oplus {\cal C}^{2+}_q \oplus {\cal C}^{2-}_q.$

\spd

\n Then, $d_q (t)$ can be written as

\[ d_q (t) = \bl {cc|cc}
0 & \underline{d}_{1,q} & \alpha_q & \beta_q \\
0 & 0  & \varepsilon_q & \gamma_q \\
\cline{1-4}
0 & 0 & 0 & \underline{d}_{2,q} \\
0 & 0 & 0 & 0 \er \]

\n where $\alpha_q \equiv \alpha_q (t) = t \dot{\alpha}_q, 
\beta_q = t \dot{\beta}_q, 
\varepsilon_q = t \dot{\varepsilon}_q$ and $\gamma_q = t 
\dot{\gamma}_q$.
From 
$d_{q+1} (t) d_q (t) = 0$ we deduce

\spf

\n {\bf (A.2)} \hsps $\underline{d}_{1,q+1} \varepsilon_q = 0; \
\varepsilon_{q+1} \underline{d}_{2,q} = 0; $

\spf

\n {\bf (A.3)} \hsps $\underline{d}_{1,q+1} \gamma_q + \alpha_{q+1} 
\underline{d}_{2,q} = 0.$

\spd

\n As ${\cal C}^1$ and ${\cal C}^2$ are acyclic, $\underline{d}_{1,q}$
and $\underline{d}_{2,q}$ are isomorphisms, and (A.2) and (A.3) imply

\spf

\n {\bf (A.4)} \hsps $\varepsilon_q = 0$

\spf

\n {\bf (A.5)} \hsps $\gamma_q = - \underline{d}^{-1}_{1,q+1} \alpha_{q+1}
\underline{d}_{2,q}.$

\spd

\n Next, let us describe the Hodge decomposition of ${\cal C}_{\ast} 
= : {\cal C}^+_{\ast} (t) \oplus {\cal C}^-_{\ast} (t)$ 
\newline of $({\cal C}_{\ast}, 
d_{\ast} (t)).$

\spf

\n {\bf Lemma A.1}

\spd

\n (i) ${\cal C}^+_q = Ker d_q = \Big \{ (x_+ - \underline{d}^{-1}_{1,q} 
\alpha_q y_+, y_+, 0) \Big| x_+ \in {\cal C}^{1+}_q; \ 
y_+ \in {\cal C}^{2+}_q \Big \}$

\spd

\n (ii) ${\cal C}^-_q = Ker (d^{\ast}_{q-1} ) = \Big\{ (0, x_{-}, (
\underline{d}_{1,q}^{-1} \alpha_q)^{\ast} x_{-}, y_{-}) \Big| 
x_{-} \in {\cal C}^{1-}_q; \ y_{-} \in {\cal C}^{2-}_q \Big\}.$

\spf

\n {\bf Proof} The statements follow from a straightforward verification.
\ \ \ \carre
\spd

\n We want to compute the $t$-derivative of $\log T (t) = \frac{1}{2} 
\sum (-1)^q \log \det (\underline{d}^{\ast}_q (t) \underline{d}_q (t)).$ Notice
that

{\arraycolsep0.5mm
\[ \begin{array}{lll}
\frac{d}{dt} \Big( d_q^{\ast}(t) d_q(t)\Big) & = &
\bl {cc} 0 & 0 \\ f_q^{\ast} & 0 \er 
\bl{cc} d_{1,q} & tf_q \\ 0 & d_{2,q} \er
+ \bl {cc} d_{1,q}^{\ast} & 0 \\ tf_q^{\ast} & d^{\ast}_{2,q} 
\er \bl {cc} 0 & f_q \\
0 & 0 \er \\
&& \\
& = & \bl {cc} 0 & 0 \\ f^{\ast}_{q} d_{1,q} & tf^{\ast}_q f_q \er + \bl {cc}
0 & d^{\ast}_{1,q} f_q \\
0 & tf^{\ast}_q f_q \er  \\
&& \\
& = &  \bl {cc} 0 & d^{\ast}_{1,q} f_q \\ f^{\ast} d_{1,q} & 2t f^{\ast}_q 
f_q \er. \end{array} \]}

\n As $f$ is a morphism of trace class (cf Definition 1.10) we conclude that
$\frac{d}{dt} d^{\ast}_q d_q (t)$ and therefore $\frac{d}{dt} 
\underline{d}_q^{\ast}
(t) \underline{d}_q (t)$ are of trace class. In view of Proposition 1.3 and the
fact that $(\underline{d}^{\ast}_q (t) \underline{d}_q (t))^{-1}$ is bounded
we deduce that $\frac{d}{dt} (\underline{d}^{\ast}_q (t) \underline{d}_q (t))
(\underline{d}^{\ast}_{q} (t) \underline{d}_q (t))^{-1}$ is of trace class. 
Therefore

\[\frac{d}{dt} \log \det (\underline{d}^{\ast}_q (t) \underline{d}_q (t)) =
Tr (\frac{d}{dt} (\underline{d}^{\ast}_q (t) \underline{d}_q (t)) 
(\underline{d}_q^{\ast} (t) \underline{d}_q (t))^{-1}).\]

\spf

\n {\bf (A.6)} \hsps $
\begin{array}{l}
Tr (\frac{d}{dt} (\underline{d}^{\ast}_q (t) 
\underline{d}_q (t)) (\underline{d}^{\ast}_q (t) \underline{d}_q 
(t))^{-1} ) = \\
= Tr (( \frac{d}{dt} \underline{d}^{\ast}_q (t)) \underline{d}_q^{\ast} (t)^{-1})
+ Tr (\underline{d}_q (t)^{-1} \frac{d}{dt} \underline{d}_q (t)) \\
= 2 Re Tr ((\frac{d}{dt} \underline{d}_q (t)) \underline{d}_q (t)^{-1}).
\end{array}$

\spd

\n It is convenient to introduce $D_q := d_q (t)$ and denote by $D^{-1}_q$ the
operator given by $D^{-1}_q |_{{\cal C}_{q+1}^{\cdot}} = 0$ and $D^{-1}_q 
|_{{\cal C}^+_{q+1}} = \underline{d}_q (t)^{-1}.$

\spf

\n {\bf (A.7)} \hsps $Tr (D^{-1}_q \frac{d}{dt} D_q) = Tr (\underline{d}_q (t)^{-1}
\frac{d}{dt} \underline{d}_q (t)).$

\spd

\n With respect to the decomposition (A.1), $D_q^{-1}$ and 
$\dot{D}_q \equiv
\frac{d}{dt} D_q$  ($\dot{D}_q$ is in\-de\-pen\-dent of $t$) take the form

\[ D_q^{-1} = \bl{llll}
A_{11} & A_{12} & A_{13} & A_{14} \\
A_{21} & A_{22} & A_{23} & A_{24} \\
A_{31} & A_{32} & A_{33} & A_{34} \\
A_{41} & A_{42} & A_{43} & A_{44} \er ; \dot{D}_q = \bl {llll}
0 & 0 & \dot{\alpha}_q & \dot{\beta}_q \\
0 & 0 & 0 & \dot{\gamma}_q \\
0 & 0 & 0  & 0 \\
0 & 0 & 0 & 0 \er . \]

\n Thus 

\[\dot{D}_q D_q^{-1} = \bl{cccc}
\dot{\alpha}_q A_{31} + \dot{\beta}_q A_{41} & \ast & \ast & \ast \\
\ast & \dot{\gamma}_q A_{42} & \ast & \ast \\
0 & 0 & 0 & 0 \\
0 & 0 & 0 & 0 \er\]

\n and

\spf

\n {\bf (A.8)} \ $Tr (\dot{D}_q D_q^{-1}) = Tr (\dot{\alpha}_q A_{31} + 
\dot{\beta}_q
A_{41}) + Tr (\dot{\gamma}_q A_{42}).$

\spf

\n {\bf Lemma A.2} (i) $A_{41} = 0;$

\spd

\n (ii) $(-A_{42} \ \underline{d}^{-1}_{1,q+1} \alpha_{q+1} + A_{43}) = 
\underline{d}_{2,q}^{-1} ;$

\spd

\n (iii) $A_{42} = - A_{43} (\underline{d}_{1,q+1}^{-1} \alpha_{q+1})^{\ast};$

\spd

\n (iv) $\underline{d}_{2,q}^{-1} = A_{43} ((\underline{d}_{1,q+1}^{-1} 
\alpha_{q+1})^{\ast} (\underline{d}^{-1}_{1,q+1} \alpha_{q+1}) + Id);$

\spd

\n (v) $A_{42} = - \underline{d}_{2,q}^{-1} (Id + (\underline{d}^{-1}_{1,q+1}
\alpha_{q+1})^{\ast} (\underline{d}^{-1}_{1,q+1} \alpha_{q+1}))^{-1} 
(\underline{d}^{-1}_{1, q+1} \alpha_{q+1})^{\ast};$

\spd

\n (vi) $\dot{\gamma}_q A_{42} = \underline{d}_{1, q+1}^{-1} 
\dot{\alpha}_{q+1} (Id + (\underline{d}^{-1}_{1,q+1} \alpha_{q+1} )^{\ast} 
(\underline{d}^{-1}_{1,q+1} \alpha_{q+1}))^{-1} (\underline{d}^{-1}_{1,q+1}
\alpha_{q+1})^{\ast}.$

\spf

\n {\bf Proof} (i) Take $x = (x_+, 0, 0, 0) \in {\cal C}_q.$ By Lemma $A.1, x \in 
{\cal C}^+_{q+1}$. Then

{\arraycolsep0.5mm
\[x = D_q D^{-1}_q  X = \bl {cccc}
0 & \underline{d}_{1,q} & \alpha_q & \beta_q \\
0 & 0 & 0 & \gamma_q \\
0 & 0 & 0 &  \underline{d}_{2,q} \\
0 & 0 & 0 & 0 \er \bl {cc}
A_{11} & x_+ \\
A_{21} & x_+ \\
A_{31} & x_+ \\
A_{41} & x_+ \er = \bl {ccc}
\ast &&  \\
\ast && \\
\underline{d}_{2,q} & A_{41} & x_+ \\
0 && \er \]}

\spd

\n and therefore $0 = \underline{d}_{2,q} A_{41} x_+.$ As $\underline{d}_{2,q}$
is an isomorphism we conclude that $A_{41} x_+ = 0.$

\spd

\n (ii) Take $x (0, 0, 0, y_-) \in {\cal C}_q.$ By Lemma A.1, $x \in {\cal C}^-_q.$
Then 

\[\begin{array}{l}
x = D^{-1}_q D_q x = D^{-1}_q (\beta_q y_-, \gamma_-, \underline{d}_{2,q}
y_{-}, 0) = \\
= (\ast, \ast, \ast, A_{41} \beta_q y_- + A_{42} \gamma_q y_- + A_{43} 
\underline{d}_{2,q} y_-).
\end{array}\]

\spd

\n In view of (i), $A_{41} = 0$ and thus 

\spf

\n {\bf (A.9)} \hsps $(A_{42} \ \gamma_q + A_{43} \ 
\underline{d}_{2,q}) \ y_- = y_-.$

\spd

\n Recall that $\gamma_q = - \underline{d}^{-1}_{1, q+1} \alpha_{q+1} 
\underline{d}_{2,q}$ and substitute into (A.9) to obtain (ii).

\spd

\n (iii) As $D^{-1}_q |_{{\cal C}^-_{q+1}} = 0$ and in view of Lemma A.1, 

\[\begin{array}{l}
0 = D^{-1}_q (0, x_{-}, (\underline{d}^{-1}_{1, q+1} \alpha_{q+1})^{\ast} 
x_{-}, 0) = \\
= (\ast, \ast, \ast, A_{42} \ x_- + A_{43} \ (\underline{d}^{-1}_{1,q+1}
\alpha_{q+1} )^{\ast} x_-).
\end{array}
\]

\n Thus $A_{42} + A_{43} (\underline{d}^{-1}_{1,q+1} \alpha_{q+1})^{\ast} = 0$
which proves (iii) 

\spd

\n (iv) follows from substituting (iii) into (ii).

\spd

\n (v) Notice that $Id + (\underline{d}^{-1}_{1,q+1} \alpha_{q+1})^{\ast}
(\underline{d}^{-1}_{1, q+1} \alpha_{q+1}) \ge Id$ is invertible and there\-fore,
we can solve (iv) for $A_{43}.$ Substituted into (iii) we obtain (v).

\spd

\n (vi) Substitute $\dot{\gamma}_q = - \underline{d}^{-1}_{1,q+1} 
\dot{\alpha}_{q+1} \underline{\dot{d}}_{2,q}$ into (v) to get (vi). \ \ \ \carre

\spf

\n {\bf Lemma A.3}

\spd

\n (i) $\alpha_q A_{31} + \underline{d}_{1,q} A_{21} = 
Id \ (\mbox{on} \ {\cal C}^+_{1,q}),$

\spd

\n (ii) $A_{31} = (\underline{d}^{-1}_{1,q} \alpha_q)^{\ast} A_{21};$

\spd

\n (iii) $A_{21} = (Id + (\underline{d}^{-1}_{1,q} \alpha_q) 
(\underline{d}_{1,q}^{-1} \alpha_q)^{\ast})^{-1} \underline{d}^{-1}_{1,q};$

\spd

\n (iv) $\dot{\alpha}_{q} A_{31} = \dot{\alpha}_q 
(\underline{d}^{-1}_{1,q} \alpha_q)^{\ast} (Id + (\underline{d}^{-1}_{1,q} 
\alpha_q) (\underline{d}_{1,q}^{-1} \alpha_q)^{\ast})^{-1} 
\underline{d}^{-1}_{1,q}.$

\spf

\n {\bf Proof} The proof of Lemma A.3 is similar to the one of Lemma A.2. 
\ \ \ \carre

\n {\bf Lemma A.4}

\spd

\n (i) $\sum_q (-1)^q Tr(\dot{D}_q D^{-1}_q) = 0;$

\spd

\n (ii) $\frac{d}{dt} \log \det T(t) = 0.$

\spd

\n {\bf Proof} (i) Substituting Lemma A.2 (i) and (vi) as well as Lemma A.3
(iv) into (A.8) leads to

{\arraycolsep0.1mm
\[
\begin{array}{lll}
Tr (\dot{D}_q D_q^{-1}) & =  & Tr (\dot{\alpha}_q A_{31}) + 0 + Tr 
(\dot{\gamma}_q A_{42}) \\
& = & Tr (\dot{\alpha}_q (\underline{d}^{-1}_{1,q} \alpha_q)^{\ast} 
(Id + (\underline{d}^{-1}_{1,q} \alpha_q)
(\underline{d}^{-1}_{1,q} \alpha_q)^{\ast})^{-1} \underline{d}^{-1}_{1,q})\\
& + & Tr (\underline{d}^{-1}_{1,q+1} \dot{\alpha}_{q+1} (Id + 
(\underline{d}^{-1}_{1,q+1} \alpha_{q+1})^{\ast} (\underline{d}^{-1}_{1,q+1}
\alpha_{q+1}))^{-1} 
(\underline{d}^{-1}_{1,q+1} \alpha_{q+1})^{\ast}) \\
& = & Tr ( \dot{\alpha}_q (\underline{d}^{-1}_{1,q} \alpha_q)^{\ast} (Id +
(\underline{d}^{-1}_{1,q} \alpha_q) (\underline{d}^{-1}_{1,q} \alpha_q)^{\ast}
)^{-1} \underline{d}^{-1}_{1,q}\\
& + & Tr (\dot{\alpha}_{q+1} (\underline{d}^{-1}_{1,q} \alpha_q)^{\ast} 
(Id + (\underline{d}^{-1}_{1, q+1} \alpha_{q+1})^{\ast} 
(\underline{d}^{-1}_{1,q+1} \alpha_{q+1}))^{-1} 
\underline{d}^{-1}_{1,q+1})
\end{array}\]}

\n where for the last equality we have used the fact that $Tr (AB) = Tr (BA)$
and $(Id + AA^{\ast})^{-1} A^{\ast} = A^{\ast} (Id + AA^{\ast})^{-1}.$
Thus with $\eta_q := \underline{d}^{-1}_{1,q} \alpha_q$

\[
\begin{array}{lll}
\sum_q (-1)^q Tr(\dot{D}_q D_q^{-1}) & = & \sum (-1)^q Tr 
(\dot{\alpha}_q
\eta^{\ast}_q (Id + \eta_q \eta_q^{\ast})^{-1} \underline{d}^{-1}_{1,q}) \\
& - & \sum (-1)^{q+1} Tr (\dot{\alpha}_{q+1} \eta^{\ast}_{q+1} 
(Id + \eta_{q+1} \eta^{\ast}_{q+1}) \underline{d}^{-1}_{1,q+1}) \\
& = & 0. \end {array} \]

\n (ii) follows from (i) and (A.6). \ \ \ \carre

\vspace{1cm}

\n {\large{\bf Appendix B \ Determinant class property}}

\spd

\n In this appendix, we discuss the concept of determinant class and 
provide examples of pairs $(M,\rho )$ which are not of 
determinant class.

Given an $\cal A$-Hilbert module 
${\cal {W}}$ and  $\varphi \in {\cal L}_{\cal A}(\cal {W})$ an operator 
which satisfies $Op 1- Op 6,$
one says that $\varphi $ is of determinant class  (cf [BFKM]) iff

\n {\bf (B.1)} \hsps $\int_{0+}^1 \log \lambda  dF_{|\varphi |}(\lambda) > - \infty $

\n where $|\varphi|= (\varphi^{\ast}\cdot \varphi)^{1/2}.$

It is not difficult to provide morphisms $\varphi $ which are  not 
of determinant class.
For example for $\cal A = 
\cal N(\bbz)$ and  ${\cal {W}}= l_2(\bbz)=L^2(S^1;\bbc)$, 
$S^1= \bbr/ \bbz,$ the multiplication  by 
a function $\alpha \in L^\infty (S^1;
\bbc)$ defines an element $\varphi= M_{\alpha} 
\in {\cal L}_{\cal A}(\cal {W})$ 
which satisfies $Op1- Op 6.$
Take $\alpha :[ 0, 1)\to \bbc$ defined by 
$\alpha (x)=\exp (-1/{x^2}).$ As 
$F_{|M_{\alpha}|}(\lambda)= \mu(\{x\in [0, 1], 
|\alpha(x)| \leq \lambda\})$, where $\mu (X)$ denotes the Lebesgue measure 
of the set $X,$  one can see that  
$F_{|M_{\alpha}|}(\lambda)= - (\log \lambda)^{1/2}$ for $0< \lambda< 1/e $ 
and the  integral (B.1)  is not convergent. 

Consider $(K,\tau, \rho, \mu)$ where  
$(K,\tau)$  is a  CW complex 
with finitely many cells in each dimension, $\rho$ is a representation of
the group $\pi=\pi_1(K)$ on the $\cal A$- Hilbert module of finite type $\cal {W},$  
and $\mu$ is a Hermitian structure in the 
flat bundle $\cal E\to K$ induced by 
$\rho.$ 

We say that  $(K,\tau, \rho, \mu)$ is of c-determinant 
class (cf [BFKM])
iff the associated cochain complex of $\cal A$-Hilbert modules of 
finite type $C^*(K, \tau,\rho,\mu)$
is of determinant class, i.e.  $\delta_i,$ or equivalently the 
combinatorial
Laplacians $\Delta^{comb}_i,$   are of determinant class for all $i$.

We also say that  $ (M,g,\rho, \mu),$ with $ (M,g)$ a smooth 
closed Riemannian manifold and  $(\rho, \mu)$  as above, is of a- determinant class
if the deRham complex $\Omega^*(M,\rho)$ of $\cal A$- Hilbert modules 
whose Hilbert
module structure is given by the scalar products induced by $g$ and $\mu,$
is of determinant class,  i.e.  $d_i,$ or equivalently the Laplacians $\Delta_i,$ are 
of determinant class.

It was shown in [BFKM]  (actually only for $\mu$  parallel,  but the same 
arguments remain valid in the generality presented here) that the c-determinant class 
property is independent of $\tau$ and $\mu,$
the  a-determinant class property is independent of  $g$ and $\mu$  and both
properties are
homotopy invariant.
Moreover for a compact manifold,  possibly with boundary, the  a-determinant 
class property
holds iff the  c-determinant class property holds. Therefore the determinant 
class property 
is a homotopy invariant for a pair $(K,\rho)$ with $K$ a compact  space
of the homotopy type of a CW complex, which can be defined 
both analytically and combinatorially.

If $\rho$ is a representation of $\Gamma,$ induced by a homomorphism 
$\pi: \Gamma\to G ,$  $G$ a countable discrete group, on the ${\cal N} (G)$- Hilbert 
module
${\cal {W}}=l_2(\Gamma),$
then $(K,\rho)$ is of determinant class when  $G$ is residually finite or 
$G$ is amenable.
Recently  B.Clair [C] has verified the determinant class property 
when $G$ is  a residually amenable group.

A representation $\rho\colon\pi_1(S^1)=
\bbz\to{\cal L}_{{\cal N}(\bbz)}(L^2(S^1;\bbc))$ is 
determined by $\rho (1)$.  Assume that $\rho (1)$ is 
the operator ${\cal M}_{1+
\alpha}$ given by multiplication by $1+\alpha$ where 
$\alpha\in L^\infty (S^1;
\bbc)$.  Consider the cochain complex induced by $\rho$ and the generalized 
triangulation $\tau =(h,g)$ with $h(t)={\cos (2\pi t)+1 \over 2}$ and $g$ the 
standard metric.  
\n As the triangulation $\tau$ has one $0$-cell $E_0$ and one 
$1$-cell $E_1$, the cochain complex is of the form
\[
0\to W\mathop{\to}\limits^\delta W\to 0
\]
where $\delta (E_0)=E_1-{\cal M}_{1+\alpha}E_1 =
-{\cal M}_\alpha E_1$. Therefore, 
\[\log T_{comb} = {1 \over 2}\log\det \delta^\ast\delta = 
{1 \over 2} \log\det
\alpha^\ast\alpha.\]

\n Observe  that $(M,\rho )$ is of determinant class iff 
${\cal M}_\alpha$ is 
of determinant class, thus  for $\alpha(x)= \exp (-1/{x^2}),$ 
as introduced above,  
$(S^1, \rho)$ is NOT of determinant class.  
 We point out that the regular 
representation of $\bbz$ corresponds to the function $\alpha(t)=
\exp (2\pi it)-1$.

\newpage

\n {\large{\bf Bibliography}}

\spd

\begin{description}

\item[{[B1]} D.Burghelea:] \hfill  ``Lectures on Witten-  Helffer- Sj\"ostrand theory'',  E.S.I.  preprint ,
Vienna, 1998.

\item[{[B2]} D.Burghelea:] \hfill  ``Removing Metric Anomalies from Ray- Singer Torsion'', Letters in Mathematical Physics 47 (1999), p.149-158.

\item[{[BZ]} J.P. Bismut, W. Zhang:] ``An extension of a theorem by
Cheeger and M\"uller'', Ast\'erisque 205 (1992), p. 1-223.

\item[{[BFK1]} D. Burghelea, L. Friedlander, T. Kappeler:] \hfill ``Witten 
deforma-
\linebreak
tion
of the analytic torsion and the Reidemeister torsion'', AMS Transl. (2), vol 184 (1998),
p. 23-34

\item[{[BFK2]} D. Burghelea, L. Friedlander, T. Kappeler:] ``Torsions for
\linebreak manifolds with boundary and gluing formulas'', IHES preprint, 1996,
(to appear in Math. Nach.).

\item[{[BFK3]} D. Burghelea, L. Friedlander, T. Kappeler:] \hfill``Asymptotic ex- 
\linebreak
pan\-sion of the Witten deformation of the analytic torsion'', J. Funct. Anal.
137 (1996), p. 320-363.

\item[{[BFK4]} D. Burghelea, L. Friedlander, T. Kappeler:] ``Mayer-Vietoris
\linebreak type formula for determinants of elliptic differential operators'', J. Funct.
Anal. 107 (1992), p. 34-66.

\item[{[BFK5]} D. Burghelea, L. Friedlander, T. Kappeler:] ``Relative torsion for homotopy  triangulation''
,to appear in Contemporary Mathematics, AMS 1999.

\item[{[BFKM]} D. Burghelea, L. Friedlander, 
T. Kappeler,  P. McDonald: ]
 ``Analytic and Reidemeister torsion for representations in finite 
type 
\newline Hilbert
modules'', GAFA 6 (1996), p. 751-859.

\item[{[C]} B. Clair:]  
"Residual amenability and the approximation 
 of $L^2$ 
 \newline -invariants,"
 preprint Univ. Chicago ,1997.

\item[{[CFM]} \hfill A.L. Carey, M.Farber, \hfill V. Mathai:] 
"Determinant lines, von Neumann algebras and $L^2$ torsion,"
J.reine angew. Math.484 (1997), p. 153-181.

\item[{[CMM]} \hfill A.L. Carey, V. Mathai, \hfill A. Mishchenko:] 
``On \hfill analytic \hfill torsion  
\linebreak
over
$C^{\ast}$-algebras'', preprint.

\item[{[Di]} J. Dixmier:] $C^{\ast}$-algebras, North-Holland, Math. Library, 
Vol 15, 1997.

\item[{[HS]} B. Helffer, J. Sj\"ostrand:] ``Puits multiples en m\'ecanique
semi-classique IV - Etude du complexe de Witten'', Comm. in PDE 10 (1985), 
p. 245-340.

\item[{[Gi]} P. Gilkey:] ``Invariance theory, the heat equation and the
Atiyah-Singer index theorem'', Publish or Perish, Wilmington,  1984.

\item[{[GS]} M. Gromov, M.A. Shubin:] ``Von Neumann \hfill 
spectrum \hfill near zero'', 
\linebreak
GAFA 1 (1991), p. 375-404.

\item[{[L]} J.Lott:] ``Heat kernels on covering spaces and topological
\linebreak
invariants, J. of Diff. Geo. 35 (1992), p 471-510.

\item[{[RS]} M. Reed, B. Simon:] ``Methods in modern mathematical physics'',
Vol 1, $2^{\mbox{nd}}$ ed., Academic Press, New York, 1980.

\item[{[Tu]} Turaev,V] :  "Euler structures, nonsingular vector fields and torsion of Reidemeister type",
Izv.Acad. Sci. USSR. 53   (1989), p 130-146.  

\item[{[Wi]} E. Witten:] ``Super symmetry and Morse theory'', J. of Diff.
Geom. 17 (1982), p. 661-692.
\end{description}

\end{document}